\PassOptionsToPackage{unicode}{hyperref}
\PassOptionsToPackage{hyphens}{url}
\PassOptionsToPackage{dvipsnames,svgnames,x11names,table}{xcolor}
\documentclass[
  paper=a4,
  ,captions=tableheading
]{scrartcl}
\usepackage{amsmath,amssymb}
\usepackage{setspace}
\usepackage{iftex}
\ifPDFTeX
  \usepackage[T1]{fontenc}
  \usepackage[utf8]{inputenc}
  \usepackage{textcomp} % provide euro and other symbols
\else % if luatex or xetex
  \usepackage{unicode-math} % this also loads fontspec
  \defaultfontfeatures{Scale=MatchLowercase}
  \defaultfontfeatures[\rmfamily]{Ligatures=TeX,Scale=1}
\fi
\usepackage{lmodern}
\ifPDFTeX\else
\fi
\IfFileExists{upquote.sty}{\usepackage{upquote}}{}
\IfFileExists{microtype.sty}{% use microtype if available
  \usepackage[]{microtype}
  \UseMicrotypeSet[protrusion]{basicmath} % disable protrusion for tt fonts
}{}
\makeatletter
\@ifundefined{KOMAClassName}{% if non-KOMA class
  \IfFileExists{parskip.sty}{%
    \usepackage{parskip}
  }{% else
    \setlength{\parindent}{0pt}
    \setlength{\parskip}{6pt plus 2pt minus 1pt}}
}{% if KOMA class
  \KOMAoptions{parskip=half}}
\makeatother
\usepackage{xcolor}
\definecolor{default-linkcolor}{HTML}{A50000}
\definecolor{default-filecolor}{HTML}{A50000}
\definecolor{default-citecolor}{HTML}{4077C0}
\definecolor{default-urlcolor}{HTML}{4077C0}
\usepackage[margin=2.5cm,includehead=true,includefoot=true,centering,]{geometry}
\usepackage{footnotebackref}
\providecommand{\tightlist}{%
  \setlength{\itemsep}{0pt}\setlength{\parskip}{0pt}}
\ifLuaTeX
\usepackage[bidi=basic]{babel}
\else
\usepackage[bidi=default]{babel}
\fi
\babelprovide[main,import]{american}
\def\languageshorthands#1{}
\usepackage{amsfonts, amsthm, amsmath, amssymb, bbm}
\usepackage{mathtools, graphicx, hyperref}
\usepackage{biblatex}
\usepackage{multicol, float, caption, faktor, caption}

\usepackage{mathabx}
\usepackage{stmaryrd}
\usepackage{tikz-cd, float, suffix, graphicx, ebproof, logicproof, cancel, mathdots, diagbox}
\usepackage[T1]{fontenc}
\usepackage[multiple]{footmisc}
\usepackage{etoolbox}
\usepackage{pdfpages}
\usetikzlibrary{calc}
\usetikzlibrary{shapes}
\usetikzlibrary{shapes.multipart}
\nocite{*}

\newcommand{\ostar}{ 
  {\mathbin{
    \mathchoice
      {\buildcirclestar{\displaystyle}}
      {\buildcirclestar{\textstyle}}
      {\buildcirclestar{\scriptstyle}}
      {\buildcirclestar{\scriptscriptstyle}}
  }} 
}
\newcommand\buildcirclestar[1]{%
  \begin{tikzpicture}[baseline=(X.base), inner sep=0, outer sep=0]
    \node[draw,circle] (X)  {$#1\star$};
  \end{tikzpicture}%
}

\newcounter{countT}
\newtheorem{The}[countT]{Theorem}
\newcounter{countT2}
\newtheorem{The2}[countT2]{Theorem}
\newcounter{countP}
\newtheorem{Pro}[countP]{Proposition}
\newcounter{countL}
\newtheorem{Lem}[countL]{Lemma}
\newcounter{countC}
\newtheorem{Cor}[countC]{Corollary}
\newcounter{countD}
\newtheorem{Def}[countD]{Definition}
\newcounter{countE}
\newtheorem{Ex}[countE]{Example}
\newtheorem*{Ex*}{Ex}
\newcommand\bgC[1][]{\begin{Cor}[#1]}
\newcommand\ndC{\end{Cor}}
\newcommand\bgD[1][]{\begin{Def}[#1]}
\newcommand\ndD{\end{Def}}
\newcommand\bgL[1][]{\begin{Lem}[#1]}
\newcommand\ndL{\end{Lem}}
\newcommand\bgP[1][]{\begin{Pro}[#1]}
\newcommand\ndP{\end{Pro}}
\newcommand\bgT[1][]{\begin{The}[#1]}
\newcommand\ndT{\end{The}}
\newcommand\bgE[1][]{\begin{Ex}[#1]}
\newcommand\ndE{\end{Ex}}
\WithSuffix\newcommand\bgE*[1][]{\begin{Ex*}[#1]}
\WithSuffix\newcommand\ndE*{\end{Ex*}}
\newcommand\bgpr[1][]{\begin{proof}[Proof]}
\newcommand\ndpr{\end{proof}}

\let\mbb\mathbbm

\let\mca\mathcal
\newcommand\mtt[1]{\text{\fontfamily{qtm}\selectfont{#1}}}
\newcommand\bgmat{\begin{matrix}}
\newcommand\ndmat{\end{matrix}}
\newcommand\bgpmat{\begin{pmatrix}}
\newcommand\ndpmat{\end{pmatrix}}
\newcommand\sepList[2]{\foreach \n [count=\ni] in {#1} {
  \ifnum\ni=1{\n}
  \else{#2\n}
  \fi}
}

\newcommand{\nat}{{\mbb{N}}}
\newcommand{\ent}{{\mbb{Z}}}

\newcommand{\id}{\mtt{id}}

\newcommand{\bgdiagram}{
  \begin{figure}[H]
  \centering
  \begin{tikzcd}[column sep = small]}
\newcommand{\nddiagram}{
  \end{tikzcd}
  \end{figure}}

\mathchardef\mathcomma=\mathcode`+

\let\logicsbox\Box
\renewcommand\Box{\mathop{\logicsbox}}
\newcommand\rres{\mathbin{\backslash}}
\newcommand\lres{\mathbin{/}}

\renewenvironment{abstract}[1]{%
\begin{center}\normalfont\textbf{Abstract}\end{center}
\begin{quotation} #1 \end{quotation}
}{%
\vspace{1cm}
}

\patchcmd{\logicproof}{\center}{\flushleft}{}{}
\patchcmd{\endlogicproof}{\endcenter}{\endflushleft}{}{}
\ifLuaTeX
  \usepackage{selnolig}  % disable illegal ligatures
\fi
\IfFileExists{bookmark.sty}{\usepackage{bookmark}}{\usepackage{hyperref}}
\IfFileExists{xurl.sty}{\usepackage{xurl}}{} % add URL line breaks if available
\hypersetup{
  pdftitle={A Study on Actions for Atomic Logics},
  pdfauthor={Raül Espejo Boix},
  pdflang={en-US},
  hidelinks,
  breaklinks=true,
  pdfcreator={LaTeX via pandoc with the Eisvogel template}}
\title{A Study on Actions for Atomic Logics}
\author{Raül Espejo Boix}
\date{}
\usepackage{pagecolor}
\usepackage{afterpage}
\usepackage[margin=2.5cm,includehead=true,includefoot=true,centering]{geometry}

\PassOptionsToPackage{hyphens}{url}

\usepackage{csquotes}

\definecolor{caption-color}{HTML}{777777}
\usepackage[font={stretch=1.2}, textfont={color=caption-color}, position=top, skip=4mm, labelfont=bf, singlelinecheck=false, justification=raggedright]{caption}
\setcapindent{0em}

\definecolor{blockquote-border}{RGB}{221,221,221}
\definecolor{blockquote-text}{RGB}{119,119,119}
\usepackage{mdframed}
\newmdenv[rightline=false,bottomline=false,topline=false,linewidth=3pt,linecolor=blockquote-border,skipabove=\parskip]{customblockquote}
\renewenvironment{quote}{\begin{customblockquote}\list{}{\rightmargin=0em\leftmargin=0em}%
\item\relax\color{blockquote-text}\ignorespaces}{\unskip\unskip\endlist\end{customblockquote}}

\ifnum 0\ifxetex 1\fi\ifluatex 1\fi=0 % if pdftex
    \usepackage[default]{sourcesanspro}
  \usepackage{sourcecodepro}
  \else % if not pdftex
    \usepackage[default]{sourcesanspro}
  \usepackage{sourcecodepro}

  \ifxetex
    \makeatletter
    \defaultfontfeatures[\ttfamily]
      { Numbers   = \sourcecodepro@figurestyle,
        Scale     = \SourceCodePro@scale,
        Extension = .otf }
    \makeatother
  \fi
  \fi

\definecolor{heading-color}{RGB}{40,40,40}
\addtokomafont{section}{\color{heading-color}}
\usepackage{titling}
\title{A Study on Actions for Atomic Logics}
\author{Raül Espejo Boix}
\date{}

\usepackage[headsepline,footsepline]{scrlayer-scrpage}

\newpairofpagestyles{eisvogel-header-footer}{
  \clearpairofpagestyles
  \ihead*{A Study on Actions for Atomic Logics}
  \chead*{}
  \ohead*{}
  \ifoot*{Raül Espejo Boix}
  \cfoot*{}
  \ofoot*{\thepage}
  \addtokomafont{pageheadfoot}{\upshape}
}
\begin{document}

%%
%% begin titlepage
%%
\begin{titlepage}
\newgeometry{left=6cm}
\newcommand{\colorRule}[3][black]{\textcolor[HTML]{#1}{\rule{#2}{#3}}}
\begin{flushleft}
\noindent
\\[-1em]
\color[HTML]{5F5F5F}
\makebox[0pt][l]{\colorRule[435488]{1.3\textwidth}{4pt}}
\par
\noindent

{
  \setstretch{1.4}
  \vfill
  \noindent {\huge \textbf{\textsf{A Study on Actions for Atomic
Logics}}}

  {A Thesis for the Master’s Degree in Fundamental Mathematics of the University
of Rennes}
    \vskip 3em
  \noindent {\Large \textsf{Raül Espejo Boix}}
  \vskip 2em
  {\textit{Supervisors}}
  
  {\Large \textsf{Dr. Guillaume Aucher (Université de Rennes)}\\
  \textsf{Prof. Rajeev P. Goré (Australian National University)}}
  \vfill
}

\textsf{}
\end{flushleft}
\end{titlepage}
\restoregeometry
\pagenumbering{arabic} 

%%
%% end titlepage
%%

% \maketitle

\renewcommand\bgpr[1][]{\begin{proof}}
\floatstyle{boxed}
\newfloat{rules}{H}{lop}
\floatname{rules}{Rules}

\clearpage
\null
\thispagestyle{empty}
\clearpage

\maketitle

\thispagestyle{empty}

\abstract{
    Nowadays there is a large number of non-classical logics, each one best suited for reasoning about some issues in abstract fields, such as linguistics or epistemology, among others.
    Proving interesting properties for each one of them supposes a big workload for logicians and computer scientists.
    We want an approach into this problematic that is modular.
    To adress this issue, the report shows new insights in the construction of Atomic Logics introduced by Guillaume Aucher. 
    Atomic Logics let us represent very general left and right introduction rules and they come along a new kind of rules based on display logics and residuation.
    A new approach is taken into the definition of Atomic Logics, which is now built on a class of actions for which we prove cut-elimination.
    We show that some of them are equivalent to Aucher's Atomic Logics and we prove cut-elimination and Craig Interpolation for a class of them.
    The introduced theory is applied to the non-associative Lambek Calculus throughout the report.
    It is accompanied by a computer-checked formalisation of the original syntax in the proof assistant Coq.
}

\vfill

\hypertarget{acknowledgements}{%
\subsection*{Acknowledgements}\label{acknowledgements}}

I would like to thank my supervisors Rajeev Goré and Guillaume Aucher
for the support and guidance I received during this master thesis, thank
them for the opportunity this stage has been and also thank anyone who
has taken part by giving advice or counsel throughout the stage.

\vfill
\newpage
\tableofcontents
\thispagestyle{empty}
\newpage

\hypertarget{introduction}{%
\section*{Introduction}\label{introduction}}

\setcounter{page}1

Logics originated as the area of philosophy which studies the
\emph{forms of valid reasoning}. What we now consider Classical Logics
arises from Aristotelian syllogisms, in particular it is the historical
antecedent of predicate logics, designed during the end of the 19th
century to study mathematical arguments. The patterns of reasoning which
can be expressed with Classical Logics can not properly study most of
our daily life experiences nor our knowledge or uncertainties. For this
reason the first non-classical logics appeared, studying some phenomena
which falls outside of the propositional logics and for which predicate
logic was too expressive.

For this Master Thesis I have been presented a class of Logics, Atomic
Logics, which aim to provide an uniform and generic way to explore and
study non-classical logics. The overall objective was to develop a Coq
library for it, with the possible goal of expanding it into a Coq
library for Non-Classical Logics, so that logicians can prove general
properties for them in a common setting.

The Master Thesis research has worked with the \emph{mathcomp} library
for formalizing Atomic Logics in Coq. On
\href{https://github.com/The-busy-man/Universal-FO-Logics/commit/fe6a68645030adee9f51bae8d10e212f9300004b}{\color{blue} the \textit{proof of action} commit}
we can find its calculus without boolean negation. The difficulties
found in the process have redirected the objective into the addition of
some changes into the theory. We have had to work with group theory to
provide a finite action for Atomic Logics, we have had to introduce new
definitions for formally providing a concept of Connective Family and
Structural Connective and we have managed to provide display calculi for
a much bigger class of actions than the original one. This report is the
result of my intership at IRISA and it is structured as follows.

This \hyperref[sect1]{\textbf{first section}} presents some of the
logical concepts that are necessary to understand what Atomic Logics
are. This includes a presentation of a case study, using Lambek Calculus
from Linguistics, an informal presentation of logics and a formal
presentation of proof-theoretical concepts. In the
\hyperref[sect2]{\textbf{second section}} we briefly present the
mathematical concepts we are going to work with, mainly those of groups
and actions. In the \hyperref[sect3]{\textbf{third section}} we give out
of context Guillaume Aucher's action for Atomic Logics \(\alpha*\beta\)
and we show two ways in which we can redefine it. In the
\hyperref[sect4]{\textbf{fourth section}} we give the lacking context by
formally defining the algebraic objects that represent Connectives and
Structural Connectives. In the \hyperref[sect5]{\textbf{fifth section}}
this is used to redefine Atomic Logics on general Connective Families,
no matter what action they are using. A condition for the system to be
properly displayable is given in terms of its action (we can say that it
must be transversal \emph{enough}). We also prove Craig Interpolation
for Connective Families. This formalisation only accounts for the Atomic
Logics and further research is required for extending it into Aucher's
Molecular Logics and the addition of structural rules. In the
\hyperref[sect6]{\textbf{sixth section}} we get back Aucher's Atomic
Logics and we prove some interesting results about them.

During the report we will present small snippets of the
\href{https://github.com/The-busy-man/Universal-FO-Logics/commit/fe6a68645030adee9f51bae8d10e212f9300004b}{\color{blue}\textbf{Coq code}},
in the proof of action commit.\\
The present version of the code can be found at
\href{https://github.com/The-busy-man/Universal-FO-Logics/tree/test}{\color{blue} https://github.com/The-busy-man/Universal-FO-Logics/tree/test},
it is in an unfinished transition into using connective families.

\hypertarget{logical-notions}{%
\section{Logical Notions}\label{logical-notions}}

\label{sect1}

\hypertarget{an-example-from-linguistics}{%
\subsection{An Example from
Linguistics}\label{an-example-from-linguistics}}

Traditionally linguistics presents the following elementary units, which
constitute the atoms of a language's grammar.

\begin{itemize}
\tightlist
\item
  \(S\): Sentence.
\item
  \(N\): Noun.
\item
  \(NP\): Noun Phrase.
\item
  \(AP\): Adjective Phrase.
\item
  \(VP\): Verb Phrase.
\item
  \(PP\): Adpositional Phrase.
\end{itemize}

Using this elementary elements we can build new expressions just by
concatenation. Then, the product connective \(\otimes\) is the
concatenation operator. This operator can be given left and right
divisions: for any two categories \footnote{Note that we are using
  categories in a loose sense, following the structuralist language of
  Moortgat's book.} \(X\) and \(Y\) we say that the functor \(X/Y\)
combines with an argument \(Y\) to the right to form an \(X\) and the
functor \(Y\backslash X\) combines with an argument \(Y\) to the left to
form an \(X\). \cite{moortgat_1992}

\begin{quote}
We can see that \(\otimes\) works as a sort of conjunction by
representing plain and simple concatenation of concepts and
\(/,\,\backslash\), called its left and right residuals, work as
functors between the continguous categories.
\end{quote}

Syntactic trees can now be represented as a proposition in the language
formed by the connectives \(\{\otimes,\,/,\,\backslash\}\), which we
call \(\mbb{FL}\) on propositional letters
\(\{S,\,N,\,NP,\,AP,\,VP,\,PP\}\).

\bgE

\emph{(Jo) presento la lògica de la gramàtica}, where:

\begin{itemize}
\tightlist
\item
  \emph{Jo} is an \(NP\).\\
  This category can be understood as a personal pronoun.
\item
  \emph{Presento} is an \((NP\backslash S)/NP\).\\
  This category can be understood as a transitive verb.
\item
  \emph{La} is an \(NP/N\).\\
  This can be understood as a determiner.
\item
  \emph{Lògica} is an \(N\).\\
  This can be understood as a common noun.
\item
  \emph{De} is an \((N\backslash N)/NP\).\\
  This can be understood as a preposition.
\item
  \emph{Gramàtica} is an \(N\).\\
  This can be understood as a common noun.
\end{itemize}

To represent what in this context we have written as \(x\) is an \(X\)
we will write \(x:X\).

\begin{itemize}
\tightlist
\item
  This example would then correspond in \(\mbb{FL}\) to:
  \begin{multline*}
  \textit{(Jo) presento la lògica de la gramàtica}:\\
  NP
  \otimes
  (((NP\backslash S)
    /
    NP)
   \otimes
   ((NP/N)
    \otimes
    (N
     \otimes
     (((N\backslash N)
       /
       NP)
      \otimes
      ((NP/N)
       \otimes
       N)))))
  \end{multline*}
\item
  By the rules informally presented previously this reduces to
  \[\textit{(Jo) presento la lògica de la gramàtica}:S\] This can be
  interpreted as the sentence being well-formed.
\end{itemize}

\ndE

The particular formal rules of the logics we just presented can be
encaptured in the following proof system, which we will call
Non-Associative Lambek Calculus: \cite{lambek61}

\begin{rules}
\begin{multicols}{3}
    \underline{\textbf{Axiom}}:

    $$\begin{prooftree}
        \infer0[Id]{
            A\vdash A}
    \end{prooftree}$$

    \columnbreak
    \underline{\textbf{For $\rres$}}:
    $$\begin{prooftree}
        \hypo{
            B\vdash A\rres C}
        \infer1[\(L1\)]{
            A\otimes B\vdash C
        }
    \end{prooftree}$$
    $$\begin{prooftree}
        \hypo{
            A\otimes B\vdash C}
        \infer1[\(L2\)]{
            B\vdash A\rres C
        }
    \end{prooftree}$$

    \columnbreak
    \underline{\textbf{For $\lres$}}:  
    $$\begin{prooftree}
        \hypo{
            A\vdash C\lres B}
        \infer1[\(L3\)]{
            A\otimes B\vdash C
        }
    \end{prooftree}$$
    $$\begin{prooftree}
        \hypo{
            A\otimes B\vdash C}
        \infer1[\(L4\)]{
            A\vdash C\lres B
        }
    \end{prooftree}$$
\end{multicols}
\caption*{Rules of $\mathbb{FL}$}
\end{rules}

\begin{figure}\begin{center}

\includegraphics{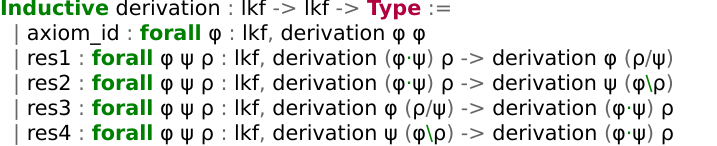}

\caption{Implementation of Non-Associative Lambek Calculus}
\label{non_assoc_lamb}

\end{center}\end{figure}

We have written it in Coq in Figure \ref{non_assoc_lamb}. Throughout the
report we will be showing some Coq implementations from the repository
related to what is being discused.

\hypertarget{logics}{%
\subsection{Logics}\label{logics}}

The development of the field gave birth to a common methodology, being
so that a Logic is now understood as an object of study itself, with
multiple possible approaches into its construction.

Commonly, a logical system is constituted by three parts: syntax,
semantics and a proof system.

\begin{itemize}
\tightlist
\item
  The syntax expresses how the concepts are written, usually through
  objects called formulas given by a grammar.
\item
  The semantics gives meaning to syntax, usually through objects called
  models.
\item
  The proof systems gives appropriate and desired ways of reasoning with
  the syntax (not the semantics), usually through objects called rules
  arranged in sets called calculus. In this case the calculi are usually
  divided between Hilbert/Sequent systems, depending on how much does it
  focuses on the axioms.
\end{itemize}

As commonly understood, both semantics and proof systems alone are
enough to characterize a Logic. For example, we can begin by giving a
class of models and then deriving the validities from them. It happens
frequently that they are defined simultaneously, then it is of our
interest showing that both methods are talking about the same Logic. It
is not common to discover a Logic's syntax and semantics at the same
time. For example, modal logics' Kripke semantics were not found until
1959, while the proof system was developed in the forties and sixties.

\begin{quote}
We will fix the meaning of logics on its semantical approach.
\end{quote}

We will say that a \textit{\textbf{logics}} is a triple
\((\mca{L},\,\mca{E},\,\vDash)\) where

\begin{itemize}
\tightlist
\item
  \(\mca{L}\) is a logical language built as a set of well-formed
  expressions from a set of connectives and a set of propositional
  letters.
\item
  \(\mca{E}\) is a class of pointed models.
\item
  \(\vDash\) is a satisfaction relation relating in a compositional
  manner elements of \(\mca{L}\) to models of \(\mca{E}\).
\end{itemize}

The set of formulas in \(\mca{L}\) such that they are satisfied in all
pointed models is called the \textit{\textbf{set of validities}}.

\bgE

We define semantics for Lambek Calculus.

\label{routley-meyer}

We define the \emph{Routley-Meyer models} for Lambek Calculus as
structures \((W,\,R,\,V)\), where \(W\) is a non-empty set of worlds,
\(R\) is a ternary relation over \(W\) and \(V\) is a valuation from
\(\mbb{A}\times W\) to the booleans, where \(\mbb{A}\) is the set of
propositional letters in our language.

Then we define inductively on formulas its
\textit{\textbf{interpretation function}} through:

\begin{itemize}
\tightlist
\item
  \(\llbracket p\rrbracket := \{w\in V(p)\}\).
\item
  \(\llbracket \varphi\otimes\psi\rrbracket := \{w\in W\mid\exists u,\,v\in W,\,u\in\llbracket\varphi\rrbracket\mtt{ and }v\in\llbracket\psi\rrbracket\mtt{ and }R\ u\ v\ w\}\).
\item
  \(\llbracket \varphi\rres\psi\rrbracket := \{w\in W\mid\forall u,\,v\in W,\,R\ u\ w\ v\mtt{ and }u\in\llbracket\varphi\rrbracket\mtt{ imply }v\in\llbracket\psi\rrbracket\}\).
\item
  \(\llbracket \varphi\lres\psi\rrbracket := \{w\in W\mid\forall u,\,v\in W,\,R\ w\ v\ u\mtt{ and }v\in\llbracket\psi\rrbracket\mtt{ imply }u\in\llbracket\varphi\rrbracket\}\).\cite{dunn_bimbo}
\end{itemize}

We define \((M,\,w)\vDash\varphi\) to hold if and only if
\(w\in\llbracket M\rrbracket\).

\begin{figure}\begin{center}

\includegraphics{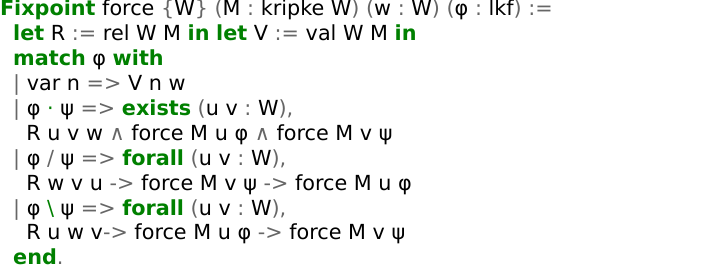}

\caption{Implementation of Routley-Meyer semantics for Lambek Calculus with Internalised Forcing}

\label{forcing}

\end{center}\end{figure}

\ndE

\hypertarget{sequent-calculi}{%
\subsubsection*{Sequent Calculi}\label{sequent-calculi}}

In the early 20th century a replacement for Hilbert systems was sought
in logics. They needed a system that adapted better to the usual proof
methodology of mathematicians, while working in the full formality of
proof systems. One of the logicians who developed the new methods was
Gerhard Gentzen, who found the first natural deduction system.

In his 1969 article \cite{gentzen_cut} Gentzen expanded his rule-set for
classical and intuitionistic logics by introducing sequent deduction
systems. It required a metatheoretical implication (which we note
\(\vdash\)) bonding two formulas into a new object called sequent, but
it made possible to work proofs without opening hypothesis (in contrast
with natural deduction).

\begin{quote}
Compare a possible Gentzen's sequent system and a Hilbert system proofs of the law of the excluded middle:

\textit{\textbf{Gentzen}},

$$\begin{prooftree}
    \infer0[Id]{
        A\vdash A
    }
    \infer1[\(\vee\)-IS]{
        A\vdash A\vee\neg A
    }
    \infer1[\(\neg\)-IS]{
        \vdash A\vee\neg A,\,\neg A
    }
    \infer1[\(\vee\)-IS]{
        \vdash A\vee\neg A,\, A\vee\neg A
    }
    \infer1[Contraction]{
        \vdash A\vee\neg A
    }
  \end{prooftree}$$
  
\textit{\textbf{Hilbert}},

    {\footnotesize\begin{logicproof}{1}
    A\to A\vee(A\to\bot) & Axiom \\
    ((A\vee(A\to\bot))\to\bot)\to A\to\bot & N 1 \\
    (A\to\bot)\to A\vee(A\to\bot) & Axiom \\
    ((A\vee(A\to\bot))\to\bot)\to (A\to\bot)\to\bot & N 3 \\
    ((A\to\bot)\to\bot)\to A & Axiom \\
    ((A\vee(A\to\bot))\to\bot)\to A & B 4 5 \\
    (((A\vee(A\to\bot))\to\bot)\to A\to\bot)\to(((A\vee(A\to\bot))\to\bot)\to A)\to((A\vee (A\to\bot))\to\bot)\to\bot & Axiom \\
    (((A\vee(A\to\bot))\to\bot)\to A)\to((A\vee (A\to\bot))\to\bot)\to\bot & MP 7 2 \\
    ((A\vee (A\to\bot))\to\bot)\to\bot & MP 8 6 \\
    (((A\vee (A\to\bot))\to\bot)\to\bot)\to A\vee(A\to\bot) & Axiom \\
    A\vee(A\to\bot) & MP 10 9
  \end{logicproof}}

Where we call $B$ the lemma

    {\footnotesize\begin{logicproof}{1}
    A \to B & Premise \\
    B\to C & Premise \\
    (B\to C)\to A\to B\to C & Axiom \\
    A\to B\to C & MP 2 3 \\
    (A\to B\to C)\to (B\to C)\to A\to C & Axiom \\
    (B\to C)\to A\to C & MP 5 4\\
    A\to C & MP 6 1
    \end{logicproof}}

And we call $N$ the lemma

{\footnotesize\begin{logicproof}{1}
    A \to B & Premise \\
    (A\to B\to\bot)\to(A\to B)\to A\to\bot & Axiom \\
    (B\to\bot)\to A\to B\to\bot & Axiom \\
    (B\to\bot)\to(A\to B)\to A\to\bot & B 3 2 \\
    ((B\to\bot)\to(A\to B)\to A\to\bot)\to((B\to\bot)\to(A\to B))\to(B\to\bot)\to A\to\bot & Axiom \\
    ((B\to\bot)\to(A\to B))\to(B\to\bot)\to A\to\bot & MP 5 4 \\
    (A\to B)\to (B\to\bot)\to(A\to B) & Axiom \\
    (B\to\bot)\to(A\to B) & MP 7 1 \\
    (B\to\bot)\to A\to\bot & MP 6 8
    \end{logicproof}}

Both proof styles always produce tree structures.
\end{quote}

We begin with a language of formulas, given inductively by a set of
variables and connectives (including constants). Each formula is
completely characterised by a labeled tree structure, called
\textit{\textbf{formation tree}}. Indeed, we can consider atoms as
leaves and composition of formulas through a connective as a common
parent of the nodes corresponding to each of its components, we then
label each node with the full formula used to construct it.

We define subformulas of \(\varphi\) as the labels appearing in
\(\varphi\)'s tree. We define a formula occurrence \(\psi\) in
\(\varphi\) as a subformula \(\psi\) in \(\varphi\) along with a
particular node labeled \(\psi\) in the tree structure of \(\varphi\).

For each of the elementary connectives we define its structural
connective as another connective with the same arity. For example, in
Belnap's article\cite{belnap_display} the structure of \(\wedge\) is
noted \(\circ\), another common notation is the semicolon \(;\) but in
our case we will work with Gentzen's comma \(,\,\)\cite{gentzen_cut}. We
define structures as the closure of formulas under application of
structural connectives. Substructures and structure occurrences are
defined similarly to the formula case but now labels are structural
variables, structural connectives and full formulas. For each pair of
structures \(X\), \(Y\) we may form a sequent \(X\vdash Y\).

Let \(\mca{L}\) be a formal language. From any set of formulas and an
extra formula in \(\mca{L}\) we may form an \textit{\textbf{inference}}.
We call the former \textit{\textbf{the premises}} and the latter
\textit{\textbf{the conclusion}}.

A \textit{\textbf{formal system}} or \textit{\textbf{proof system}}
based on a language \(\mca{L}\) is defined by a set of rules, called
\textit{\textbf{inference rules}}. Each rule is a set of inferences with
a fixed number of premises.

We call the inference rules with no premises axioms.

\begin{quote}
The previous definition's reference to formulas has to be understood in
the broader context of a formal language. In this sense we can also
consider sequents as formulas in the formal language of sequents.
\end{quote}

\bgE

In the introductory section we worked with the formal language given by
\[\mca{L}_L ::= V\ |\ \mca{L}_L\otimes \mca{L}_L\ |\ \mca{L}_L\lres \mca{L}_L\ |\ \mca{L}_L\rres \mca{L}_L\]
Where \(V = \{S,\,N,\,NP,\,AP,\,VP,\,PP\}\).

This notation must be understood as follows:

\begin{itemize}
\tightlist
\item
  Each \(p\in V\) is a well-formed formula.
\item
  For any \(p,\,q\) well-formed formulas, \(p\otimes q\) is also a
  well-formed formula.
\item
  For any \(p,\,q\) well-formed formulas, \(p\lres q\) is also a
  well-formed formula.
\item
  For any \(p,\,q\) well-formed formulas, \(p\rres q\) is also a
  well-formed formula.
\end{itemize}

Another common formal language we will work with is the one for
propositional logics:
\[\mca{L}_B ::= V\ |\ \bot\ |\ \top\ |\ \mca{L}_B\wedge \mca{L}_B\ |\ \mca{L}_B\vee \mca{L}_B\ |\ \mca{L}_B\rightarrow \mca{L}_B\ |\ \neg \mca{L}_B\]

\ndE

\begin{quote}
Rules for a proof system are to be understood as transformations on
derivable formulas. Sequents are to be understood as meta implications
between formulas, implications with an extra abstraction level. The
meaning of sequents is common in logical systems we are working with.

For the translation from structures to formulas there is some more
detail which we will comment further, as commonly structures are to be
understood as syntactical clauses such as parenthesis, shaping the way
we can work with sequents (so that they might be seen as sets,
sequences, etc.). In spite of that, we will usually see structures as
meta logical connectives. The meaning of structures might vary depending
on the calculi rule-set (these rules will be accordingly named
structural rules).
\end{quote}

We will call \textit{\textbf{Hilbert systems}} the proof systems based
on formulas.\\
We will call \textit{\textbf{Gentzen systems}} or sequent systems the
proof systems based on sequents.

\begin{quote}
Although, as mentioned before, the difference between Hilbert and
Gentzen systems is commonly considered as given by the proportion of
axioms over rules we will follow the Francesca Poggiolesi's criteria in
\emph{Gentzen Calculi for Modal Propositional Logic} for this formal
approach.\cite{poggi_gentzen}
\end{quote}

Let \(\mca{L}\) be a formal language. A \textit{\textbf{derivation}} in
a proof system \(\mca{P}\) based on \(\mca{L}\) is a finite, upward
growing tree. The nodes of the tree are labeled with formulas of
\(\mca{L}\). For each intermediate node the label must correspond to the
conclusion of an inference rule with all its premises corresponding to
the labels of its immediate predecessors. If all leaves are axioms then
we say that the derivation is a \textit{\textbf{proof}}. \footnote{In
  the proof of Theorem \ref{hilbert_seq} we can see two examples of
  derivations which are not proofs.}

The root of the tree is the \textit{\textit{conclusion}} of the whole
derivation. The set of formulas derivables in \(\mca{P}\) are its
\textit{\textbf{theorems}}. The turnstile \(\vdash\) is commonly used to
denote provability of formulas.

For a logics \((\mca{L},\,\mca{E},\,\vDash)\) and a proof system
\(\mca{P}\) defined on \(\mca{L}\), we say that \(\mca{P}\) is
\textit{\textbf{sound}} with respect to \(\mca{E}\) if every theorem is
a validity and we say that \(\mca{P}\) is \textit{\textbf{complete}}
with respect to \(\mca{E}\) if every validity is a theorem.

\hypertarget{proof-theory}{%
\subsection{Proof Theory}\label{proof-theory}}

\begin{quote}
We now assume a language \(\mca{L}\) and a proof system \(\mca{P}\) we
are interested in studying.
\end{quote}

The study of logical deductions themselves grounded the area of proof
theory, where the objects of study are formal proofs. Hilbert and
Ackermann 1928 \emph{Grundzüge der Theorischen Logik} can be seen a
departing point for Proof Theory (following Ono), while still only
working with Hilbert systems.

A rule \(\mca{R}\in\mca{P}\) is said to be \textit{\textbf{eliminable}}
if whenever there is a derivation of the premises of \(\mca{R}\) then
there is also a derivation of the conclusion not containing any
application of \(\mca{R}\).

A rule \(\mca{R}\notin\mca{P}\) is said to be
\textit{\textbf{admissible}} if whenever there is a derivation of the
premises of \(\mca{R}\) then there is also a derivation of the
conclusion.

A rule \(\mca{R}\in\mca{P}\) with a set of premises \(T\) and conclusion
\(S'\) is said to be \textit{\textbf{invertible}} when for any
\(S\in T\) the inference rule with premise \(S'\) and conclusion \(S\)
is admissible in \(\mca{P}\) or is a rule of \(\mca{P}\). \footnote{For
  example, in the proof system \(GGL_C\) we present in section 6, left
  and right introduction rules are invertible, in a similar fashion to
  the logical variant of Gentzen Calculus given by Poggiolesi
  \cite{poggi_gentzen}.} The later are called
\textit{\textbf{the inverse rules}} of the rule \(\mca{R}\).

\begin{quote}
Commonly, with sequent systems rule-sets one can not derive transitivity
of \(\vdash\). This property is a version of the cut rule, which can be
given as: For each \(X\), \(Z\) structures and \(\varphi\) formula, from
\(X\vdash \varphi\) and \(\varphi\vdash Z\) derive \(X\vdash Z\). In
spite of cut rule (which we will assume as being the last presented) not
being derivable, Gentzen showed that it was admissible, so that the
sequents derivable in the rule-set with and without cut rule were the
same.
\end{quote}

\hypertarget{display-logics}{%
\subsubsection{Display Logics}\label{display-logics}}

In his 1987 article \emph{Display Logics} Nuel D. Belnap extracted 8
sufficient conditions for replicating Gentzen's Cut-Elimination and
Subformula Property theorems. We follow the process, after introducing
some necessary extra definitions.

We call the structure on the left of a sequent its
\textit{\textbf{antecedent}} and the structure on the right of a sequent
its \textit{\textbf{consequent}}. Depending on each syntax we will have
to introduce a definition for antecedent part and consequent part. In
the particular case of atomic logics we formally introduce it in
definition \ref{ant_parts}.

\textit{\textbf{Constituents}} of an inference are the structure
occurences (substructures and its tree-node) in structures in sequents
participating in the inference, along with its position in the inference
(whether it is antecedent or consequent and which sequent is it
substructure of).

We say that a constituent \(X\) of a sequent \(S\) is
\textit{\textbf{displayed}} if \(X\) is the full antecedent or
consequent of \(S\).

We say that sequents \(S\) and \(S'\) are
\textit{\textbf{structurally equivalent}} when there is a derivation
with conclusion \(S'\) and unique premise \(S\) and there is also a
derivation with \(S'\) as the only premise and conclusion \(S\).

A calculus satisfies the \textit{\textbf{display property}} if for every
constituent \(X\) of a sequent \(S\) there is a structurally equivalent
sequent \(S'\) where \(X\) is displayed.

\textit{\textbf{Parameters}} are those constituents which occur as
substructures of structures assigned to the structure variables in the
statement of a rule.

\begin{quote}
Like the parameters \(U\), \(V\), \(W\), \((U,\,V)\) and \(X\) when
\(X:=(U,\,V),\,W\) in the rule: \[\begin{prooftree}
\hypo{X,\,X\vdash Y}
\infer1{X\vdash Y}
\end{prooftree}\] Giving the particular inference: \[\begin{prooftree}
\hypo{((U,\,V),\,W),\,((U,\,V),\,W)\vdash Y}
\infer1{((U,\,V),\,W)\vdash Y}
\end{prooftree}\]
\end{quote}

\textit{\textbf{Principal Constituents}} are non-parametric constituents
which are formulas and appear only in the conclusion of a rule.

We call a calculus \textit{\textbf{properly displayable}} if it
satisfies the display property and the following 8 properties.

\begin{itemize}
\tightlist
\item
  \(C1\) \textbf{Preservation of formulas}:\\
  All formulas occurring in premises are subformulas of some formula in
  the conclusion.
\item
  \(C2\) \textbf{Shape-alikeness of parameters}:\\
  We have an equivalence relation between parameters in an inference,
  which we call \emph{Congruence}. \emph{Congruent parameters} are at
  least ocurrences of the same structure. Conditions C3 and C4 will
  restrict the properties this relation must fulfill.
\item
  \(C3\) \textbf{Non-proliferation of parameters}:\\
  Any parameter is congruent to maximum one parameter in the conclusion
  of any inference rule.
\item
  \(C4\) \textbf{Position-alikeness of parameters}:\\
  Parameters are congruent only to parameters in a single side of
  sequents (note that this is only restricts apparition of parameters
  inside of each rule), in the sense that parameters in antecedent parts
  are only congruent to parameters in antecedent parts.
\item
  \(C5\) \textbf{Display of principal constituents}:\\
  Principal constituents in conclusions are either the entire consequent
  or the entire antecedent of the conclusion.
\item
  \(C6\) \textbf{Closure under substitution for consequent
  parameters}:\\
  When substituting simultaneously all congruent parameters of a
  consequent part in a rule for another structure the resulting
  inference is also in the rule set. Furthermore, constituents of the
  substituted parameters are also parametric and constituents not
  substituted are parametric or not depending on whether they originally
  were so.
\item
  \(C7\) \textbf{Closure under substitution for antecedent
  parameters}:\\
  When substituting simultaneously all congruent parameters of an
  antecedent part in a rule for \(M\), where \(M\) is a non-parametric
  formula consequent of the conclusion of some rule of the rule-set (in
  particular a principal constituent), the resulting inference is also
  in the rule set.\\
  Furthermore, constituents of the substituted parameters are also
  parametric and constituents not substituted are parametric or not
  depending on whether they originally were so.
\item
  \(C8\) \textbf{Eliminability of principal consituents}:\\
  For any pairs of conclusions from the rule set, \(X\vdash\varphi\) and
  \(\varphi\vdash Z\) of inferences \(I\) and \(J\), where \(\varphi\)
  is a principal constituent, either \(\varphi\) equals \(X\) or \(Z\)
  or we can infer \(X\vdash Z\) from the premises of \(I\) and \(J\),
  using the cut rule only for strict subformulas of \(\varphi\).
\end{itemize}

\bgT

In any calculus satisfying the rules \(C2\) to \(C8\) the cut rule,
expressed as follows, is admissible:

For any pair of structures \(X\) and \(Z\) and any formula \(\varphi\),
from the sequents \(X\vdash \varphi\) and \(\varphi\vdash Z\) we can
infer \(X\vdash Z\).

\ndT

\begin{proof}

We prove the statement by strong induction on the complexity of
\(\varphi\). Therefore, we want to show that if the cut rule is
admissible for all proper subformulas of \(\varphi\) then it is
admissible for \(\varphi\).

We proceed by cases, depending on whether \(\varphi\) is parametric or
not:

\begin{itemize}
\tightlist
\item
  If \(\varphi\) is non-parametric in both derivations from the premise
  the resulting derivation corresponds to the one given by \(C8\).
\item
  If \(\varphi\) were to be parametric only in the derivation
  \(\varphi\vdash Z\), we already know that whenever we have a
  derivation \(\varphi\vdash Y\) where \(\varphi\) is non-parametric
  then we can also derive \(X\vdash Y\).\\
  We recursively build a set out of the constituents congruent to the
  parametric \(\varphi\) of the derivation of \(\varphi\vdash Z\) (which
  we call set of \textit{\textbf{parametric ancestors}} of \(\varphi\)
  in the derivation of \(\varphi\vdash Z\)).\\
  We build a new inference tree by substituting in the derivation of
  \(\varphi\vdash Z\) every parametric ancestor of \(\varphi\) by \(X\).
  If suffices to show that this inference tree with conclusion
  \(X\vdash Z\) can be derived.\\
  We proceed by induction on the derivation \(\varphi\vdash Z\). Let us
  take the inference \(Inf\) of premises \(S_1,\,\ldots,\,S_l\) and
  conclusion \(U\vdash V\), falling in the set of rules. By induction
  hypothesis we take the derivations \(S_1',\,\ldots,\,S_l'\) with \(X\)
  in place of \(\varphi\). We take the sequent \(U'\vdash V'\) by
  replacing \(\varphi\) with \(X\) in all of its parametric ancestors
  which are also parametric in the conclusion \(Inf\). By definition of
  parametric ancestor and \(C3\) we have that if there is any
  constituent of \(Inf\) which is a parametric ancestor of \(\varphi\)
  then all of the constituents congruent to it in \(Inf\) also are
  parametric ancestors of \(\varphi\). By \(C4\), all the parametric
  ancestors of \(\varphi\) are antecedent parts, we can then use \(C7\)
  to show that the inference with premises \(S_1',\,\ldots,\,S_l'\) and
  conclusion \(U'\vdash V'\) falls also into the set of rules (thanks to
  the derivation \(X\vdash\varphi\), where we already know \(\varphi\)
  to be non-parametric). If all the parametric ancestors in
  \(U\vdash V\) were to be parametric also in the conclusion of \(Inf\)
  then we are done. If not, by \(C5\) and \(C4\), we have that
  \(\varphi\) covers all of \(U\) and therefore it is equal to \(U'\).
  By the second part of \(C7\), as \(U\) has not been substituted we
  know that it is non-parametric in the new inference, case which we
  have already shown.
\item
  If \(\varphi\) were to be parametric in the derivation
  \(X\vdash\varphi\), we already know that whenever we have a derivation
  \(Y\vdash\varphi\) where \(\varphi\) is non-parametric then we can
  also derive \(Y\vdash Z\).\\
  We build a new inference tree by substituting in the derivation of
  \(X\vdash\varphi\) every parametric ancestor of \(\varphi\) by \(Z\).
  If suffices to show that this inference tree with conclusion
  \(X\vdash Z\) can be derived.\\
  We proceed by induction on the derivation \(X\vdash\varphi\). Let us
  take the inference \(Inf\) of premises \(S_1,\,\ldots,\,S_l\) and
  conclusion \(U\vdash V\), falling in the set of rules. By induction
  hypothesis we take the derivations \(S_1',\,\ldots,\,S_l'\) with \(Z\)
  in place of \(\varphi\). We take the sequent \(U'\vdash V'\) by
  replacing \(\varphi\) with \(Z\) in all of its parametric ancestors
  which are also parametric in the conclusion \(Inf\). By definition of
  parametric ancestor and \(C3\) we have that if there is any
  constituent of \(Inf\) which is a parametric ancestor of \(\varphi\)
  then all the constituents congruent to it in \(Inf\) also are
  parametric ancestors of \(\varphi\). By \(C4\) all the parametric
  ancestors of \(\varphi\) are consequent parts, we can then use \(C6\)
  to show that the inference with premises \(S_1',\,\ldots,\,S_l'\) and
  conclusion \(U'\vdash V'\) falls also in the set of rules. If all the
  parametric ancestors in \(U\vdash V\) were to be parametric also in
  the conclusion of \(Inf\) then we are done. If not, by \(C5\) and
  \(C4\), we have that \(\varphi\) covers all of \(V\) and therefore it
  is equal to \(V'\). By the second part of \(C6\), as \(V\) has not
  been substituted we know that it is non-parametric in the new
  inference, case which we have already shown.
\end{itemize}

\begin{quote}
In case our rule-set satisfies rule \(C6\) also in the antecedent part,
the second step of the proof can be replaced with the same proof as in
the third step.
\end{quote}

\ndpr

\hypertarget{properties-of-display-logics}{%
\subsubsection{Properties of Display
Logics}\label{properties-of-display-logics}}

\begin{quote}
We have already seen Cut-Elimination as a consequence of \(\mca{P}\)
being properly displayable, now we are going to see what else can be
deduced from it.
\end{quote}

We say that \(\mca{P}\) satisfies the
\textit{\textbf{subformula property}} if every derivable sequent has a
proof where all formulas appearing are subformulas of some formula in
the conclusion.

\bgT

If \(\mca{P}\) satisfies \(C1\) then it has the subformula property.

\ndT

\begin{quote}
Note that the previously presented cut rule doesn't satisfy \(C1\), so
that if we want cut in a proof system and the subformula property we are
going to need also the cut-elimination theorem.
\end{quote}

\begin{proof}

We procede by induction on the deduction. Each step is a call on \(C1\).

\ndpr

We call the following rule \textit{\textbf{analytic cut}}:

If we can derive \(X\vdash\varphi\) and \(\varphi\vdash Z\) and
\(\varphi\) is a subformula of some formula in \(X\) or \(Z\), then we
can derive \(X\vdash Z\).

\begin{quote}
Cut-elimination is very handy, but not quite the unique option for
proving desired properties in proof systems. For this reason, an
alternative to cut-elimination is proving that the cut rule is
admissible in systems with analytic cut. In this context it can be seen
that the subformula property still holds.
\end{quote}

Let \(J\) be a subset of the logical connectives of \(\mca{L}\). We call
the \(J\)-\textit{\textbf{fragment}} of \(\mca{P}\) to the language only
containing the formulas for logical connectives in \(J\) (the
\(J\)-formulas) and the proof system only containing the rules for
logical connectives in \(J\).

Let \(\mca{P}\) and \(\mca{P}'\) be proof systems on languages
\(\mca{L}\) and \(\mca{L}'\supseteq\mca{L}\), respectively. We say that
\(\mca{P}'\) is a \textit{\textbf{conservative extension}} of
\(\mca{P}\) if any formula from \(\mca{L}\) derivable in \(\mca{P}'\) is
also derivable in \(\mca{P}'\).

\bgT

Whenever \(\mca{P}\) has the subformula property, \(\mca{P}\) is a
conservative extension of any \(J\)-fragment of \(\mca{P}\).

\ndT

\begin{proof}

We just take the proof given by the subformula property.

\ndpr

Another interesting property that can be obtained by cut-elimination is
Craig's Interpolation.

Let \(X\) be a sequent. We note by \(V(X)\) the variables in formulas
appearing in \(X\).

We say that \(\mca{P}\) has the \emph{Craig's Interpolation property} if
any pair of formulas \(\varphi\) and \(\rho\) such that
\(\varphi\vdash \rho\) is derivable in \(\mca{P}\), then there is a
formula \(\psi\) such that \(\varphi\vdash\psi\) and \(\psi\vdash\rho\)
are derivable and \(V(\psi)\subseteq V(\varphi)\cap V(\rho)\).

\begin{quote}
We later see a new version of it, slightly stronger, for Atomic Logics.
\end{quote}

\hypertarget{algebraic-preliminaries}{%
\section{Algebraic Preliminaries}\label{algebraic-preliminaries}}

\label{sect2}

A \textit{\textbf{monoid}} \((C,\,\cdot,\,1)\) is a set \(C\) equipped
with an associative binary operation with identity 1 such that
\(\forall a\in C,\,a\cdot 1 = 1\cdot a = a\).

A \textit{\textbf{group}} is a monoid such that every element in \(C\)
is invertible, that is, \(\forall x\in C,\,\exists y\in C\) such that
\(x\cdot y = 1\). If its operation is commutative,
\(\forall x,\,y\in C,\,x\cdot y = y\cdot x\), we say that \(C\) is an
\textit{\textbf{abelian}} group. A \textit{\textbf{subgroup}} of \(C\)
is a subset \(D\subseteq C\) such that the group operation of \(C\)
restricted to \(D\) is also a group operation.

\begin{quote}
Groups have two common notations for its operators, namely
\((\cdot,\,1)\) and \((+,\,0)\), while the later is always reserved for
abelian groups, the former can be used for both. In our case we must
take care, as we still use \(\cdot\) for Aucher's notation.
\end{quote}

Let \(C\) and \(D\) be monoids. A \textit{\textbf{morphism of monoids}}
\(\varphi:C\to D\) is a function such that \(\varphi(1) = 1\) and
\(\forall x,\,y\in C,\,\varphi(x\cdot y) = \varphi(x)\cdot\varphi(y)\).
A \textit{\textbf{morphism of groups}} \(\varphi:C\to D\) is a function
such that
\(\forall x,\,y\in C,\,\varphi(x\cdot y) = \varphi(x)\cdot\varphi(y)\).
An \textit{\textbf{isomorphism of groups}} is a bijective morphism.

Let \(C\) be a monoid and \(S\) a set. We call a
\textit{\textbf{monoid action}} \(a:C\times S\to S\) such that
\(\forall x\in S,\,a(1,\,x) = x\) and
\(\forall x\in S,\,g,\,h\in C,\,a(g,\,a(h,\,x)) = a(g\cdot h,\,x)\). We
call a \textit{\textbf{group action}} to a monoid action of a group.

Let \(C\) be a group. We note by \(\mtt{Aut}(C)\) the group of
isomorphisms \(C\to C\) with product its composition. We call it the
\textit{\textbf{automorphisms}} group.

Let \(C\) and \(D\) be groups with group operations \(\cdot_C\) and
\(\cdot_D\) respectively. The \textit{\textbf{semi-direct product}} of
\(C\) and \(D\) for a morphism \(\varphi:D\to\mtt{Aut}(C)\) is composed
of:

\begin{itemize}
\tightlist
\item
  The set \(C\times D\),
\item
  the group operation
  \((x,\,y),\,(z,\,w)\mapsto (x\cdot_C\varphi(y)(z),\,y\cdot_D w)\). It
  is noted \(C\rtimes D\).
\end{itemize}

The \textit{\textbf{direct product}} of \(C\) and \(D\) is the
semi-direct product on \(C\) and \(D\) given by morphism
\(\varphi:g\mapsto\id\).

\(\mtt{Sym}(n+1)\) is called the \textit{\textbf{symmetric group}} on
\(n+1\) elements. Its elements are the bijections on the set
\(\{1,\,\ldots,\,n+1\}\). We use for \(\mtt{Sym}(n+1)\) the operation
\(\sigma\cdot\rho = \rho\circ\sigma\), where \(\circ\) is the
composition of \(\rho\) and \(\sigma\) as functions.

Let \(C\) and \(D\) be two groups, the \textit{\textbf{free product}} of
\(C\) and \(D\), noted \(C*D\), is composed of:

\begin{itemize}
\tightlist
\item
  The set of words \(W\) on the cartesian product \(C\times D\) with
  identity \(\id_W\) the empty word, quotiented by
  \(v(g,\,h)(\id_C,\,h')w\sim v(g,\,h\cdot h')w\),
  \(v(g,\,\id_D)(g',\,h)w\sim v(g\cdot g',\,\id_D)w\) and
  \(v(\id_C,\,\id_D)w = vw\), where
  \(g,\,g'\in C,\,h,\,h'\in D,\,v,\,w\in C*D\),
\item
  the operation \(w,\,v\mapsto w\cdot v\), the concatenation of both
  words.
\end{itemize}

If \(C\) is a group and \(D\subseteq C\) is a subgroup we call a
\textit{\textbf{left}} \(D\)-\textit{\textbf{coset}} an equivalence
class of the equivalence relation given by
\(\forall x,\,y\in C,\,x\sim y\) if and only if \(x^{-1}\cdot y\in D\).
In the precedent definition, \(y\) would be a member of the \(xD\) left
coset. Similarly on the right.

If \(D\) is a \textit{\textbf{normal subgroup}}, \emph{i.e.}
\(\forall x\in C,\,xD = Dx\), its cosets form a group with product
\(xD,\,yD\mapsto xD\cdot yD = (x\cdot y)D\). We call it the
\textit{\textbf{quotient group}}.

\bgE

The subgroup
\(\mtt{Sym}(n)\simeq\{\sigma\in\mtt{Sym}(n)\mid \sigma(n+1) = n+1\}\)
provides an example for cosets which do not form a group with
element-wise multiplication. It is not a normal subgroup of
\(\mtt{Sym}(n+1)\), in accordance to what we just said. This is the
reason we need to introduce left/right cosets.

\ndE

\bgP[First Isomorphism Theorem]

Let \(C\) and \(D\) be two groups, and let \(f:C\to D\) be a morphism.
Then \(\overline{f}:f\lres\mtt{Ker}\ f \to \mtt{Im}\ f\) is an
isomorphism.

\ndP

We say that an element \(g\) of a set \(C\) with operation \(\cdot_C\)
is \textit{\textbf{involutive}} whenever \(g\cdot_C g = \id\).

\hypertarget{actions}{%
\subsection{Actions}\label{actions}}

Let \(S\) be a set and \(a:C\times S\to S\) an action. The action \(a\)
is called \textit{\textbf{free}} whenever for any \(x\in S\),
\(a(g,\,x) = x\) implies \(g = \id\). The action \(a\) is called
\textit{\textbf{faithful}} whenever \((\forall x\in S,\,a(g,\,x) = x)\)
implies \(g = \id\). The action \(a\) is called
\textit{\textbf{transitive}} whenever \(\forall x,\,y\in S\) there is
some \(g\in C\) such that \(a(g,\,x) = y\). A \textit{\textbf{regular}}
action is a free and transitive action.

\bgP

\label{lem_reg_act}

Let \(G\) be a group, \(S\) be a set and \(a:G\times S\to S\) be an
action. An action is regular if and only if
\(\forall x\in S,\,a(\,.\,,\,x)\) is a bijective function \(G\to S\).

Through this function we can define a group structure on \(S\) with
neutral \(x\).

\ndP

\begin{proof}

Indeed, we send identity to \(x\). Any other \(y\in S\) has some
\(g\in G\) such that \(a(g,\,x) = y\) so that we can send \(g\) to \(y\)
and the aplication is surjective, we call it \(\phi\). Finally, for any
two \(a(g,\,x) = y = a(h,\,x)\) we have \(g = h\) so that the
application is injective. The other implication is trivial. We define
the group operation on \(S\) by
\(y\cdot z = \phi(\phi^{-1}(y)\cdot\phi^{-1}(z))\).

\ndpr

\bgP

Free actions are regular when restricted on their orbits.

\ndP

\begin{proof}

Let \(a\) be an action of the group \(G\) on the set \(S\) and let
\(x,\,y\in S\). As \(x\) and \(y\) are in the same orbit there is some
\(g\in G\) such that \(a(g,\,x) = y\). Therefore, \(a\) restricted to
\(G\times a(G,\,x)\) is transitive.

\ndpr

Given a group \(G\), a set \(A\) and a group action
\(\alpha:G\times A\to A\), the \textit{\textbf{orbit}} of \(\alpha\)
over \(B\subseteq A\) is the image \(\alpha(G,\,B)\).

\begin{quote}
Whenever we can assume safely the group we are talking about we will
note the orbit \(\alpha(B)\). When \(B = \{x\}\) for some \(x\in S\) we
write \(\alpha(x)\).
\end{quote}

Let \(A\) be a set. The \textit{\textbf{semi-direct product for}}
\(\varphi\) \textit{\textbf{of two actions}} \(a:G\times A\to A\),
\(b:H\times A\to A\) such that
\(\forall x\in A,\,g\in G,\,h\in H,\,b(h,\,a(g,\,x)) = a(\varphi(h,\,g),\,b(h,\,x))\)
is the action
\(b\times a:G\rtimes H\times A\to A:((g,\,h),\,x)\mapsto a(g,\,b(h,\,x))\).

The \textit{\textbf{composition action}} of two commuting actions
\(a:G\times A\to A\), \(b:H\times A\to A\) is the semi-direct product of
its actions over the morphism \(g\mapsto\id\). It is noted \(a\circ b\).

\begin{quote}
\(b\times a\) is an action thanks to
\((b\times a)((g,\,h),\,(b\times a)((g',\,h'),\,x)) = a(g,\,b(h,\,a(g',\,b(h',\,x)))) = a(g\cdot_G\varphi(h,\,g),\,b(h\cdot_H h',\,x)) = (b\times a)((g,\,h)\cdot_{G\rtimes H},\,x)\).

\(a\circ b\) equals the free product action of \(b\) and \(a\).
\end{quote}

Let \(G\) and \(H\) be two groups. If \(\alpha\) and \(\beta\) are
actions of \(G\) and \(H\) on a set \(X\), then the
\textit{\textbf{free product}} action \(\alpha * \beta\) is the mapping
\(\alpha * \beta : G * H \times X \to X\) given by
\((\alpha * \beta)(g,x) = \alpha(g_1,\,\beta(h_1,\,\ldots,\,\alpha(g_n,\,\beta(h_n,x))))\),
where \(\forall i\in\{1,\,\ldots,\,n\},\,g_i\in G\), \(h_i\in H\) and
\(g = g_1h_1\ldots g_nh_n\).

Let \(C\) be a group and \(S\) be a set. Let \(a:C\times S\to S\) be an
action. The \textit{\textbf{stabilizer}} of an element \(x\in C\), noted
\(C_x\), is the set of \(g\in C\) such that \(a(g,\,x) = x\).

\bgP

\label{stab_quot}

Stabilizers are normal subgroups. Stabilizers of two elements in the
same orbit are isomorphic through conjugation. Therefore, the quotient
of the group with its stabilizers is well-defined (although it depends
on the orbit). Free actions have trivial stabilizers (isomorphic to the
trivial group, with only the identity).

\ndP

\begin{proof}

Let \(C\) be a group, \(S\) a set, \(a:C\times S\to S\) an action. Let
\(x,\,y\in S\) such that \(a(x) = a(y)\), then we take a \(g\in C\) such
that \(a(g,\,x) = y\) and we define the morphism
\(\varphi:C_x\to C_y:h\mapsto g\cdot h\cdot g^{-1}\).

\ndpr

\bgL

\label{morph_stab}

Let \(C\) and \(D\) be two groups, \(\varphi:C\to D\) and
\(\psi:D\to C\) be morphisms. Let \(S\) be a set, \(x\in S\) and
\(a:C\times S\to S\) and \(b:D\times S\to S\) be two actions such that
\(\forall g\in C,\,b(\varphi(g),\,x) = a(g,\,x)\) and
\(\forall g\in D,\,a(\psi(g),\,x) = b(g,\,x)\). Then
\(C\lres C_x\cong D\lres D_x\).

\ndL

\begin{proof}

Let \(x\in S\). We define morphisms
\(\overline\varphi:C\to D\lres D_x:g\mapsto\overline{\varphi(g)}\) and
\(\overline\psi:D\to C\lres C_x:g\mapsto\overline{\psi(g)}\).

We just have to show that the induced
\(\overline\varphi:C\lres C_x\to D\lres D_x\) and
\(\overline\psi:D\lres D_x\to C\lres C_x\) are inverses. Let \(g\in C\).
As \(b(\varphi(\psi(g)),\,x) = a(\psi(g),\,x) = b(g,\,x)\), we know that
\(\varphi(\psi(g))\in C_x\) if and only if \(g\in C_x\), so that
\(\overline\varphi\circ\overline\psi = \id\). The other order follows
from the other equality.

\ndpr

\hypertarget{modules}{%
\subsection{Modules}\label{modules}}

We call a \textit{\textbf{ring}} \((R,\,+,\,\cdot)\) an abelian group
\((R,\,+)\) along with a monoid \((R,\,\cdot)\) such that \(\cdot\) is
distributive over \(+\).

We call a (left) \(R\)-\textit{\textbf{module}} for a ring \(R\) an
abelian group \((M,\,+)\) together with a left action for the group
\(\cdot:R\times M\to M\) such that \(\cdot\) is distributive on both
sides.

A \textit{\textbf{morphism of}} \(R\)-\textit{\textbf{modules}}
\(\varphi:M\to N\) is a morphism of groups \(M\to N\) such that
\(\forall x\in R,\,m\in M,\,\varphi(x\cdot f) = x\cdot\varphi(f)\).
Again, inversible morphisms are called isomorphisms and \(M\to M\)
isomorphisms are called automorphisms.

We say that the ring \((R,\,+,\,\cdot)\) is a \textit{\textbf{field}}
whenever \((R-\{0\},\,\cdot)\) is a group. When \(R\) is a field
\(R\)-modules are called \(R\)-\textit{\textbf{spaces}} or vector
spaces. Morphisms of vector spaces are usually called
\textit{\textbf{linear transformations}}.

We will note by \(\mtt{GL}(S)\) the
\textit{\textbf{general linear group}} over \(S\) a module, which is the
group of the automorphisms on \(S\) with operation the composition.

A regular action of an \(R\)-module (with the addition group operation)
over a set will be called an \textit{\textbf{affine}}
\(R\)-\textit{\textbf{module}}. A regular action of an \(R\)-space (with
the addition group operation) over a set will be called an
\textit{\textbf{affine}} \(R\)-\textit{\textbf{space}}.

\begin{quote}
There is an isomorphism between \(n\)-matrices on a ring \(R\) and
\(R^n\) linear transformations, so that we can use them interchangeably.

Whenever we have a morphism \(\varphi:G\to\mtt{GL}(M)\) from a group
\(G\) into the linear group of an \(R\)-module \(M\) we can see
\(\varphi\) as an action of \(G\) on \(M\). For this reason we will use
the notations
\(\forall g\in G,\,m\in M,\,\varphi(g,\,m) = \varphi(g)(m)\).
\end{quote}

\hypertarget{residuation-and-negation-on-atomic-logics}{%
\section{Residuation and Negation on Atomic
Logics}\label{residuation-and-negation-on-atomic-logics}}

\label{sect3}

In this section we will introduce both Guillaume Aucher's Atomic Logics
syntax and a new way of presenting it. Although the new syntax will have
to wait until sections 4 and 5 to be fully introduced, we will compare
in this section the actions he has used and the action I now introduce.

\hypertarget{original-syntax}{%
\subsection{Original Syntax}\label{original-syntax}}

\bgD

We define the signs as elements of the unique two elements group (modulo
isomorphism), \(\mbb{B}\), where we note by \(+\) and \(\forall\) the
identity and by \(-\) and \(\exists\) the unique other number. We note
this group \(\mbb{B}\).\\
This will let us operate with signatures.

\ndD

\begin{quote}
$\ent\lres2\ent$ is to be interpreted as the ring with two elements.

We will use the following different notations for the same two-elements group (or rather isomorphic groups) and/or ring.
In this list the relation $\sim$ wants to loosely mark a correspondence in the isomorphism.
\begin{multicols}{2}
    The set is:
    \begin{itemize}
    \item \begin{tabular}{m{3em} m{.5em} m{3em} m{.5em} m{3em}}\centering$\mbb{B}$ &\centering $\sim$ &\centering $\ent\lres2\ent$ &\centering $\sim$ &\centering \texttt{Bool}\end{tabular}
    \end{itemize}
    The two elements are:
    \begin{itemize}
    \item \begin{tabular}{m{3em} m{.5em} m{3em} m{.5em} m{3em}}\centering$+$ &\centering $\sim$ &\centering $0$ &\centering $\sim$ &\centering \texttt{false}\end{tabular}
    \item \begin{tabular}{m{3em} m{.5em} m{3em} m{.5em} m{3em}}\centering$-$ &\centering $\sim$ &\centering $1$ &\centering $\sim$ &\centering \texttt{true}\end{tabular}
    \end{itemize}
    The operations are:
    \begin{itemize}
    \item \begin{tabular}{m{3em} m{.5em} m{3em} m{.5em} m{3em}}\centering$\cdot$ &\centering $\sim$ &\centering $+$ &\centering $\sim$ &\centering \texttt{XOR}\end{tabular}
    \item \begin{tabular}{m{3em} m{.5em} m{3em} m{.5em} m{3em}}& &\centering $\cdot$ &\centering $\sim$ &\centering \texttt{AND}\end{tabular}
    \end{itemize}

    \columnbreak
    \setlength{\tabcolsep}{2pt}
    \begin{center}\begin{tabular}{|c||c|c|}
    \hline
    \diagbox[width=\dimexpr \textwidth/8+2\tabcolsep]{\texttt{XOR}}{$+$} & \diagbox[width =\dimexpr \textwidth/8+2\tabcolsep, linewidth = 0pt, linecolor = white]{\texttt{false}}{$ 0$} & \diagbox[width =\dimexpr \textwidth/8+2\tabcolsep, linewidth = 0pt, linecolor = white]{\texttt{true}}{$1$}\\
    \hline
    \hline
    \diagbox[width =\dimexpr \textwidth/8+2\tabcolsep, linewidth = 0pt, linecolor = white]{\texttt{false}}{$ 0$} & \diagbox[width =\dimexpr \textwidth/8+2\tabcolsep, linewidth = 0pt, linecolor = white]{\texttt{false}}{$ 0$} & \diagbox[width =\dimexpr \textwidth/8+2\tabcolsep, linewidth = 0pt, linecolor = white]{\texttt{true}}{$1$}\\
    \hline
    \diagbox[width =\dimexpr \textwidth/8+2\tabcolsep, linewidth = 0pt, linecolor = white]{\texttt{true}}{$1$} & \diagbox[width =\dimexpr \textwidth/8+2\tabcolsep, linewidth = 0pt, linecolor = white]{\texttt{true}}{$1$} & \diagbox[width =\dimexpr \textwidth/8+2\tabcolsep, linewidth = 0pt, linecolor = white]{\texttt{false}}{$ 0$}\\
    \hline
    \end{tabular}
    
    \begin{tabular}{|c||c|c|}
    \hline
    \diagbox[width =\dimexpr \textwidth/8+2\tabcolsep]{\texttt{AND}}{$\cdot$} & \diagbox[width =\dimexpr \textwidth/8+2\tabcolsep, linewidth = 0pt, linecolor = white]{\texttt{false}}{$ 0$} & \diagbox[width =\dimexpr \textwidth/8+2\tabcolsep, linewidth = 0pt, linecolor = white]{\texttt{true}}{$ 1$}\\
    \hline
    \hline
    \diagbox[width =\dimexpr \textwidth/8+2\tabcolsep, linewidth = 0pt, linecolor = white]{\texttt{false}}{$ 0$} & \diagbox[width =\dimexpr \textwidth/8+2\tabcolsep, linewidth = 0pt, linecolor = white]{\texttt{false}}{$ 0$} & \diagbox[width =\dimexpr \textwidth/8+2\tabcolsep, linewidth = 0pt, linecolor = white]{\texttt{false}}{$ 0$}\\
    \hline
    \diagbox[width =\dimexpr \textwidth/8+2\tabcolsep, linewidth = 0pt, linecolor = white]{\texttt{true}}{$1$} & \diagbox[width =\dimexpr \textwidth/8+2\tabcolsep, linewidth = 0pt, linecolor = white]{\texttt{false}}{$ 0$} & \diagbox[width =\dimexpr \textwidth/8+2\tabcolsep, linewidth = 0pt, linecolor = white]{\texttt{true}}{$ 1$}\\
    \hline
    \end{tabular}\end{center}
\end{multicols}

The operation for $\mathtt{Bool}$ is the one used in the code.
We also note $\forall$ for $+\in\mbb{B}$ and $\exists$ for $-\in\mbb{B}$.

I use $\nat^* = \nat-\{0\}$.
\end{quote}

Let \(n\in\mbb{N}^*\), we say that the elements of
\(\mbb{C}_n := \mtt{Sym}(n+1)\times\{+,\,-\}\times\{\forall,\,\exists\}\times(\nat^*)^{n+1}\times\{+,\,-\}^{n}\)
are the \(n\)-\textit{\textbf{ary connective skeletons}}. The \(0\)-ary
connectives are also called \textit{\textbf{atom skeletons}} as they'll
be fundamental in the semantics. The connective skeleton set is
\(\bigcup_{n\in\nat}\mbb{C}_n\).

\begin{figure}\begin{center}

\includegraphics{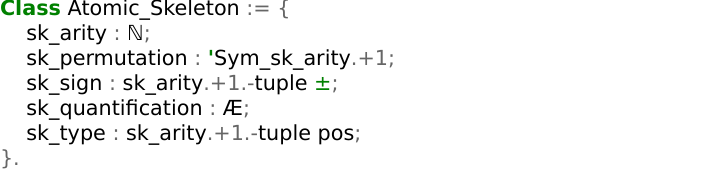}

\caption{In Coq we have introduced this version of the Skeleton Connectives}

\end{center}\end{figure}

So, the elementary objects of Atomic Logics are composed of two parts:

\begin{itemize}
\tightlist
\item
  A \textbf{signature}, consisting of:

  \begin{itemize}
  \tightlist
  \item
    a \textbf{semantic sign}; its second component, assigned to a \(+\)
    or a \(-\).\\
    The projection to the second component is noted
    \(\pm:\mbb{C}\to\mbb{B}\).
  \item
    a \textbf{quantification sign}; its third component, assigned to a
    \(\forall\) or an \(\exists\).\\
    The projection to the third component is noted
    \(\mtt{\AE}:\mbb{C}\to\mbb{B}\).
  \item
    a \textbf{tonicity signature}; its fifth component, assigned to a
    tuple of \(+\) and \(-\).\\
    For any and \(i\in\{1,\,\ldots,\,n\}\), the projections to the
    \(i\)-essime component of the fifth component are noted
    \(\pm_i:\mbb{C}_n\to\mbb{B}\).
  \item
    a \textbf{type signature}; its fourth component, assigned to a tuple
    of positive numbers.\\
    For any \(i\in\{1,\,\ldots,\,n+1\}\) the projections to the
    \(i\)-essime component of the fourth component are noted
    \(k_i:\mbb{C}_n\to\nat^*\). For connective skeletons of arity \(n\),
    we will note \(k\) for \(k_{n+1}\) and call it the \textbf{output
    type} of the connective skeleton.
  \end{itemize}
\item
  A \textbf{permutation}; its first component, assigned to a
  permutation.
\end{itemize}

\begin{quote}
These parameters will allow us to define the truth conditions of the
connectives of atomic logics in section \ref{semantics}.

\begin{itemize}
\tightlist
\item
  The semantic sign describes the assertivity of the connective
  skeleton: does it refer to the accessible worlds or to the
  complementary relation's accessible worlds?
\item
  The quantification sign defines the quantification of the truth value:
  is it true in all worlds or in at least some of them?
\item
  The tonicity signature is composed of the signs of each component (is
  it true on some worlds or on the complementary's - in boolean logic,
  is it true or is it false required in this component to get a true
  result?).
\item
  The type signature is composed of two parts, the output type or simply
  type of the connective skeleton itself and the types of the
  components. It defines the nature of the truth values.
\item
  The permutation will mark the positions of the components in the
  accessibility relation for the semantics.
\end{itemize}
\end{quote}

\bgD

\label{alpha1}

We define
\(\alpha_n:\mtt{Sym}(n+1)\times\mbb{C}_n\to\mbb{C}_n:(\rho,\,C)\mapsto\alpha_n(\rho,\,C)\):

\begin{itemize}
\item If $\rho(n+1) = n+1$, then the action maps:
      $$\begin{matrix}
      (\sigma,\,\pm,\,
          \text{\AE},\,
          (k_1,\,\ldots,\,k_n,\,k_{n+1}),\,
          (\pm_1,\,\ldots,\,\pm_n))\\
          \rotatebox[origin=c]{-90}{$\mapsto$}\\
      (\rho\cdot\sigma,\,\pm,\,
          \text{\AE},\,
          (k_{\rho(1)},\,\ldots,\,k_{\rho(n)},\,k_{\rho(n+1)}),\,
          (\pm_{\rho(1)},\,\ldots,\,\pm_{\rho(n)}))
      \end{matrix}$$
\item When $\rho(n+1)\neq n+1$, then the action maps:
      
      {\noindent\begin{minipage}{\textwidth}
      $$\begin{matrix}
      (\sigma,\,\pm,\,
          \text{\AE},\,
          (k_1,\,\ldots,\,k_n,\,k_{n+1}),\,
          (\pm_1,\,\ldots,\,\pm_n))\\
          \rotatebox[origin=c]{-90}{$\mapsto$}
      \end{matrix}$$
      \vspace{-2.5em}
      \begin{align*}(\rho\cdot\sigma,\,
          -\pm_{\rho(n+1)}\pm,\,
          -&\pm_{\rho(n+1)}\text{\AE},\,
          (k_{\rho(1)},\,\ldots,\,k_{\rho(n)},\,k_{\rho(n+1)}),\,\\
        (-&\pm_{\rho(n+1)}\pm_{(n+1\ \rho(n+1))(\rho(1))},\,\ldots,\,-\pm_{\rho(n+1)}\pm_{(n+1\ \rho(n+1))(\rho(i-1))},\,\\
        &\pm_{(n+1\ \rho(n+1))(\rho(\rho^-(n+1)))},\,\\
        -&\pm_{\rho(n+1)}\pm_{(n+1\ \rho(n+1))(\rho(i+1))},\,\ldots,\,-\pm_{\rho(n+1)}\pm_{(n+1\ \rho(n+1))(\rho(n))}))
      \end{align*}
      \end{minipage}}
\end{itemize}

We will note the image \(\alpha_n(\sigma,\,\star)\) as \(\sigma\star\).

\ndD

\begin{figure}\begin{center}

\includegraphics{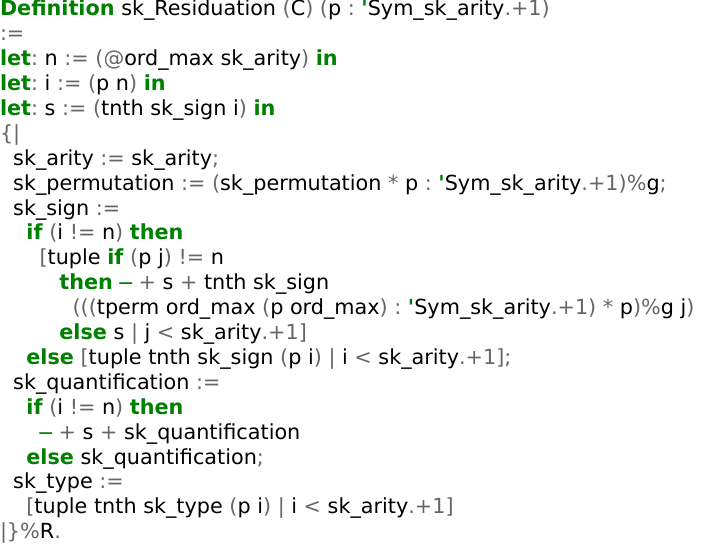}

\caption{Implementation of $\alpha$ in Coq.\\
  In Coq we have used the group operation $\circ$ for $\mtt{Sym}(n+1)$, this is because our group actions are left actions and in the libraries all group actions are right actions.
  The only change this supposes is the order of the group operation.
}

\end{center}\end{figure}

\begin{itemize}
\item
  We now introduce the action \(\alpha\) as presented in Aucher's
  \emph{Display and Hilbert Calculi for Atomic and Molecular Logics}.

  We define the function
  \(\alpha'_n:\mtt{Sym}(n+1)\times\mbb{C}_n\to\mbb{C}_n:(\tau,\,\star)\mapsto\tau\star\)
  inductively as follows. We refer to Aucher's articles for its
  well-definedness. Let
  \(\star = (\sigma,\,\pm,\,\mtt{\AE},\,\overline{k},\,(\pm_1,\,\ldots,\,\pm_n))\in\mbb{C}_n\).

  \begin{itemize}
  \tightlist
  \item
    If \(\tau\) is the transposition \((j\ n+1)\), then
    \[\alpha'_n(\tau,\,\star) = ((j\ n+1)\cdot\sigma,\,-\pm_j\pm,\,-\pm_j\mtt{\AE},\,(k_1,\,\ldots,\,k_{n+1},\,\ldots,\,k_{j}),\,(-\pm_j\pm_1,\,\ldots,\,-\pm_j\pm_n))\]
  \item
    If \(\tau\) is the cycle \((a_1\ \ldots\ a_k\ n+1)\), then
    \[\alpha'_n(\tau,\,\star) = \alpha'_n((a_1\ n+1),\,\ldots\alpha'_n((a_k\ n+1),\,\star))\]
  \item
    If \(\tau\) is a cycle fixing \(n+1\), then
    \[\alpha'_n(\tau,\,\star) = (\tau\cdot\sigma,\,\pm,\,\mtt{\AE},\,(k_{\tau(1)},\,\ldots,\,k_{\tau(n)}),\,(\pm_{\tau(1)},\,\ldots,\,\pm_{\tau(n)}))\]
  \item
    Finally, if \(\tau\) is an arbitrary permutation of
    \(\mtt{Sym}(n+1)\), it can be factorized into a product of disjoint
    cycles \(\tau = c_1\cdot\ldots\cdot c_k\) and this factorization is
    unique (modulo reordering). So, we define
    \(\alpha'_n(\tau,\,\star) = \alpha'_n(c_1,\,(\ldots\alpha_n'(c_k,\,\star)))\).
  \end{itemize}
\item
  Secondly, we need to define boolean negation as the action
  \(\beta_n:\ent\lres2\ent\times\mbb{C}_n\to\mbb{C}_n:(b,\,C)\mapsto\beta(b,\,C)\)
  defined by:

  \begin{itemize}
  \tightlist
  \item
    if \(b = +\), then the action maps \(\star\) to \(\star\).
  \item
    if \(b = -\), then the action maps
    \((\sigma,\,\pm,\,\text{\AE},\,(k,\,k_1,\,\ldots,\,k_n),\,(\pm_1,\,\ldots,\,\pm_{n}))\)
    to:
    \[(\sigma,\,-\pm,\,-\text{\AE},\,(k,\,k_1,\,\ldots,\,k_n),\,(-\pm_1,\,\ldots,\,-\pm_{n}))\]
  \end{itemize}

  We will note the image \(\beta_n(\pm,\,\star)\) as \(\pm\star\).
\item
  Closely related there is the dual action as
  \(\delta_n:\ent\lres2\ent\times\mbb{C}_n\to\mbb{C}_n:(b,\,C)\mapsto\delta(b,\,C)\)
  defined by:

  \begin{itemize}
  \tightlist
  \item
    if \(b = +\), then the action maps \(\star\) to \(\star\).
  \item
    if \(b = -\), then the action maps
    \((\sigma,\,\pm,\,\mtt{\AE},\,(k,\,k_1,\,\ldots,\,k_n),\,(\pm_1,\,\ldots,\,\pm_{n}))\)
    to:
    \[(\sigma,\,-\pm,\,-\text{\AE},\,(k,\,k_1,\,\ldots,\,k_n),\,(\pm_1,\,\ldots,\,\pm_{n}))\]
  \end{itemize}
\item
  Finally, we present the switches as
  \(\forall i\in\{1,\,\ldots,\,n\},\,s_i:(\ent\lres2\ent)^{n+1}\times\mbb{C}_n\to\mbb{C}_n:(b,\,C)\mapsto s_i(b,\,C)\)
  defined by:

  \begin{itemize}
  \tightlist
  \item
    if \(b = +\), then the action maps \(\star\) to \(\star\).
  \item
    if \(b = -\), then the action maps
    \((\sigma,\,\pm,\,\mtt{\AE},\,(k,\,k_1,\,\ldots,\,k_n),\,(\pm_1,\,\ldots,\,\pm_n))\)
    to:
    \[(\sigma,\,\pm,\,\mtt{\AE},\,(k,\,k_1,\,\ldots,\,k_n),\,(\pm_1,\,\ldots,\,-\pm_i,\,\ldots,\,\pm_{n}))\]
  \end{itemize}
\end{itemize}

\bgP

\(\alpha_n\) is an action.

\ndP

\begin{proof}

Checking \(\alpha_n(\id,\,\star) = \star\) is trivial.

We now want to prove that
\(\alpha_n(\sigma_1,\,\alpha_n(\sigma_2,\,\star)) = \alpha_n(\sigma_1\cdot\sigma_2,\,\star)\).
We work on cases on the value of \(\sigma_1(n+1)\), \(\sigma_2(n+1)\)
and \(\sigma_2(\sigma_1(n+1))\).

\begin{itemize}
\tightlist
\item
  We suppose that none of them equals \(n+1\).

  \begin{itemize}
  \tightlist
  \item
    To illustrate the method we show that both sides of the equality
    have the same sign \(\pm\). \begin{align*}
    &\pm(\alpha_n(\sigma_1,\,\alpha_n(\sigma_2,\,\star)))\\
    =&-\pm_{\sigma_1(n+1)}(\alpha_n(\sigma_2,\,\star))\pm(\alpha_n(\sigma_2,\,\star))\\
    =&-(-\pm_{\sigma_2(n+1)}(\star)\pm_{(n+1\ \sigma_2(n+1))(\sigma_2(\sigma_1(n+1)))}(\star))(-\pm_{\sigma_2(n+1)}(\star)\pm(\star))\\
    =&-\pm_{\sigma_2(\sigma_1(n+1))}(\star)\pm(\star)\\
    =&-\pm_{(\sigma_1\cdot\sigma_2)(n+1)}(\star)\pm(\star)\\
    =&\pm(\alpha_n(\sigma_1\cdot\sigma_2,\,\star))
    \end{align*}
  \item
    Now we do the proof for \(\pm_i\) for \(i\leq n\).\\
    If \(\sigma_1(i)\neq n+1\) and \(\sigma_2(\sigma_1(i))\neq n+1\).
    \begin{align*}
    &\pm_i(\alpha_n(\sigma_1,\,\alpha_n(\sigma_2,\,\star)))\\
    =&-\pm_{\sigma_1(n+1)}(\alpha_n(\sigma_2,\,\star))\pm_{(n+1\ \sigma_1(n+1))(\sigma_1(i))}(\alpha_n(\sigma_2,\,\star))\\
    =&-(-\pm_{\sigma_2(n+1)}(\star)\pm_{(n+1\ \sigma_2(n+1))(\sigma_2(\sigma_1(n+1)))}(\star))\\
    &(-\pm_{\sigma_2(n+1)}(\star)\pm_{(n+1\ \sigma_2(n+1))\sigma_2((n+1\ \sigma_1(n+1))(\sigma_1(i)))}(\star))\\
    =&-\pm_{\sigma_2(\sigma_1(n+1))}(\star)\pm_{(n+1\ \sigma_2(n+1))\sigma_2((n+1\ \sigma_1(n+1))(\sigma_1(i)))}(\star)\\
    =&-\pm_{(\sigma_1\cdot\sigma_2)(n+1)}(\star)\pm_{(n+1\ \sigma_2(\sigma_1(n+1)))\sigma_2(\sigma_1(i))}(\star)\\
    =&\pm_i(\alpha_n(\sigma_1\cdot\sigma_2,\,\star))
    \end{align*} The other cases are similar, for example when
    \(\sigma_1(i) = n+1\) we use that
    \[\pm_{(n+1\ \sigma_1(n+1))(\sigma_1(i))}(\alpha_n(\sigma_2,\,\ast)) = \pm_{\sigma_1(n+1)}(\alpha_n(\sigma_2,\,\ast))\]
  \end{itemize}
\item
  Whenever \(\sigma_1(\sigma_2(n+1)) = n+1\) we should note that
  \[\pm_{\sigma_1(n+1)}(\alpha_n(\sigma_2,\,\star)) = \pm_{(n+1\ \sigma_2(n+1))(\sigma_2(\sigma_1(n+1)))}(\star)\]
\item
  Showing the other equalities is a similar procedure.
\end{itemize}

The full proof can be seen in the appendix and the code available in
GitHub, written in the SSReflect extension of Coq.

\ndpr

\begin{figure}\begin{center}

\includegraphics{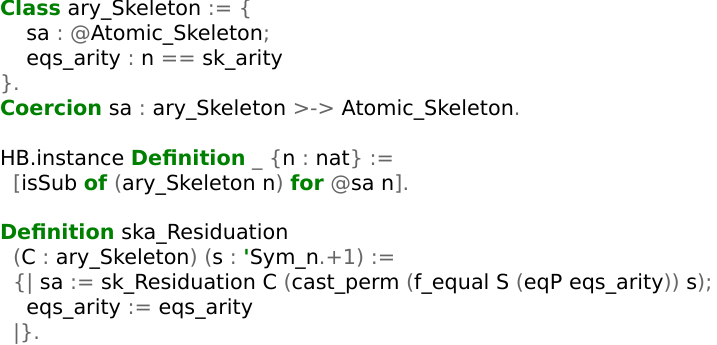}
\includegraphics{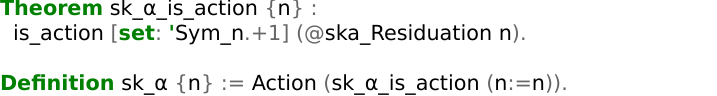}

\caption{We have provided a proof of $\alpha_n$ being an action}

\end{center}\end{figure}

\bgP

\label{origin_alph}

\(\alpha_n\) and \(\alpha'_n\) are equal on
\(\mtt{Sym}(n+1)\times\mbb{C}_n\).

\ndP

\begin{proof}

Now we check that
\(\forall\sigma\in\mtt{Sym}(n+1),\,\star\in\mbb{C}_n,\,\alpha_n(\sigma,\,\star) = \alpha'_n(\sigma,\,\star)\).
We will proceed by induction on the decomposition of \(\sigma\) by
transpositions.

\begin{itemize}
\tightlist
\item
  For the base cases, involving transpositions with \(n+1\) and cycles
  fixing \(n+1\), it is trivial
\item
  If \(\sigma\in\mtt{Sym}(n+1)\) satisfies the induction hypothesis, we
  take \(j\in\{1,\,\ldots,\,n\}\): \begin{align*}
    &\alpha_n((j\ n+1)\cdot\sigma,\,\star)\\
    =&\alpha_n((j\ n+1),\,\alpha_n(\sigma,\,\star))\\
    =&\alpha_n((j\ n+1),\,\alpha_n'(\sigma,\,\star))\\
    =&\alpha_n'((j\ n+1),\,(\alpha_n'(\sigma,\,\star)))\\
    =&\alpha_n'((j\ n+1)\cdot\sigma,\,\star)
  \end{align*} Where the first equality comes from the new \(\alpha_n\)
  being an action, the second and third equalities come from the
  induction hypothesis and finally, the fourth action comes from
  Auchers's \(\alpha'_n\) also being an action.
\end{itemize}

\ndpr

We can also work on consecutive applications of actions \(\alpha_n\) and
\(\beta_n\) through its free product action \(\alpha_n*\beta_n\),
introduced in the Algebraic Preliminaries.

\begin{quote}
We sometimes note \(\sigma\ostar\) for \(\alpha_n(\sigma,\,\ostar)\) and
\(s\ostar\) for \(\beta_n(s,\,\ostar)\).
\end{quote}

\newpage

\hypertarget{the-main-result-a-finite-version-of-the-group-actions}{%
\subsection{The Main Result: A Finite Version of the Group
Actions}\label{the-main-result-a-finite-version-of-the-group-actions}}

\begin{quote}
The main motivation for this detour is that the free product action
\(\alpha_n*\beta_n\) used for Atomic Logics works with infinite groups,
which is undesirable for working on it in Coq. We are going to redefine
the free product action \(\alpha_n*\beta_n\) into another group action
\(\alpha_n\times\varsigma_n\) which resorts only to finite groups, using
the notion of semi-direct product. While the original syntax worked
through all the possible sequences of compositions of the pair of
actions, our new syntax will describe all connectives without repeating
itself.
\end{quote}

In this section we will introduce a new morphism \(Q\), in subsection
\textit{Morphism $Q$ for tonicity signatures}, which has some of the
interesting properties of \(\alpha\) while being a morphism into the
general linear group of tonicity signatures, we get back our \(\alpha\)
in subsection \textit{A new action $\alpha$ from $Q$}, in subsection
\textit{Correspondence between the $\alpha$} we check that they are all
the same and, finally, in subsection
\textit{The free product and the semi-direct product} we provide a
translation between the connectives given by \(\alpha*\beta\), defined
over a free group, and some action \(\alpha\times\varsigma\), defined
over a finite group. For details on how Theorem \ref{quot_sign} implies
equivalence for the resulting atomic logics, we invite the reader to go
to definition \ref{equiv_fam}, proposition \ref{isom_equi} and corollary
\ref{semi_free_equi}.

\bgD

For any \(n\in\nat^*\), the \textit{\textbf{simplified}}
\(n\)-\textit{\textbf{ary connective skeletons}} are the elements of
\(\mbb{D}_n := \mtt{Sym}(n+1)\times(\nat^*)^{n+1}\times\{+,\,-\}^{n+1}\times\{+,\,-\}\).
We note the projections to the \(i\)-essime component of the fourth
component \(\pm_i\).

\ndD

The last component is not going to correspond to the sign \AE{} but
rather to the sum \(\pm+\mtt{\AE}\), where \AE{} is meant to represent
the quantification sign and \(\pm\) the semantical sign from
\(\mbb{C}_n\).

\bgD

\label{skel_trans}

The translations back and forth from the previous data structure
\(\mbb{C}_n\) are:

\begin{itemize}

\item $\iota^{-1}:(\sigma,\,k,\,\pm,\,s)\mapsto(\sigma,\,\pm_{n+1},\,\pm_{n+1}+ s,\,k,\,(\pm_1,\,\ldots,\,\pm_n))$.
\item $\iota:(\sigma,\,\pm,\,\mtt{\AE},\,k,\,\pm)\mapsto (\sigma,\,k,\,(\pm_1,\,\ldots,\,\pm_n,\,\pm),\,\pm+\mtt{\AE})$.

\end{itemize}

\ndD

\begin{quote}
From now on we will use \(\mbb{C}\) also for refering to \(\mbb{D}\). On
the follow, unless explicitly stated we will assume to be working with
\(\mbb{D}\) with signs \(\pm_i\), \(\mtt{\AE}\), \(\pm\) living on
\(\ent\lres2\ent\). Furthermore, all actions we will work with do not
change \(\pm+\mtt{\AE}\), which lets us define them in a simpler way.
\end{quote}

\bgD

On the new structure skeleton structures, we define the switch action as
\[\begin{matrix}\varsigma:&(\ent\lres2\ent)^{n+1}\times\mbb{D}_n&\to&\mbb{D}_n\\&(b,\,\begin{bmatrix}\sigma'\\ k\\ v\\ s\end{bmatrix})&\mapsto&b + \begin{bmatrix}\sigma'\\ k\\ v\\ s\end{bmatrix} = \begin{bmatrix}\sigma'\\ k\\ b + v\\ s\end{bmatrix}\end{matrix}\]
It is again a free action, which makes its orbits along the connectives
skeletons affine \((\ent\lres2\ent)^{n+1}\)-spaces. In fact, most
actions on connective skeletons we will work with are free.

\ndD

We define the \textit{\textbf{Kronecker delta}} (although we use the
same greek letter for the dual action, note the difference in the
notation)
\[\delta_{i,\,j} := \left\{\begin{matrix}1 & \mtt{if }i = j\\ 0 & \mtt{if }i \neq j\end{matrix}\right.\]

\bgD

We define a function
\(P_{\nat^*}:\mtt{Sym}(n+1)\hookrightarrow\mtt{Sym}({\nat^*}^{n+1})\):
\[P_{\nat^*}(\sigma)(n_1,\,\ldots,\,n_m) := (n_{\sigma(1)},\,\ldots,\,n_{\sigma(m)})\]
We define a function
\(P:\mtt{Sym}(n+1)\hookrightarrow\mtt{Sym}((\ent\lres2\ent)^{n+1})\),
given by the matrix: \[P(\sigma)_{i,\,j} := \delta_{\sigma(i),\,j}\]

\ndD

\bgP

\label{p_act}

\(P\) is an injective morphism
\(\mtt{Sym}(n+1)\hookrightarrow\mtt{GL}(\ent\lres2\ent)^{n+1}\).
\(P_{\nat^*}\) is an injective morphism.

\ndP

\begin{proof}

Checking the requirement for the identity is trivial.
\begin{align*}(P(\sigma)P(\sigma'))_{i,\,j}
&= \sum_k P(\sigma)_{i,\,k}P(\sigma')_{k,\,j} = \sum_k \delta_{\sigma(i),\,k}\delta_{\sigma'(k),\,j} = \sum_k \delta_{\sigma(i),\,k}\delta_{k,\,\sigma'^{-1}(j)}\\
&= \delta_{\sigma(i),\,\sigma'^{-1}(j)} = \delta_{\sigma'(\sigma(i)),\,j} = \delta_{(\sigma\cdot\sigma')(i),\,j} = P(\sigma\cdot\sigma')_{i,\,j}\end{align*}

For any \(\sigma\in\mtt{Sym}(n+1)\) we see that \(P(\sigma)\) is
bijective by showing \(P(\sigma)P(\sigma^{-1}) = P(\id) = \id\). Indeed
\(P(\id)_{i,\,j} = \delta_{\id(i),\,j} = \id_{i,\,j}\).

For any \(\sigma\in\mtt{Sym}(n+1),\,P(\sigma) = \id\) implies that
\(\delta_{\sigma(i),\,j} = \delta_{i,\,j}\) so that \(\sigma(i) = i\)
and \(\sigma = \id\).

\ndpr

\hypertarget{morphism-q-for-tonicity-signatures}{%
\subsubsection*{\texorpdfstring{Morphism \(Q\) for tonicity
signatures}{Morphism Q for tonicity signatures}}\label{morphism-q-for-tonicity-signatures}}

We will introduce a morphism which will later be useful. It can be seen
as a conjugation of \(P\).

\bgD

\label{q_t}

We define the morphism
\(Q:\mtt{Sym}(n+1)\to\mtt{Sym}(\ent\lres2\ent)^{n+1}\), given by the
matrix \begin{align*}
Q(\sigma)_{i,\,j} :=&
  \delta_{(n+1\ \sigma(n+1))(\sigma(i)),\,j}
    + \delta_{\sigma(n+1),\,j}
        (1 + \delta_{n+1,\,\sigma(i)})
        (1 + \delta_{\sigma(n+1),\,n+1})\\
=& \left\{\begin{matrix}\delta_{(n+1\ \sigma(n+1))(\sigma(i)),\,j}&\mtt{if }j\notin\{n+1,\,\sigma(n+1)\}\\1 & \mtt{if }j\neq n+1,\,j = \sigma(n+1)\\\delta_{i,\,j} & \mtt{else}\end{matrix}\right.
\end{align*}

We define the linear transformation
\(T:(\ent\lres2\ent)^{n+1}\to(\ent\lres2\ent)^{n+1}\), given by the
matrix \[
T_{i,\,j} := \delta_{i,\,j} + \delta_{j,\,n+1}(1 + \delta_{i,\,n+1})
\]

\ndD

\begin{quote}
We recall the action notation \(Q(\sigma,\,v) = Q(\sigma)(v)\)
introduced in the algebraic preliminaries.
\end{quote}

\bgP

\label{perm_q_t}

Let \(n\in\nat\). For any \(\sigma\in\mtt{Sym}(n+1)\), the equality
\(T\cdot Q(\sigma) = P(\sigma)\cdot T\) holds.

\ndP

\begin{proof}

The proof makes use of Kronecker delta's properties and that sum is
involutive in \(\ent\lres2\ent\). \begin{align*}
&(T\cdot Q(\sigma))_{i,\,j}=\sum_k T_{i,\,k}Q(\sigma)_{k,\,j}\\
=&\sum_k(\delta_{i,\,k} + \delta_{k,\,n+1}(1 + \delta_{i,\,n+1}))(\delta_{(n+1\ \sigma(n+1))(\sigma(k)),\,j} + \delta_{\sigma(n+1),\,j}(1 + \delta_{n+1,\,\sigma(k)})(1 + \delta_{\sigma(n+1),\,n+1}))\\
=&\sum_k\delta_{i,\,k}(\delta_{(n+1\ \sigma(n+1))(\sigma(k)),\,j} + \delta_{\sigma(n+1),\,j}(1 + \delta_{n+1,\,\sigma(k)})(1 + \delta_{\sigma(n+1),\,n+1}))\\
&\hspace{1em}+\delta_{k,\,n+1}(1 + \delta_{i,\,n+1})(\delta_{(n+1\ \sigma(n+1))(\sigma(k)),\,j} + \delta_{\sigma(n+1),\,j}(1 + \delta_{n+1,\,\sigma(k)})(1 + \delta_{\sigma(n+1),\,n+1}))\\
=&\delta_{(n+1\ \sigma(n+1))(\sigma(i)),\,j} + \delta_{\sigma(n+1),\,j}(1 + \delta_{n+1,\,\sigma(i)})(1 + \delta_{\sigma(n+1),\,n+1})\\
&\hspace{1em}+ (1 + \delta_{i,\,n+1})(\delta_{(n+1\ \sigma(n+1))(\sigma(n+1)),\,j} + \delta_{\sigma(n+1),\,j}(1 + \delta_{n+1,\,\sigma(n+1)})(1 + \delta_{\sigma(n+1),\,n+1}))
\end{align*} Now we will use that
\(\forall i,\,j,\,k,\,l\in\{1,\,\ldots,\,n+1\},\, \delta_{(i\ j)(k),\,l} = \delta_{k,\,(i\ j)(l)} = \delta_{k,\,l}(1 + \delta_{i,\,l})(1 + \delta_{j,\,l})+ \delta_{k,\,j}\delta_{i,\,l}(1 + \delta_{j,\,l}) + \delta_{k,\,i}\delta_{j,\,l}\),
which can be seen by doing case analysis on \(\delta_{i,\,l}\) and
\(\delta{j,\,l}\). \begin{align*}
&(T\cdot Q(\sigma))_{i,\,j}
=\delta_{\sigma(i),\,j}(1 + \delta_{n+1,\,j})(1 + \delta_{\sigma(n+1),\,j})
+ \delta_{\sigma(i),\,\sigma(n+1)}\delta_{n+1,\,j}(1 + \delta_{\sigma(n+1),\,j})\\
&\hspace{1em}+ \delta_{\sigma(i),\,n+1}\delta_{\sigma(n+1),\,j}+ \delta_{\sigma(n+1),\,j}(1 + \delta_{n+1,\,\sigma(i)})(1 + \delta_{\sigma(n+1),\,n+1})\\
&\hspace{1em}+ (1 + \delta_{i,\,n+1})(\delta_{n+1,\,j} + \delta_{\sigma(n+1),\,j}(1 + \delta_{n+1,\,\sigma(n+1)}))\\
=&\delta_{\sigma(i),\,j}+\delta_{\sigma(i),\,j}\delta_{n+1,\,j}+\delta_{\sigma(i),\,j}\delta_{\sigma(n+1),\,j}+\delta_{\sigma(i),\,j}\delta_{n+1,\,j}\delta_{\sigma(n+1),\,j}\\
&\hspace{1em}+ \delta_{i,\,n+1}\delta_{n+1,\,j}(1 + \delta_{\sigma(n+1),\,j})
+ \delta_{\sigma(i),\,n+1}\delta_{\sigma(n+1),\,j}\\
&\hspace{1em}+ (1 + \delta_{i,\,n+1})\delta_{n+1,\,j} + (1 + \delta_{n+1,\,\sigma(i)} + 1 + \delta_{i,\,n+1})\delta_{\sigma(n+1),\,j}(1 + \delta_{n+1,\,\sigma(n+1)})\\
=&\delta_{\sigma(i),\,j}+\delta_{\sigma(i),\,j}\delta_{n+1,\,j}+\delta_{\sigma(i),\,\sigma(n+1)}\delta_{\sigma(n+1),\,j}+\delta_{\sigma(i),\,\sigma(n+1)}\delta_{n+1,\,j}\delta_{\sigma(n+1),\,j}\\
&\hspace{1em}+ \delta_{i,\,n+1}\delta_{n+1,\,j}\delta_{\sigma(n+1),\,j} %a
+ \delta_{\sigma(i),\,n+1}\delta_{\sigma(n+1),\,j}\\
&\hspace{1em}+ \delta_{n+1,\,j} + (\delta_{n+1,\,\sigma(i)}+\delta_{n+1,\,\sigma(i)}\delta_{n+1,\,\sigma(n+1)} + \delta_{i,\,n+1} + \delta_{i,\,n+1}\delta_{n+1,\,\sigma(n+1)})\delta_{\sigma(n+1),\,j}\\
=&\delta_{\sigma(i),\,j}+\delta_{\sigma(i),\,j}\delta_{n+1,\,j}+\delta_{i,\,n+1}\delta_{\sigma(n+1),\,j}
+ \delta_{\sigma(i),\,n+1}\delta_{\sigma(n+1),\,j}\\
&\hspace{1em}+ \delta_{n+1,\,j} + (\delta_{n+1,\,\sigma(i)}+\delta_{n+1,\,\sigma(i)}\delta_{\sigma(i),\,\sigma(n+1)} + \delta_{i,\,n+1} + \delta_{i,\,n+1}\delta_{n+1,\,\sigma(i)})\delta_{\sigma(n+1),\,j}\\
=&\delta_{\sigma(i),\,j}+\delta_{\sigma(i),\,n+1}\delta_{n+1,\,j}+\delta_{i,\,n+1}\delta_{\sigma(n+1),\,j}
+ \delta_{\sigma(i),\,n+1}\delta_{\sigma(n+1),\,j}\\
&\hspace{1em}+ \delta_{n+1,\,j} + (\delta_{n+1,\,\sigma(i)} + \delta_{i,\,n+1})\delta_{\sigma(n+1),\,j}
=\delta_{\sigma(i),\,j} + \delta_{j,\,n+1}(1 + \delta_{\sigma(i),\,n+1})
\end{align*} On the other side we have, as desired: \begin{align*}
(P(\sigma)\cdot T)_{i,\,j} =&\sum_k P(\sigma)_{i,\,k}T_{k,\,j} = \sum_k \delta_{\sigma(i),\,k}(\delta_{k,\,j} + \delta_{j,\,n+1}(1 + \delta_{k,\,n+1}))\\
=&\delta_{\sigma(i),\,j} + \delta_{j,\,n+1}(1 + \delta_{\sigma(i),\,n+1})
\end{align*}

\ndpr

\bgP

\label{q_wd}

The linear transformation \(T\) is involutive with the composition and
therefore bijective.

From proposition \ref{perm_q_t} we deduce that
\(Q(\sigma) = T\cdot P(\sigma)\cdot T\) and that
\(\forall\sigma\in\mtt{Sym}(n+1),\,Q(\sigma)\in\mtt{GL}(\ent\lres2\ent)^{n+1}\).

\ndP

\begin{proof}

We just have to see that \(T^2 = \id\). \begin{align*}
&(T\cdot T)_{i,\,j}
=\sum_kT_{i,\,k}T_{k,\,j}=\sum_k(\delta_{i,\,k} + \delta_{k,\,n+1}(1 + \delta_{i,\,n+1}))(\delta_{k,\,j} + \delta_{j,\,n+1}(1 + \delta_{k,\,n+1}))\\
=&\delta_{i,\,j} + \delta_{j,\,n+1}(1 + \delta_{i,\,n+1})
+(1 + \delta_{i,\,n+1})(\delta_{n+1,\,j} + \delta_{j,\,n+1}(1 + \delta_{n+1,\,n+1}))=\delta_{i,\,j} = (\id)_{i,\,j}
\end{align*}

As for any \(\sigma\in\mtt{Sym}(n+1)\)
\(T,\,P(\sigma)\in\mtt{GL}(\ent\lres2\ent)^{n+1}\), then
\(Q(\sigma)\in\mtt{GL}(\ent\lres2\ent)^{n+1}\).

\ndpr

\begin{quote}
The following tells us that \(Q\) is an action, as the product of
\(\mtt{GL}(\ent\lres2\ent)^{n+1}\) is composition.
\end{quote}

\bgP

\label{q_act}

\(Q:\mtt{Sym}(n+1)\to\mtt{GL}(\ent\lres2\ent)^{n+1}\) is an injective
morphism.

\ndP

\begin{proof}

Let \(\sigma,\,\sigma'\in\mtt{Sym}(n+1)\). By using proposition
\ref{q_wd} we can deduce
\(Q(\sigma\cdot\sigma') = T\cdot P(\sigma\cdot\sigma')\cdot T = T\cdot P(\sigma)\cdot P(\sigma')\cdot T = T\cdot P(\sigma)\cdot T\cdot T\cdot P(\sigma')\cdot T = Q(\sigma)\cdot Q(\sigma')\).

Let us suppose now for some \(\sigma\in\mtt{Sym}(n+1)\) that
\(Q(\sigma) = \id\). As
\(P(\sigma) = T\cdot Q(\sigma)\cdot T = T\cdot T = \id\) and \(P\) is
injective, we deduce \(\sigma = \id\).

\ndpr

\begin{quote}
\begin{multicols}{2}

\begin{raggedright}
As $Q$ is an injective morphism into $\mtt{GL}$ we can see straightforward that it is a faithful action.

In spite of that, $Q$ is a not free action, as $\forall n\geq2,\,Q((1\ 2),\,(0,\,\ldots,\,0)) = (0,\,\ldots,\,0)$.
\end{raggedright}

\columnbreak
$Q(\sigma)$ is a permutation of a matrix like this:
$$
\begin{pmatrix}
\id_{\sigma(n+1)-1} & 1+\delta_{\sigma(n+1),\,n+1} & 0\\
0 & 1 & 0\\
0 & 1+\delta_{\sigma(n+1),\,n+1} & \id_{n+1-\sigma(n+1)}\\
\end{pmatrix}
$$
\end{multicols}
\end{quote}

\hypertarget{a-new-action-alpha-from-q}{%
\subsubsection*{\texorpdfstring{A new action \(\alpha\) from
\(Q\)}{A new action \textbackslash alpha from Q}}\label{a-new-action-alpha-from-q}}

We will redefine a third time \(\alpha_n\), using \(Q\). The action
\(R\) represents the action \(\alpha\) on tonicity signatures.

\bgD

\label{r_alpha}

We define an action on tonicity signatures: \[
  \begin{matrix}R&:&\mtt{Sym}(n+1)\times(\ent\lres2\ent)^{n+1}&\to&(\ent\lres2\ent)^{n+1}\\
  &&(\sigma,\,v)&\mapsto&\begin{pmatrix}d_1 \\\vdots\\d_{n+1}\end{pmatrix}+Q(\sigma,\,v)\end{matrix}\]
where
\(d_j = (1 + \delta_{n+1,\,\sigma(j)})(1 + \delta_{n+1,\,\sigma(n+1)})\).
Note that \(d_{\sigma^{-1}(n+1)} = 0\).

And an action on connective skeletons: \[
  \begin{matrix}\alpha''_n&:&\mtt{Sym}(n+1)\times\mbb{C}_n&\to&\mbb{C}_n\\
  &&(\sigma,\,\begin{bmatrix}\sigma'\\ k\\ v\\ s\end{bmatrix})&\mapsto&\begin{bmatrix}\sigma\cdot\sigma'\\ P_{\nat^*}(\sigma,\,k)\\ R(\sigma,\,v)\\ s\end{bmatrix}\end{matrix}\]

\ndD

\bgP

\label{alpha_act}

\(R\) and \(\alpha_n''\) are actions.

\ndP

\begin{proof}

\(R\) being an action entails \(\alpha''_n\) being an action. Indeed, by
using proposition \ref{p_act}:

{\vspace{-2.5em}\small
$$\alpha''_n(\sigma,\,\alpha''_n(\sigma',\,\begin{bmatrix}\rho\\ k\\ v\\ s\end{bmatrix})) = \alpha''_n(\sigma,\,\begin{bmatrix}\sigma'\cdot\rho\\ P_{\nat^*}(\sigma',\,k)\\ R(\sigma',\,v)\\ s\end{bmatrix}) = \begin{bmatrix}\sigma\cdot\sigma'\cdot\rho\\ P_{\nat^*}(\sigma,\,P_{\nat^*}(\sigma',\,k))\\ R(\sigma,\,R(\sigma',\,v))\\ s\end{bmatrix} = \begin{bmatrix}\sigma\cdot\sigma'\cdot\rho\\ P_{\nat^*}(\sigma\cdot\sigma',\,k)\\ R(\sigma\cdot\sigma',\,v)\\ s\end{bmatrix} = \alpha''_n(\sigma\cdot\sigma',\,\begin{bmatrix}\rho\\ k\\ v\\ s\end{bmatrix})$$}

We state the following equality, which we will prove afterwards.
\begin{equation}
\begin{pmatrix}(1 + \delta_{n+1,\,\sigma(1)})(1 + \delta_{n+1,\,\sigma(n+1)})\\\vdots\\(1 + \delta_{n+1,\,\sigma(n+1)})(1 + \delta_{n+1,\,\sigma(n+1)})\end{pmatrix} + Q(\sigma,\,b) = \begin{pmatrix}1\\\vdots\\1\end{pmatrix} + Q(\sigma,\,\begin{pmatrix}1\\\vdots\\1\end{pmatrix} + b)
\end{equation} Now we can prove the main result

{\vspace{-2.5em}\scriptsize
\begin{align*}R(\sigma,\,R(\sigma',\,v))
  &= R(\begin{pmatrix}1\\\vdots\\1\end{pmatrix}+Q(\sigma',\,\begin{pmatrix}1\\\vdots\\1\end{pmatrix} + v))
  = \begin{pmatrix}1\\\vdots\\1\end{pmatrix}+Q(\sigma,\,\begin{pmatrix}1\\\vdots\\1\end{pmatrix}+\begin{pmatrix}1\\\vdots\\1\end{pmatrix}+Q(\sigma',\,\begin{pmatrix}1\\\vdots\\1\end{pmatrix}+v))\\
  &= \begin{pmatrix}1\\\vdots\\1\end{pmatrix}+Q(\sigma,\,Q(\sigma',\,\begin{pmatrix}1\\\vdots\\1\end{pmatrix}+v))
  = \begin{pmatrix}1\\\vdots\\1\end{pmatrix}+Q(\sigma\cdot\sigma',\,\begin{pmatrix}1\\\vdots\\1\end{pmatrix}+v)
  = R(\sigma\cdot\sigma',\,v)
\end{align*}}

where we used (1) in the first two steps, we used that the additive
group \((\ent\lres2\ent)^{n+1}\) has characteristic 2 in the third step
and we used proposition \ref{q_act}.

The equality can be shown by case analysis on \(\sigma(n+1) = n+1\).

\begin{itemize}
\tightlist
\item
  If \(\sigma(n+1) = n+1\), then
  \(Q(\sigma,\,(1,\,\ldots,\,1)) = (1,\,\ldots,\,1)\) and therefore the
  right hand side is \(0 + Q(\sigma,\,b)\) like the left hand side
  (because \(1 + \delta_{\sigma(n+1),\,n+1}) = 0\)).
\item
  If \(\sigma(n+1)\neq n+1\), then
  \(Q(\sigma,\,(1,\,\ldots,\,1)) = (\delta_{n+1,\,\sigma(1)},\,\ldots,\,\delta_{n+1,\,\sigma(n+1)})\)
  for \(\sigma(n+1)\neq n+1\) and therefore we get
  \[R(\sigma,\,b) = \begin{pmatrix}1+\delta_{n+1,\,\sigma(1)}\\\ldots\\1+\delta_{n+1,\,\sigma(n+1)}\end{pmatrix} + Q(\sigma,\,b)\]
\end{itemize}

\ndpr

\bgL

\label{alpha2_free}

\(\alpha_n''\) is a free action.

\ndL

\begin{proof}

Checking it is as simple as observing that in its first component
\(\sigma\cdot\sigma' = \sigma'\) if and only if \(\sigma = \id\).

\ndpr

\hypertarget{correspondence-between-the-actions-alpha}{%
\subsubsection*{\texorpdfstring{Correspondence between the actions
\(\alpha\)}{Correspondence between the actions \textbackslash alpha}}\label{correspondence-between-the-actions-alpha}}

In this proposition we show that \(\alpha_n\) and \(\alpha''_n\) are the
same action, modulo translation between the different syntax \(\mbb{C}\)
and \(\mbb{D}\). We already knew that \(\alpha_n\) and \(\alpha'_n\)
were equal.

\bgP

\label{alpha2_alpha}

For any connective skeleton \(\ostar\) and permutation \(\sigma\),
\(\alpha_n''(\sigma,\,\ostar) = \iota(\alpha_n(\sigma,\,\iota^{-1}(\ostar)))\),
where \(\iota:\mbb{C}_n\to\mbb{D}_n\) is the translation given in
definition \ref{skel_trans}.

Along with proposition \ref{origin_alph} we can now tell that all the
presented alpha actions are equal.

\ndP

\begin{proof}

We provide the proof for
\[R(\sigma,\,\begin{pmatrix}\pm_1\\\vdots\\\pm_{n+1}\end{pmatrix}) = \begin{pmatrix}1+\pm_{\sigma(n+1)}+\pm_{(n+1\ \sigma(n+1))(\sigma(1))}\\\vdots\\1+\pm_{\sigma(n+1)}+\pm_{(n+1\ \sigma(n+1))(\sigma(i-1))}\\\pm_{\sigma(n+1)}\\1+\pm_{\sigma(n+1)}+\pm_{(n+1\ \sigma(n+1))(\sigma(i+1))}\\\vdots\\1+\pm_{\sigma(n+1)}+\pm_{(n+1\ \sigma(n+1))(\sigma(n+1))}\end{pmatrix}\]

Indeed, \begin{align*}
&(R(\sigma,\,(\pm_1,\,\ldots,\,\pm_n)))_{i}\\
=& 1 + \delta_{n+1,\,\sigma(i)}+\sum_j Q(\sigma)_{i,\,j}\pm_j\\
=& 1 + \delta_{n+1,\,\sigma(i)} + \sum_j (\delta_{(n+1\ \sigma(n+1))(\sigma(i)),\,j}
    + \delta_{\sigma(n+1),\,j}
        (1 + \delta_{n+1,\,\sigma(i)})
        (1 + \delta_{\sigma(n+1),\,n+1}))\pm_j\\
=& 1 + \delta_{n+1,\,\sigma(i)} + \pm_{(n+1\ \sigma(n+1))(\sigma(i))}
    + \pm_{\sigma(n+1)}
        (1 + \delta_{n+1,\,\sigma(i)})
        (1 + \delta_{\sigma(n+1),\,n+1})\\
\end{align*}

\begin{itemize}
\tightlist
\item
  If \(n+1 = \sigma(i)\) we have that
  \[(R(\sigma,\,(\pm_1,\,\ldots,\,\pm_n)))_{i} = \pm_{(n+1\ \sigma(n+1))(\sigma(i))} = \pm_{\sigma(n+1)}\]
\item
  If \(n+1\neq \sigma(i)\) but \(n+1 = \sigma(n+1)\) we have that
  \[(R(\sigma,\,(\pm_1,\,\ldots,\,\pm_n)))_{i} = \pm_{(n+1\ \sigma(n+1))(\sigma(i))} = \pm_{\sigma(i)}\]
\item
  In any other case we have
  \[(R(\sigma,\,(\pm_1,\,\ldots,\,\pm_n)))_{i} = 1 + \delta_{n+1,\,\sigma(i)} + \pm_{(n+1\ \sigma(n+1))(\sigma(i))}
    + \pm_{\sigma(n+1)}\]
\end{itemize}

Now we can easily check that for \(\sigma(n+1)\neq n+1\)
\[\alpha''_n(\sigma,\,\begin{bmatrix}\sigma'\\ k\\ v\\ s\end{bmatrix}) =
\begin{bmatrix}\sigma\cdot\sigma'\\ P_{\nat^*}(\sigma,\,k)\\ R(\sigma,\,v)\\ s\end{bmatrix} =
\begin{bmatrix}\sigma\cdot\sigma'\\ (k_{\sigma(1)},\,\ldots,\,k_{\sigma(n+1)})\\ \begin{pmatrix}-\pm_{\sigma(n+1)}\pm_{(n+1\ \sigma(n+1))(\sigma(1))}\\\vdots\\-\pm_{\sigma(n+1)}\pm_{(n+1\ \sigma(n+1))(\sigma(i-1))}\\\pm_{\sigma(n+1)}\\-\pm_{\sigma(n+1)}\pm_{(n+1\ \sigma(n+1))(\sigma(i+1))}\\\vdots\\-\pm_{\sigma(n+1)}\pm_{(n+1\ \sigma(n+1))(\sigma(n+1))}\end{pmatrix}\\ s\end{bmatrix}\]
Where the last component of the equality is written in multiplicative
notation of \(\mbb{B}\).

\ndpr

\bgP

\label{q_from_r}

Let
\(\ostar\in\mbb{C}_n,\,\sigma\in\mtt{Sym}(n+1),\,b,\,v\in(\ent\lres2\ent)^{n+1}\).
Then

\begin{enumerate}
\def\labelenumi{\arabic{enumi}.}
\tightlist
\item
  We have \(R(\sigma,\,b+v) = Q(\sigma,\,b)+R(\sigma,\,v)\),
\item
  We have
  \(\alpha_n''(\sigma,\,\varsigma(b,\,\ostar)) = \varsigma(Q(\sigma,\,b),\,\alpha_n''(\sigma,\,\ostar))\),
\item
  We have
  \(\forall\sigma\in\mtt{Sym}(n+1),\,b,\,v\in(\ent\lres2\ent)^{n+1},\,R(\sigma,\,b) + R(\sigma,\,v) = Q(\sigma,\,b + v)\).
\end{enumerate}

\ndP

\begin{proof}

\begin{align*}
R(\sigma,\,b+v) &= \begin{pmatrix}(1 + \delta_{n+1,\,\sigma(1)})(1 + \delta_{n+1,\,\sigma(n+1)}) \\\vdots\\(1 + \delta_{n+1,\,\sigma(n+1)})(1 + \delta_{n+1,\,\sigma(n+1)})\end{pmatrix}+Q(\sigma,\,b+v)\\
&= \begin{pmatrix}(1 + \delta_{n+1,\,\sigma(1)})(1 + \delta_{n+1,\,\sigma(n+1)}) \\\vdots\\(1 + \delta_{n+1,\,\sigma(n+1)})(1 + \delta_{n+1,\,\sigma(n+1)})\end{pmatrix}+Q(\sigma,\,b)+Q(\sigma,\,v) = Q(\sigma,\,b) + R(\sigma,\,v)
\end{align*} where in the last equality we used commutativity of \(+\).

The points 2 and 3 are easy consequences of 1.

\ndpr

\hypertarget{the-free-product-and-the-semi-direct-product}{%
\subsubsection*{The free product and the semi-direct
product}\label{the-free-product-and-the-semi-direct-product}}

In this section we will introduce two possible ways of composing the
actions on connective skeletons. The first one, with the free product,
can be understood as recursively applying the actions for each element
of a group appearing in a word. The second one is presented below.

\bgD

\label{semi_dire_Q}

Following proposition \ref{q_from_r}, we define the product action
\(\alpha_n\times\varsigma_n:(\ent\lres2\ent)^{n+1}\rtimes\mtt{Sym}(n+1)\times\mbb{C}_n\to\mbb{C}_n\)
on the semi-direct product
\((\ent\lres2\ent)^{n+1}\rtimes\mtt{Sym}(n+1)\) for the morphism
\(Q\).\footnote{See the Algebraic Preliminaries section for details.}

\ndD

\bgD

We call \(\alpha_n\) the residuation action, \(\beta_n\) the boolean
negation action and \(\alpha_n\times\varsigma_n\) the signed residuation
action. We note them without subindex usually.

\ndD

\bgP

\label{sgn_res_free}

\(\alpha\times\varsigma\) is a free action.

\ndP

\begin{proof}

If
\((\alpha\times\varsigma)(((s_1,\,\ldots,\,s_{n+1}),\,\sigma),\,\ostar) = \ostar\)
then on the first component \(\sigma\cdot\sigma' = \sigma'\) so that
\(\sigma = \id\) and
\(\varsigma((s_1,\,\ldots,\,s_{n+1}),\,\ostar) = \ostar\). We have
already seen that \(\varsigma\) is free, so that
\((s_1,\,\ldots,\,s_{n+1}) = (0,\,\ldots,\,0)\).

\ndpr

\begin{quote}
This action translated into Aucher's notation corresponds to composing
the switches \(s_i\) and the dual \(\delta\) with \(\alpha\).
\end{quote}

\bgD

We define the morphism
\(\varphi:\ent\lres2\ent*\mtt{Sym}(n+1)\to(\ent\lres2\ent)^{n+1}\rtimes\mtt{Sym}(n+1)\):
For any
\(s_1,\,\ldots,\,s_{l}\in\ent\lres2\ent,\,\sigma_1,\,\ldots,\,\sigma_l\in\mtt{Sym}(n+1)\),
\begin{align*}\varphi((s_1,\,\sigma_1)\ldots(s_l,\,\sigma_l)) = (s_1(1,\,\ldots,\,1) + \ldots + Q(\sigma_1\cdots\sigma_{l-1},\,s_l(1,\,\ldots,\,1)),\,\sigma_1\cdots\sigma_l)\end{align*}

We define the function
\(\psi:(\ent\lres2\ent)^{n+1}\rtimes\mtt{Sym}(n+1)\to\ent\lres2\ent*\mtt{Sym}(n+1)\):
For any
\(s_1,\,\ldots,\,s_{n+1}\in\ent\lres2\ent,\,\sigma\in\mtt{Sym}(n+1)\),
\begin{align*}
  &\psi(s_1,\,\ldots,\,s_{n+1},\,\sigma)\\
  &=(0,\,(1\ n+1))(s_1,\,(1\ n+1))\ldots(0,\,(n\ n+1))(s_n,\,(n\ n+1))\\
  &\phantom{=(}((s_{n+1},\,\id)(0,\,(1\ n+1))(s_{n+1},\,(1\ n+1))\ldots(0,\,(n\ n+1))(s_{n+1},\,(n\ n+1)))(+,\,\sigma)\end{align*}

\ndD

\bgP

\(\varphi\) is a morphism.

\ndP

\begin{proof}

We first see well-definedness of \(\varphi\). We have to check which are
the values of \(\varphi(v(\sigma,\,s)(\id,\,s')w)\),
\(\varphi(v(\sigma,\,0)(\sigma',\,s))\) and \(\varphi(v(\id,\,0)w)\).

\begin{itemize}
\tightlist
\item
  In the first case, if
  \(v = (\sigma_1,\,s_1)\cdots(\sigma_{i-1},\,s_{i-1})\), we have a
  component
  \(Q(\sigma_1\cdots\sigma_{i-1},\,s(1,\,\ldots,\,1)) + Q(\sigma_1\cdots\sigma_{i-1}\cdot\id,\,s'(1,\,\ldots,\,1))\)
  which can be simplied as
  \(Q(\sigma_1\cdots\sigma_{i-1},\,(s+s')(1,\,\ldots,\,1))\). The
  associativity of \(+\) and \(\id\) being the neutral element end the
  proof.
\item
  In the second case, if
  \(v = (\sigma_1,\,s_1)\cdots(\sigma_{i-1},\,s_{i-1})\), we have a
  component \(Q(\sigma_1\cdots\sigma_{i-1},\,0(1,\,\ldots,\,1)) = 0\),
  as \(0\) is the neutral of \(+\) we can conclude.
\item
  In the third case, if
  \(v = (\sigma_1,\,s_1)\cdots(\sigma_{i-1},\,s_{i-1})\), we have a
  component \(Q(\sigma_1\cdots\sigma_{i-1},\,0(1,\,\ldots,\,1)) = 0\),
  as \(0\) is the neutral of \(+\) and \(\id\) is the neutral element of
  \(\cdot\) we can conclude.
\end{itemize}

Let \(s_1,\,\ldots,\,s_l\in \ent\lres2\ent\),
\(\sigma_1,\,\ldots,\,\sigma_l,\,\sigma,\,\sigma'\in\mtt{Sym}(n+1)\). We
see that \(\varphi\) is a morphism. \begin{align*}
&\varphi((s_1,\,\sigma_1)\ldots(s_l,\,\sigma_l)(s_1',\,\sigma_1')\ldots(s_l',\,\sigma_l'))\\
=& (s_1(1,\,\ldots,\,1) + \ldots + Q(\sigma_1\cdots\sigma_{l-1},\,s_l(1,\,\ldots,\,1)) + Q(\sigma_1\cdots\sigma_{l}\cdot\sigma_1',\,s_1'(1,\,\ldots,\,1)) + \ldots\\
&+ Q(\sigma_1\cdots\sigma_{l}\cdot\sigma_1'\cdots\sigma_{l-1}',\,s_l'(1,\,\ldots,\,1)),\,\sigma_1\cdots\sigma_l\cdot\sigma_1'\cdots\sigma_l')\\
=& (s_1(1,\,\ldots,\,1) + \ldots + Q(\sigma_1\cdots\sigma_{l-1},\,s_l(1,\,\ldots,\,1)),\,\sigma_1\cdots\sigma_l)\cdot(s_1'(1,\,\ldots,\,1) + \ldots\\
&+ Q(\sigma_1'\cdots\sigma_{l-1}',\,s_l'(1,\,\ldots,\,1)),\,\sigma_1'\cdots\sigma_l')\\
=& \varphi((s_1,\,\sigma_1)\ldots(s_l,\,\sigma_l))+\varphi((s_1',\,\sigma_1')\ldots(s_l',\,\sigma_l'))
\end{align*}

\ndpr

Now we can see that the action given by the free product of \(\alpha\)
and \(\beta\) and the semi-direct product of \(\alpha\) and
\(\varsigma\) are somehow equivalent on connective skeletons.

\bgT

\label{quot_sign}

Let \(C\) be a connective family and \(x\in\mbb{C}\). The group
\((\ent\lres2\ent)^{n+1}\rtimes\mtt{Sym}(n+1)\) defined in Definition
\ref{semi_dire_Q} is isomorphic to
\((\ent\lres2\ent*\mtt{Sym}(n+1))\lres(\ent\lres2\ent*\mtt{Sym}(n+1))_x\).
In particular we have that \(\varphi\) induces an isomorphism
\(\overline\varphi:(\ent\lres2\ent*\mtt{Sym}(n+1))\lres(\ent\lres2\ent*\mtt{Sym}(n+1))_x\simeq (\ent\lres2\ent)^{n+1}\rtimes\mtt{Sym}(n+1)\).

Furthermore,
\(\forall g\in\ent\lres2\ent*\mtt{Sym}(n+1),\,(\alpha''\times\varsigma)(\varphi(g),\,x) = (\alpha*\beta)(g,\,x)\).

\ndT

\begin{proof}

The first thing we check is the equality for the action
\(\alpha\times\varsigma\) on Connective Skeletons. We prove it by
induction on the length \(l\) of
\(g = (s_1,\,\sigma_1)\cdots(s_l,\,\sigma_l)\).

Let
\(s_1,\,\ldots,\,s_l\in\ent\lres2\ent,\,\sigma_1,\,\ldots,\,\sigma_l\in\mtt{Sym}(n+1)\)
and \(\ostar\in\mbb{C}_n\).
\((\alpha''\times\varsigma)(\varphi((s_1,\,\sigma_1)\cdots (s_l,\,\sigma_l)),\,\ostar) = (\alpha*\beta)((s_1,\,\sigma_1)\cdots (s_l,\,\sigma_l),\,\ostar)\),
we recall
\(\alpha''(\sigma,\,s + \ostar) = Q(\sigma,\,s) + \alpha''(\sigma,\,\ostar)\).

We show the induction step

\vspace{-2em}
{\small\begin{align*}&(\alpha''\times\varsigma)(\varphi((s_1,\,\sigma_1)\cdots (s_{l+1},\,\sigma_{l+1})),\,\ostar)\\
=& (\alpha''\times\varsigma)((s_1(1,\,\ldots,\,1) + \ldots + Q(\sigma_1\cdots\sigma_{l},\,s_{l+1}(1,\,\ldots,\,1)),\,\sigma_1\cdots\sigma_{l+1}),\, \ostar)\\
=& (\alpha''\times\varsigma)((s_1(1,\,\ldots,\,1) + \ldots + Q(\sigma_1\cdots\sigma_{l-1},\,s_{l}(1,\,\ldots,\,1)),\,\sigma_1\cdots\sigma_{l}),\,(\alpha\times\varsigma)((s_{l+1}(1,\,\ldots,\,1),\,\sigma_{l+1}),\,\ostar))\\
=& (\alpha*\beta)((s_1,\,\sigma_1)\cdots(s_{l},\,\sigma_{l}),\,(\alpha''\times\varsigma)((s_{l+1}(1,\,\ldots,\,1),\,\sigma_{l+1}),\,\ostar))\\
=& (\alpha*\beta)((s_1,\,\sigma_1)\cdots(s_{l},\,\sigma_{l}),\,(\alpha''*\beta)((s_{l+1},\,\sigma_{l+1}),\,\ostar))\\
=& (\alpha*\beta)((s_1,\,\sigma_1)\cdots(s_l,\,\sigma_l) (s_{l+1},\,\sigma_{l+1}),\,\ostar)\end{align*}}

where in the second step we have used the fact that
\(\alpha\times\varsigma\) is an action of the semi-direct product by
\(Q\), so that
\((\alpha''\times\varsigma)((v+Q(\sigma,\,v'),\,\sigma\cdot\sigma'),\,\ostar) = (\alpha''\times\varsigma)((v,\,\sigma)\cdot(v',\,\sigma'),\,\ostar)\),
and in the last step we have used that \(\alpha*\beta\) is an action of
the free product. In the third and fourth steps we have used induction
hypothesis.

To see that \(\overline\varphi\) is an isomorphism of groups, we will
use the first isomorphism theorem and proposition \ref{stab_quot} on
lemma \ref{alpha2_free}.

Therefore, we need to see that
\(\mtt{Ker}\ \varphi = (\ent\lres2\ent*\mtt{Sym}(n+1))_x\). Indeed
\((s_1(1,\,\ldots,\,1) + \ldots + Q(\sigma_1\cdots\sigma_{l-1},\,s_l(1,\,\ldots,\,1)),\,\sigma_1\cdot\ldots\cdot\sigma_l) = ((0,\,\ldots,\,0),\,\id)\)
implies that
\((\alpha*\beta)((s_1,\,\sigma_1)\cdots (s_{l},\,\sigma_{l}),\,\ostar) = (\alpha''\times\varsigma)(\varphi((s_1,\,\sigma_1)\cdots (s_{l},\,\sigma_{l})),\,\ostar) = (\alpha''\times\varsigma)(((0,\,\ldots,\,0),\,\id),\,\ostar) = \ostar\)
so that
\(\mtt{Ker}\ \varphi\subseteq(\ent\lres2\ent*\mtt{Sym}(n+1))_x\).

As we also have that
\((\alpha*\beta)((s_1,\,\sigma_1)\cdots (s_{l},\,\sigma_{l}),\,\ostar) = \ostar\)
implies
\((\alpha''\times\varsigma)(\varphi((s_1,\,\sigma_1)\cdots (s_{l},\,\sigma_{l})),\,\ostar) = \ostar\),
we can deduce that
\(\varphi((s_1,\,\sigma_1)\cdots(s_l,\,\sigma_l),\,\ostar) = \id\)
because \(\alpha''\times\varsigma\) is a free action, as shown in
proposition \ref{sgn_res_free}. Therefore
\(\mtt{Ker}\ \varphi\supseteq(\ent\lres2\ent*\mtt{Sym}(n+1))_x\).

Now we now know that
\(\overline\varphi:(\ent\lres2\ent*\mtt{Sym}(n+1))\lres(\ent\lres2\ent*\mtt{Sym}(n+1))_x\simeq\mtt{Im}\ \overline\varphi\).

We finish by showing that \(\psi\) is a right-inverse to \(\varphi\),
and therefore \(\varphi\) is surjective, so that \(\overline\varphi\) is
also surjective. Let
\(s_1,\,\ldots,\,s_{n+1}\in\ent\lres2\ent,\,\sigma\in\mtt{Sym}(n+1)\).
We want to prove
\(\varphi(\psi((s_1,\,\ldots,\,s_{n+1}),\,\sigma)) = ((s_1,\,\ldots,\,s_{n+1}),\,\sigma)\).
For it we unfold the definitions: \begin{align*}
&\varphi(\psi(s_1,\,\ldots,\,s_{n+1},\,\sigma)=\varphi((0,\,(1\ n+1))(s_1,\,(1\ n+1))\ldots(0,\,(n\ n+1))(s_n,\,(n\ n+1))\\
&\phantom{=(}((s_{n+1},\,\id)(0,\,(1\ n+1))(s_{n+1},\,(1\ n+1))\ldots(0,\,(n\ n+1))(s_{n+1},\,(n\ n+1)))(+,\,\sigma))\\
&= (0(1,\,\ldots,\,1) + Q((1\ n+1),\,s_1(1,\,\ldots,\,1)) + Q((1\ n+1)\cdot(1\ n+1),\,0(1,\,\ldots,\,1)) + \ldots\\
&\phantom{=(}+ Q((1\ n+1)\cdot(1\ n+1)\cdots(n\ n+1),\,s_n(1,\,\ldots,\,1))\\
&\phantom{=(}+ Q((1\ n+1)\cdot(1\ n+1)\cdots(n\ n+1)\cdot(n\ n+1),\,s_{n+1}(1,\,\ldots,\,1))\\
&\phantom{=(}+ 0(1,\,\ldots,\,1) + Q((1\ n+1),\,s_{n+1}(1,\,\ldots,\,1)) + Q((1\ n+1)\cdot(1\ n+1),\,0(1,\,\ldots,\,1)) + \ldots\\
&\phantom{=(}+ Q((1\ n+1)\cdot(1\ n+1)\cdots(n\ n+1),\,s_{n+1}(1,\,\ldots,\,1)),\,\\
&\phantom{=((}(1\ n+1)\cdot(1\ n+1)\cdots(n\ n+1)\cdot(n\ n+1)\cdot(1\ n+1)\cdot(1\ n+1)\cdots(n\ n+1)\cdot(n\ n+1)\cdot\sigma)
\end{align*} Next we use multiple times the fact that
\(\forall i\in\{1,\,\ldots,\,n+1\},\,Q((i\ n+1),\,(1,\,\ldots,\,1)) = (0,\,\ldots,\,0,\,1,\,0,\,\ldots,\,0)\),
where the 1 is located at component \(i\). \begin{align*}
&\varphi(\psi(s_1,\,\ldots,\,s_{n+1},\,\sigma)\\
&= ((0,\,\ldots,\,0) + (s_1,\,0,\,\ldots,\,0) + (0,\,\ldots,\,0) + \ldots + (0,\,\ldots,\,0,\,s_n,\,0) + s_{n+1}(1,\,\ldots,\,1)\\
&\phantom{=(}+ (0,\,\ldots,\,0) + (s_{n+1},\,0,\,\ldots,\,0) + (0,\,\ldots,\,0)) + \ldots + (0,\,\ldots,\,s_{n+1},\,0)),\,\sigma)\\
&= ((s_1,\,\ldots,\,s_n,\,0) + (0,\,\ldots,\,s_{n+1}),\,\sigma)\\
&= ((s_1,\,\ldots,\,s_n,\,s_{n+1}),\,\sigma)
\end{align*}

With it we have finished the proof of the Theorem.

\ndpr

\begin{quote}
When quotienting by the stabilizer of \(x\) we are equating all elements
in the group with equal images on \(x\), therefore providing a free
action on a smaller group. In our case the resulting group is finite and
lets us work on Atomic Logics without using the set of words on a
cartesian product (which the free product is defined on).
\end{quote}

\bgE

\label{lambek_ex1}

We now see how would the procedure used in theorem \ref{quot_sign}
transform a sequence of negations and residuations on a connective
skeleton in practice by checking the following assertion:

Let \(\rightarrow\) have skeleton
\(((2\ 3),\,-,\,\forall,\,(1,\,1,\,1),\,(-,\,+))\), then
\({-(1\ 3){-}(1\ 2\ 3)\rightarrow}\) has skeleton
\((\id,\,+,\,\exists,\,(1,\,1,\,1),\,(+,\,-))\), using the original
notation.

\ndE

\begin{proof}

We want to see \[
\iota^{-1}(\beta(1,\,\alpha''((1\ 3),\,\beta(1,\,\alpha''((1\ 2\ 3),\,\iota(\begin{bmatrix}(2\ 3)\\-\\\forall\\(1,\,1,\,1)\\(-,\,+)\end{bmatrix}))))))=\begin{bmatrix}\id\\+\\\exists\\(1,\,1,\,1)\\(+,\,-)\end{bmatrix}
\]

By definition of \(\alpha*\beta\) we have \begin{align*}
&\iota^{-1}(\beta(1,\,\alpha''((1\ 3),\,\beta(1,\,\alpha''((1\ 2\ 3),\,\iota(\begin{bmatrix}(2\ 3)\\-\\\forall\\(1,\,1,\,1)\\(-,\,+)\end{bmatrix}))))))=&\iota^{-1}((\alpha''*\beta)((1,\,(1\ 3))(1,\,(1\ 2\ 3)),\,\begin{bmatrix}(2\ 3)\\(1,\,1,\,1)\\(1,\,0,\,0)\\-\end{bmatrix}))\end{align*}
We now use \(\varphi\) and theorem \ref{quot_sign} in the first step to
get the following \begin{align*}
&\iota^{-1}((\alpha''*\beta)((1,\,(1\ 3))(1,\,(1\ 2\ 3)),\,\begin{bmatrix}(2\ 3)\\(1,\,1,\,1)\\(1,\,0,\,0)\\-\end{bmatrix}))
=(\alpha''\times\varsigma)(\varphi((1,\,(1\ 3))(1,\,(1\ 2\ 3))),\,\begin{bmatrix}(2\ 3)\\(1,\,1,\,1)\\(1,\,0,\,0)\\-\end{bmatrix})\\
=&(\alpha''\times\varsigma)(((1,\,1,\,1) + Q((1\ 3),\,(1,\,1,\,1)),\,(1\ 3)\cdot(1\ 2\ 3)),\,\begin{bmatrix}(2\ 3)\\(1,\,1,\,1)\\(1,\,0,\,0)\\-\end{bmatrix})\\
=&\iota^{-1}((\alpha''\times\varsigma)(((1,\,1,\,1) + (1,\,0,\,0),\,(2\ 3)),\,\begin{bmatrix}(2\ 3)\\(1,\,1,\,1)\\(1,\,0,\,0)\\-\end{bmatrix}))=\iota^{-1}((\alpha''\times\varsigma)(((0,\,1,\,1),\,(2\ 3)),\,\begin{bmatrix}(2\ 3)\\(1,\,1,\,1)\\(1,\,0,\,0)\\-\end{bmatrix}))\\
=&\iota^{-1}(\varsigma((0,\,1,\,1),\,\alpha''((2\ 3),\,\begin{bmatrix}(2\ 3)\\(1,\,1,\,1)\\(1,\,0,\,0)\\-\end{bmatrix})))
=\iota^{-1}(\varsigma((0,\,1,\,1),\,\begin{bmatrix}(2\ 3)\cdot(2\ 3)\\P_{\nat^*}((2\ 3),\,(1,\,1,\,1))\\R((2\ 3),\,(1,\,0,\,0))\\-\end{bmatrix}))
\end{align*}

We exemplify how would the component on the first row and second column
for \(R((2\ 3),\,(1,\,0,\,0))\) be computed:
\[\delta_{(3\ (2\ 3)(3))((2\ 3)(1)),\,2}+ \delta_{(2\ 3)(3),\,2}(1 + \delta_{3,\,(2\ 3)(1)})(1 + \delta_{(2\ 3)(3),\,3}) = \delta_{1,\,2} + \delta_{2,\,2}(1 + \delta_{3,\,1})(1 + \delta_{2,\,3}) = 1\]
Now that we have written down the linear transformations we will have to
operate, we will use the matricial expression of
\(R((2\ 3),\,(1,\,0,\,0))\). And we can conclude, \begin{align*}
&\iota^{-1}(\beta(1,\,\alpha''((1\ 3),\,\beta(1,\,\alpha''((1\ 2\ 3),\,\iota(\begin{bmatrix}(2\ 3)\\-\\\forall\\(1,\,1,\,1)\\(-,\,+)\end{bmatrix}))))))
=\iota^{-1}(\varsigma((0,\,1,\,1),\,\begin{bmatrix}(2\ 3)\cdot(2\ 3)\\P_{\nat^*}((2\ 3),\,(1,\,1,\,1))\\R((2\ 3),\,(1,\,0,\,0))\\-\end{bmatrix}))\\
=&\iota^{-1}(\varsigma((0,\,1,\,1),\,\begin{bmatrix}\id\\(1,\,1,\,1)\\\begin{pmatrix}1\\0\\1\end{pmatrix} + Q((2\ 3))\begin{pmatrix}1\\0\\0\end{pmatrix}\\-\end{bmatrix}))
=\iota^{-1}(\varsigma((0,\,1,\,1),\,\begin{bmatrix}\id\\(1,\,1,\,1)\\\begin{pmatrix}1\\0\\1\end{pmatrix} + \begin{pmatrix}1&1&0\\0&1&0\\0&1&1\end{pmatrix}\begin{pmatrix}1\\0\\0\end{pmatrix}\\-\end{bmatrix}))\\
=&\iota^{-1}(\varsigma((0,\,1,\,1),\,\begin{bmatrix}\id\\(1,\,1,\,1)\\(1+1+0,\,0+0+0,\,1+0+0)\\-\end{bmatrix}))=\iota^{-1}(\begin{bmatrix}\id\\(1,\,1,\,1)\\(0+0,\,1+0,\,1+1)\\-\end{bmatrix})=\begin{bmatrix}\id\\+\\\exists\\(1,\,1,\,1)\\(+,\,-)\end{bmatrix}
\end{align*} which ends the proof.

\ndpr

\hypertarget{group-actions-on-atomic-connectives}{%
\section{Group Actions on Atomic
Connectives}\label{group-actions-on-atomic-connectives}}

\label{sect4}

\hypertarget{connective-families}{%
\subsection{Connective Families}\label{connective-families}}

We now introduce some algebraic objects in order to be able to work and
compare different approaches on Atomic Logics. We will have to redefine
the presented Atomic Logics.

\bgD

\label{conn_fam}

A \textit{\textbf{connective}} is a symbol to which we associate a
single connective skeleton. We call \(\mtt{sk}(\ostar)\)
\textit{\textbf{the skeleton of a connective}} \(\ostar\). By abuse of
notation we will use the connective symbol to refer to its skeleton,
whenever it is not confusing.

A \textit{\textbf{family of connectives}} or connective family is

\begin{itemize}
\tightlist
\item
  a set of connectives \(C\),
\item
  a partition of it \(\mca{A}\), such that
  \(\forall D\in\mca{A},\,\forall \ostar,\,\ostar'\in D,\,\mtt{sk}(\ostar) = \mtt{sk}(\ostar')\)
  implies \(\ostar = \ostar'\),
\item
  for each \(D\in\mca{A}\) an associated group \(G_D\) and an action
  \(a_D:G_D\times\mbb{C}\to\mbb{C}\) such that
  \(\forall \ostar,\,\ostar'\in a_D(\mtt{sk}(D)),\,\ostar\) and
  \(\ostar'\) have equal arity \(n\) and have permutations
  \(\sigma,\,\sigma'\in\mtt{Sym}(n+1)\), tonicity signatures
  \(\pm,\,\pm'\in(\ent\lres2\ent)^{n+1}\) and type tuples
  \(k,\,k'\in(\nat^*)^{n+1}\) satisfying
  \(P_{\nat^*}(\sigma'\cdot\sigma^{-1},\,k) = k'\).
\end{itemize}

Whenever \(\ostar\in D\) we note \(\overline\ostar = D\).

\ndD

\begin{quote}
A simple consequence of this definition is that \(\ostar = \ostar'\) if
and only \(\overline{\ostar} = \overline{\ostar'}\) and
\(\mtt{sk}(\ostar) = \mtt{sk}(\ostar')\). It follows from injectivity of
\(\mtt{sk}\) over the classes of the partition.
\end{quote}

\bgE

\label{lamb_fam}

We define the \textit{\textbf{Lambek Connective Family}} \(L\) as:
\footnote{For more details on this example, see section \ref{sect5}.}

\begin{itemize}
\tightlist
\item
  the set \(\{\otimes,\,\rres,\,\lres\}\sqcup V\), where:

  \begin{itemize}
  \tightlist
  \item
    the tensor
    \(\mtt{sk}(\otimes)=(\id,\,(1,\,1,\,1),\,(+,\,+,\,+),\,-)\),
  \item
    the right residual
    \(\mtt{sk}(\rres)=((2\ 3),\,(1,\,1,\,1),\,(-,\,+,\,-),\,-)\),
  \item
    the left residual
    \(\mtt{sk}(\lres)=((1\ 3),\,(1,\,1,\,1),\,(+,\,-,\,-),\,-)\),
  \item
    the variables in \(V\) such that they are all propositional letters.
    There are at least a \(p\) and a \(q\) in \(V\) with
    \(\mtt{sk}(p) = (\id,\,s,\,1,\,s)\) and
    \(\mtt{sk}(q) = (\id,\,-s',\,1,\,s')\), for some
    \(s,\,s'\in\mbb{B}\).
  \end{itemize}
\item
  the partition
  \(\{\{\otimes,\,\rres,\,\lres\}\}\sqcup\{\{v\}\mid V\}\),
\item
  the action \(\alpha\times\varsigma\) on
  \(\{\otimes,\,\rres,\,\lres\}\) and the action \(\id\) on each
  \(\{v\}\), for \(v\in V\).
\end{itemize}

\ndE

\bgD

\label{atom_form}

Given a family of connectives \(C\), the
\textit{\textbf{atomic language}} of \(C\) is defined as the smallest
set \(\mca{L}_C\) such that if \(\ostar\) is a connective of arity n,
type signature \((k_1,\,\ldots,\,k_n,\,k)\) and
\(\varphi_1,\,\ldots,\,\varphi_n\in\mca{L}_C\) are of types
\(k_1,\,\ldots,\,k_n\), then \(\ostar(\varphi_1,\,\ldots,\,\varphi_n)\)
is also in \(\mca{L}_C\) and of type \(k\).\\
Its elements are called \textit{\textbf{formulas}}.

\ndD

\begin{figure}\begin{center}

\includegraphics{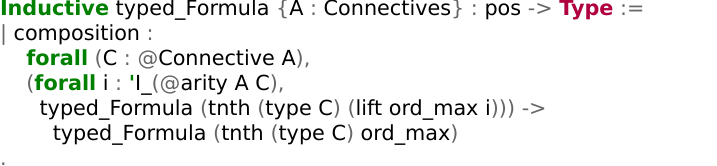}

\caption{Implementation of formulas for $\protect{(\alpha*\beta)(\protect\ostar)}$ connective families.}

\end{center}\end{figure}

\bgD

\label{equiv_fam}

Over the families of connectives we define the relation \(C \lesssim D\)
if and only if there is an injection \(\phi:C\to D\) such that for each
\(\ostar\in C\) there is a morphism
\(\varphi:G_{\overline{\phi(\ostar)}}\to G_{\overline{\ostar}}\) and
\(\forall g\in G_{\overline{\phi(\ostar)}},\,a_{\overline{\ostar}}(\varphi(g),\,\mtt{sk}(\ostar)) = a_{\overline{\phi(\ostar)}}(g,\,\mtt{sk}(\phi(\ostar)))\)
and such that
\(\forall\ostar,\,\ostar'\in C,\,\overline{\ostar} = \overline{\ostar'}\)
implies \(\overline{\phi(\ostar)} = \overline{\phi(\ostar')}\). We say
then that \(C\) is \textit{\textbf{compatible with}} the action of \(D\)
or that \(D\) extends \(C\).

Over the families of connectives we define the relation \(C \simeq D\)
if and only if there is a bijection \(\phi:C\to D\) such that for each
\(\ostar\in C\) there is a morphism
\(\varphi:G_{\overline{\ostar}}\to G_{\overline{\phi(\ostar)}}\),
\(\forall g\in G,\,a_{\overline{\ostar}}(g,\,\mtt{sk}(\ostar)) = a_{\overline{\phi(\ostar)}}(\varphi(g),\,\mtt{sk}(\phi(\ostar)))\)
and for each \(\ostar\in D\) there is a morphism
\(\psi:G_{\overline{\phi(\ostar)}}\to G_{\overline{\ostar}}\),
\(\forall g\in G,\,a_{\overline{\ostar}}(\varphi(g),\,\mtt{sk}(\ostar)) = a_{\overline{\phi(\ostar)}}(g,\,\mtt{sk}(\phi(\ostar)))\)
and such that
\(\forall\ostar,\,\ostar'\in C,\,\overline{\ostar} = \overline{\ostar'}\)
if and only if \(\overline{\phi(\ostar)} = \overline{\phi(\ostar')}\).

Over the families of connectives we define the relation \(C \sim_W' D\)
if and only if there is a bijection \(\phi:C\to D\) such that
\(\forall\ostar,\,\ostar'\in C,\,\overline{\ostar} = \overline{\ostar'}\)
implies \(\overline{\phi(\ostar)} = \overline{\phi(\ostar')}\) and for
each \(\ostar\in C\) there is a morphism
\(\varphi:G_{\overline{\ostar}}\to G_{\overline{\phi(\ostar)}}\),
\(\forall g\in G,\,a_{\overline{\ostar}}(g,\,\mtt{sk}(\ostar)) = a_{\overline{\phi(\ostar)}}(\varphi(g),\,\mtt{sk}(\phi(\ostar)))\)
and the induced morphism
\(\overline\varphi:G_{\overline{\ostar}}\lres(G_{\overline{\ostar}})_{\ostar}\to G_{\overline{\phi(\ostar)}}\lres(G_{\overline{\phi(\ostar)}})_{\phi(\ostar)}\)
is a group isomorphism. We define \(C\sim_W D\) as the
transitive-symmetric closure of \(\sim_W'\).

\ndD

\begin{quote}
For two connective families \(C\), \(D\), the \(C\lesssim D\) relation
will be used to tell that we can somehow use \(D\)'s action on \(C\)
also. The precise way we can do it will be encoded by the morphisms
\(\varphi\).
\end{quote}

\bgP

The relation \(\lesssim\) is transitive and reflexive. The relations
\(\simeq\) and \(\sim_W\) are equivalence relations.

\ndP

\begin{proof}

For the reflexivity we can take the identity as both the injection and
the morphisms for all the relations.

For the transitivity we must take the composition of both the injection
and the morphisms in both \(\lesssim\) and \(\simeq\). In \(\sim_W\) we
have it by definition.

For the symmetry of \(\simeq\) we must take the inverse of both the
bijection and the morphisms. For \(\sim_W\) we have it by definition.

Checking the conditions is routine.

\ndpr

\bgP

\label{ext_ext_sim}

Let \(C\) and \(D\) be two connective families.

\begin{itemize}
\tightlist
\item
  \(C\lesssim D\) and \(D\lesssim C\) imply \(C\simeq D\).
\item
  \(C\simeq D\) implies \(C\sim_W D\).
\end{itemize}

\ndP

\begin{proof}

We will present the proof of Schröder-Bernstein Theorem and reason over
the constructed bijection for checking the extra conditions.

\setcounter{equation}{0}

Let \(\phi:C\to D\) and \(\psi:D\to C\) be the injective functions and
\(\varphi:G_{C,\,\overline{\varphi(\,.\,)}}\to G_{D,\,\overline{\,\cdot\,}}\)
and
\(\varphi':G_{D,\,\overline{\,\cdot\,}}\to G_{C,\,\overline{\varphi(\,.\,)}}\)
be the morphisms on groups given by \(\lesssim\). We let \(\phi^-\) and
\(\psi^-\) be the partial right-inverses of \(\phi\) and \(\psi\). Then
for any \(x\in C\) we can define two sequences of images and pre-images
by
\begin{equation}\cdots\to\phi^-(\psi^-(x))\to\psi^-(x)\to x\to\phi(x)\to\psi(\phi(x))\to\cdots\end{equation}

The set of elements in \(C\) appearing in sequence 1 form a partition of
\(C\). Similarly for \(D\).

Each sequence like 1 can end in an element of \(C\) or an element of
\(D\), repeat itself forever or remain forever changing on the left.

We define a new bijection by case analysis on the previous four cases on
the partitions given by 1. Indeed, if two different elements are found
in the same sequence, both their sequences will be exactly the same as
1.

For any element of \(C\) in a sequence ending on \(C\) or never ending,
the injection \(\phi\) is a bijection onto the elements of \(D\) in the
same sequence (because for any element \(y\in D\) of the sequence 1 has
a pre-image \(\phi^-(y)\) by definition). For any element of \(C\) in a
sequence ending on \(D\), the partial right-inverse \(\psi^-\) is
well-defined and therefore it is a bijection onto the the elements of
\(D\) in the same sequence. This will be the desired bijection \(\rho\).

Now we can check on the partition given by the sequences that the new
bijection still satisfies the conditions on the actions. If
\(\ostar\in C\) we will use the equality from \(\phi\) or \(\psi^-\), if
\(\ostar\in D\) we will use the equality from \(\psi\) or \(\phi^-\).
Whenever \(\ostar\) is in a sequence ending on \(D\) we will have to
work with \(\psi^-\) and \(\psi\) as the desired inverses, which are
well-defined on \(\ostar\). For any other case we will have to work with
\(\phi\) and \(\phi^-\) as desired inverses, which are well-defined on
\(\ostar\).

The last point comes from lemma \ref{morph_stab}.

\ndpr

\begin{quote}
If the actions are free, \(C\simeq D\) is equivalent to the existence of
a bijection \(\phi:C\to D\) such that
\(\forall\ostar,\,\ostar'\in C,\,\overline\ostar = \overline{\ostar'}\)
if and only if \(\overline{\phi(\ostar)} = \overline{\phi(\ostar')}\)
and \(G_C\simeq G_{\phi(C)}\) thanks to proposition \ref{stab_quot}.
\end{quote}

\bgD

The connectives of \(n\)-ary skeletons are called
\(n\)-\textit{\textbf{ary connectives}}. The 0-ary connectives are
called \textit{\textbf{propositional letters}}.

Propositional letters whose assigned action is the identity are called
\textit{\textbf{variables}}.

\ndD

\hypertarget{families-of-structural-connectives}{%
\subsection{Families of Structural
Connectives}\label{families-of-structural-connectives}}

\bgD

Let \(C\) be a family of connectives. A family of structural
\(C\)-connectives, or simply \textit{\textbf{structural connectives}},
\([C]\) is composed of:

\begin{itemize}
\tightlist
\item
  a set \(E\),
\item
  an injective function \([\,.\,]:C\to E\),
\item
  a partial function \(p:E\to C\), right inverse of \([\,.\,]\),
\item
  an injective function \([\,.\,]:\mca{A}\to\mca{P}(E)\) such that
  \([\mca{A}]\) is a partition of \(E\), each \(D\in\mca{A}\) has an
  associated action \([a_D]:G_D\times [D]\to [D]\) such that
  \([D] = [a_D](\{[\ostar]\mid\ostar\in D\})\),
\item
  an application \(\mtt{sk}:E\to\mbb{C}\) such that
  \(\mtt{sk}([\ostar]) = \mtt{sk}(\ostar)\),
  \(\mtt{sk}([a_\ostar](g,\,[\ostar])) = a(g,\,\mtt{sk}(\ostar))\).
\end{itemize}

\ndD

\begin{quote}
Note that a family of structural connectives is also by itself a
connective family.
\end{quote}

\bgP

Let \(C\) be a connective family. There is at least one family of
structural \(C\)-connectives. Any families of structural
\(C\)-connectives \([C]\), \([C]'\) satisfy \([C]\simeq[C]'\).

\ndP

\begin{proof}

We define a family of structural \(C\)-connectives:

\begin{itemize}
\item
  The set \(E\) is
  \(\{(\overline{g},\,a_D(\mtt{sk}(\ostar)),\,D)\mid D\in\mca{A},\,a_D(\mtt{sk}(\ostar))\in\{a_D(\mtt{sk}(\ostar))\mid\ostar\in D\},\,\overline{g}\in G_D\lres (G_D)_{\mtt{sk}(\ostar)}\}\).
\item
  As each \(D\in\mca{A}\) is of fixed arity and \(\mbb{C}_n\) is a
  finite set, we can give an order \(\leq\) to \(\mbb{C}_n\). We will
  call \(\mtt{min}_\leq(J)\) the minimal connective in the set \(J\)
  from the chosen order \(\leq\).

  The injection is then
  \([\,.\,]:\ostar\mapsto(\overline{g},\,a_D(\mtt{sk}(\ostar)),\,\overline\ostar)\),
  where \(\overline{g}\in G_D\lres (G_D)_{\mtt{sk}(\ostar)}\) is such
  that
  \(a_D(g,\,\mtt{min}_\leq(\mtt{sk}(a_D(\mtt{sk}(\ostar))))) = \mtt{sk}(\ostar)\).
\item
  For every \((g,\,O,\,D)\in E\), if it exists, we note the only
  connective of skeleton \(a_D(g,\,\mtt{min}_\leq(O))\) as \(p(g,\,O)\).

  This \(p\) will be the desired partial right inverse.
\item
  The partition \([\mca{A}]\) is
  \(\{\{(\overline{g},\,a_D(\mtt{sk}(\ostar)),\,D)\mid a_D(\mtt{sk}(\ostar))\in\{a_D(\mtt{sk}(\ostar))\mid\ostar\in D\},\,\overline{g}\in G_D\lres (G_D)_{\mtt{sk}(\ostar)}\}\mid D\in\mca{A}\}\).

  We can see it is in bijection with \(\mca{A}\). We can use any of them
  interchangeably whenever they appear as an index.
\item
  The actions are
  \([a_D]:(g,\,(\overline{h},\,O,\,D))\mapsto(\overline{g\cdot_{G_D}h},\,O,\,D)\)
  on each \(D\in[\mca{A}]\).
\item
  The skeleton is
  \(\forall D\in[\mca{A}],\,(g,\,O)\mapsto a_D(g,\,\mtt{min}_\leq(O))\).
\end{itemize}

Checking all the conditions is routinary.

Now we want to see \([C]\simeq[C]'\) for any pair of families of
structural \(C\)-connectives.

\([\mca{A}]\), \([\mca{A}]'\) and \(\mca{A}\) are in bijection by
definition, as the latter is the image of the former and \([\,.\,]\) is
an injective function. Furthermore,
\(\forall[D]\in[\mca{A}],\,[D] = [a_D](\{[\ostar]\mid\ostar\in D\})\)
and
\(\forall[D]'\in[\mca{A}]',\,[D]' = [a_D]'(\{[\ostar]\mid\ostar\in D\})\).
As they are connective families, inside an orbit the skeleton is enough
to characterize each connective. Therefore,
\(\mtt{sk}([a_D](g,\,[\ostar])) = a(g,\,\mtt{sk}(\ostar)) = \mtt{sk}([a_D]'(g,\,[\ostar]))\)
tells us that there are bijections between each \([D]\) and \([D]'\).
This three results let us see that there is a bijection
\(\phi:[C]\to[C]'\).

Using the restriction on the skeleton again we can deduce that the
bijection \(\phi\) and the identity morphisms \(\id:G_D\to G_D\) satisfy
the conditions for \([C]\lesssim[C]'\). Proposition \ref{ext_ext_sim}
lets us deduce \([C]\simeq [C]'\).

\ndpr

\begin{quote}
This means that essentially there is only one possible family of
structural connectives \([C]\) for any connective family \(C\).

When families of connectives have free actions defined on them we can
use proposition \ref{lem_reg_act} to show that just by choosing an
arbitrary connective as identity we can give its orbit a group
structure.

By abuse of notation we will also write \([C]\) to refer to the set of a
family of structural \(C\)-connectives, although it doesn't correspond
to \(\{[\ostar]\mid\ostar\in C\}\). Note also that
\(\{[\ostar]\mid\ostar\in D\}\neq[D]\)!
\end{quote}

Over families of structural connectives we'll build our Sequent Calculi.
For this we introduce structures:

\bgD

\label{comp_beta}

Let \(C\) be a connective family. We define \(C\cap\beta(C)\) as the
connective family given by the partition
\(\{\{\ostar'\mid\mtt{sk}(\ostar') = \beta(s,\,\mtt{sk}(\ostar)),\,\ostar'\in D,\,s\in\mbb{B}\}\mid \ostar\in D,\,D\in\mca{A}\}\),
the group \(\mbb{B}\), the action \(\beta\) and the same set \(C\).

We say that the connective family \(C\) is compatible with \(\beta\) if
\(C\) is compatible with the action of \(C\cap\beta(C)\).

\ndD

\begin{quote}
Boolean negation plays an important role on Atomic Logics. We have to
give it a special place in the construction of structural language to
ensure it can be applied to any connective.
\end{quote}

\bgD

\label{atom_struct}

Let \(C\) be a connective family. Each formula in \(\mca{L}_C\) is a
structure, for each formula \(\varphi\) we add a new structure
\(*\varphi\) and for each structural \(n\)-ary connective \([\ostar]\)
of input type \((k_1,\,\ldots,\,k_n)\) and
\(\forall i\in \{1,\,\ldots,\,n\}\) a structure \(X_i\) of output type
\(k_i\) we have the new structure \([\ostar](X_1,\,\ldots,\,X_n)\) and,
if \(\ostar\) is not compatible with \(\beta\), the structure
\(\ast[\ostar](X_1,\,\ldots,\,X_n)\). We call structures resulting from
the two latter cases strict structures and the former two formulas.

Structures form the set \([\mca{L}_C]\).

\ndD

\begin{quote}
Note that structures are not the formulas on structural connectives, but
rather the closure of formulas with structural connectives.
\end{quote}

\bgD

We define by cases \textit{\textbf{Boolean negation on structures}} of a
family of connectives \(C\) compatible with \(\beta\) as \(*X\) through
a morphism \(\psi\):

\begin{itemize}
\tightlist
\item
  For any \(n\)-ary connective \(\ostar\in C\) compatible with \(\beta\)
  and structures \(X_i\) we define
  \(\ast[\ostar](X_1,\,\ldots,\,X_n) = [a_D](\psi(-),\,[\ostar])(X_1,\,\ldots,\,X_n)\).
\item
  For any formula \(\varphi\), \(\ast\varphi = \ast\varphi\).
\item
  For any formula \(\varphi\), \(\ast\ast\varphi = \varphi\).
\item
  For any structural connective \([\ostar]\in [C]\) non-compatible with
  \(\beta\) and structures \(X_i\),
  \(\ast[\ostar](X_1,\,\ldots,\,X_n) = \ast[\ostar](X_1,\,\ldots,\,X_n)\).
\item
  For any structural connective \([\ostar]\in [C]\) non-compatible with
  \(\beta\) and structures \(X_i\),
  \(\ast\ast[\ostar](X_1,\,\ldots,\,X_n) = [\ostar](X_1,\,\ldots,\,X_n)\).
\end{itemize}

\ndD

\hypertarget{some-useful-connective-families}{%
\subsubsection*{Some Useful Connective
Families}\label{some-useful-connective-families}}

\bgD

Let \(\ostar\) be a connective and \(a\) an action on \(\mbb{C}\)
satisfying the condition of Definition \ref{conn_fam}. We define the
point connective family, \(\hat{\{\ostar_a\}}\) as composed of:

\begin{itemize}
\tightlist
\item
  The set \(\{\ostar\}\),
\item
  The partition \(\{\{\ostar\}\}\),
\item
  The action \(a\) assigned to \(\{\ostar\}\).
\end{itemize}

Let \(C\) and \(D\) be connective families. We define a disjoint union
connective family, \(C\sqcup D\) as composed of:

\begin{itemize}
\tightlist
\item
  The set \(C\sqcup D\),
\item
  The partition \(\mca{A}_C\sqcup\mca{A}_D\),
\item
  If \(P\in\mca{A}_C\) then the assigned group and action \(G_C\),
  \(a_C\) and if \(P\in\mca{A}_D\) then the assigned group and action
  \(G_D\), \(a_D\).
\end{itemize}

Let \(C\) be a set of connectives and \(a\) be a map from \(C\) into
actions on \(\mbb{C}\), such that for any \(\ostar\in C\), the action
\(a_{\overline\ostar}\) satisfies the condition of Definition
\ref{conn_fam}. We define the discrete connective family, \(C_a\) as
composed of:

\begin{itemize}
\tightlist
\item
  The set \(C\),
\item
  The partition \(\{\{\ostar\}\mid\ostar\in C\}\),
\item
  The action \(a(\ostar)\) for each \(\ostar\in C\).
\end{itemize}

This is the disjoint union of all point connective families in \(C\).
Whenever \(a\) is the constant into the trivial action
\(\id:\{\id\}\times V\to V : (\id,\,x)\mapsto x\) assigned to \(C\) we
note \(C_\id\) and we call it the trivial connective family.

Let \(\hat{C}\) be a connective family. Let \(V\) be a set containing at
least a propositional letter of each type in the type input components
appearing in \(C\). We define a plain \(C\) connective family as the
disjoint union connective family of \(\hat{C}\) and the discrete
connective family \(V_a\) for some action \(a\).

We define a plain point connective family, \(\{\ostar_a\}\), as a plain
connective family of \(\{\hat{\ostar_a}\}\).

Let \(C\) be a connective family. We define the full connective family
\(\mca{O}(C)\) as composed of:

\begin{itemize}
\tightlist
\item
  The set of structural connectives \([C]\),
\item
  The partition
  \(\{[a_{C,\,P}](\{[\ostar]\mid\ostar\in P\})\mid P\in\mca{A}_C\}\),
\item
  The respective actions
  \((g,\,\ostar)\mapsto(g,\,\mtt{sk}^{-1}(\ostar))\mapsto\mtt{sk}([a_{C,\,P}](g,\,[\ostar]))\)
  for each \(P\in\mca{A}_{\mca{O}(C)}\).
\end{itemize}

\ndD

\begin{quote}
By abuse of notation we will use the connective notation for the
elements of structural connectives used inside connective families.

We usually assume Connective Families to be plain.
\end{quote}

\hypertarget{sequents-on-connective-families}{%
\subsection{Sequents on Connective
Families}\label{sequents-on-connective-families}}

\bgD

Let \(C\) be a connective family. We define \textit{\textbf{sequents}}
of \(\mca{L}_C\) as \(\varphi\vdash\psi\) for a pair of formulas or a
pair of structures \(\varphi,\,\psi\in\mca{L}_C,\,[\mca{L}_C]\) of the
same type. They will be the core of our proof systems.

Sequents of \(\mca{L}_C\) will be noted \(\mca{S}_C\) and sequents of
\([\mca{L}_C]\) will be noted \([\mca{S}_C]\).

\ndD

For any sign \(s\in\mbb{B}\) we note: \[S_{s}(X,\,Y) =
    \left\{
        \bgmat X\vdash Y&\mtt{if }s=-\\ Y\vdash X&\mtt{if }s=+\ndmat
    \right.\]

For any function \(s:\mbb{C}\to\mbb{B}\), \(i\in \{1,\,\ldots,\,n\}\),
formulas \(\varphi_i\in\mca{L}_C\) and structure \(X\in[\mca{L}_C]\), we
note: \[S_{s}(\ostar,\,\varphi_1,\,\ldots,\,\varphi_n,\,X)
= S_{s(\ostar)}(\ostar(\varphi_1,\,\ldots,\,\varphi_n),\,X)
\]

For any function \(s:[\mbb{C}]\to\mbb{B}\), \(i\in \{1,\,\ldots,\,n\}\)
and structures \(X_i,\,X\in[\mca{L}_C]\), we note:
\[S_{s}([\ostar],\,X_1,\,\ldots,\,X_n,\,X)
= S_{s(\ostar)}([\ostar](X_1,\,\ldots,\,X_n),\,X)
\]

\bgD

\label{ant_parts}

Let \(X\) and \(Y\) be structures or formulas. If
\(X\vdash Y\in\mca{S}_C,\,[\mca{S}_C]\), we call \(X\) the
\textit{\textbf{antecedent}} and \(Y\) the \textit{\textbf{consequent}}
of the sequent \(X\vdash Y\).

Let \(X\) be a structure and \(Z\) a substructure occurrence of \(X\).
We define by induction on \(X\), \(s(X,\,Z)\in\mbb{B}\):

\begin{itemize}
\tightlist
\item
  If \(Z\) is equal to \(X\), \(s(X,\,Z) = +\).
\item
  If \(X = *Y\) for some structure \(Y\), then \(s(X,\,Z) = -s(Y,\,Z)\).
\item
  If \(X = \ostar(X_1,\,\ldots,\,X_{n})\) and there's some
  \(i\in \{1,\,\ldots,\,n\}\) such that \(Z\) is a subtructure
  occurrence of \(X_i\), \(s(X,\,Z) = \pm_i(\ostar)s(X_i,\,Z)\).
\end{itemize}

Let \(X\vdash Y\) be a sequent,

\begin{itemize}
\tightlist
\item
  If \(Z\) is a subtructure occurrence of \(X\) we define
  \(s(X\vdash Y,\,Z) = -s(X,\,Z)\).
\item
  If \(Z\) is a subtructure occurrence of \(Y\) we define
  \(s(X\vdash Y,\,Z) = s(Y,\,Z)\).
\end{itemize}

We say that \(Z\) is an \textit{\textbf{antecedent part}} of a sequent
\(S\) whenever \(s(S,\,Z) = -\) and that \(Z\) is a
\textit{\textbf{consequent part}} of a sequent whenever
\(s(S,\,Z) = +\).

\ndD

\begin{quote}
We informally understand the antecedent parts as the ``left side'' and
the consequent parts as the ``right side'' of a sequent.
\end{quote}

\bgD

\label{act_seq_def}

Let \(C\) be a connective family, with partition \(\mca{A}\) and for
each \(D\in\mca{A}\) group \(G_D\) and action \(a_D\). Let \(S\) be a
sequent. We define a partial function
\([a_D]:G_D\times[\mca{S}_C]\to[\mca{S}_C]\) by cases on the structures
composing \(S\):

\begin{itemize}
\tightlist
\item
  Let \(s\in\mbb{B}\), \([\ostar]\in[C]\) and
  \(X_1,\,\ldots,\,X_{n+1}\in[\mca{L}_C]\). If
  \(S = S_s([\ostar],\,X_1,\,\ldots,\,X_{n+1})\) and
  \(g\in G_{\overline{\ostar}}\), let
  \(g[\ostar]:= [a_{\overline\ostar}](g,\,[\ostar])\),
  \(\pm,\,\pm'\in(\ent\lres2\ent)^{n+1}\) be the tonicity signatures,
  \(\sigma\), \(\sigma'\) be the permutations and \(b\) and \(b'\) be
  the last components of \([\ostar]\) and \(g[\ostar]\), respectively,
  then \begin{align*}[a_D](g&,\,S)\\
  =&S_{\mtt{\AE}(g\ostar) + \mtt{\AE}(\ostar) + s}(g[\ostar],\,\\
  &T(\pm')_{1}+T(R(\sigma'\cdot\sigma^{-1},\,\pm))_{1}+b'+b+X_{\sigma'\cdot\sigma^{-1}(1)},\,\ldots,\,\\
  &T(\pm')_{n+1}+T(R(\sigma'\cdot\sigma^{-1},\,\pm))_{n+1}+b'+b+X_{\sigma'\cdot\sigma^{-1}(n+1)})
  \end{align*} Where for every structure \(X\) we consider \(0+X\) as
  \(X\) and \(1+ X\) as \(*X\), \(T\) is found in definition \ref{q_t}
  and \(R\) in definition \ref{r_alpha}.
\end{itemize}

\ndD

\bgP

\label{act_seq}

Let \(C\) be a connective family with a partition \(\mca{A}\) and for
each \(D\in\mca{A}\) a group \(G_D\) and an action \(a_D\). For each
\(D\in\mca{A}\), the partial function
\([a_D]:G_D\times[\mca{S}_C]\to[\mca{S}_C]\) is well-defined and is a
partial action.

\ndP

\begin{quote}
It will be used for the definition of our display calculi.
\end{quote}

\begin{proof}

To see that \(\forall\ostar\in D\) the image of
\([a_D](g,\,\ostar(\varphi_1,\,\ldots,\,\varphi_n))\) lives in
\([\mca{S}]\) we recall the property on structure families that says
\(\forall \ostar,\,\ostar'\in a_D(\mtt{sk}(D)),\,\ostar\) and
\(\ostar'\) have equal arity \(n\) and have permutations
\(\sigma,\,\sigma'\in\mtt{Sym}(n+1)\) and type tuples
\(k,\,k'\in(\nat^*)^{n+1}\) satisfying
\(P_{\nat^*}(\sigma'\cdot\sigma^{-1},\,k) = k'\).

We can see that this action is only going to be defined in the sequents
\(S_s([\ostar],\,X_1,\,\ldots,\,X_{n+1})\) for \([\ostar]\in [D]\). On
those it is well defined thanks to the restriction on type tuples
mentioned before.

Let \(\sigma\), \(\sigma'\) and \(\sigma''\) be the permutation,
\(\pm\), \(\pm'\) and \(\pm''\) be the tonicity signature and \(b\),
\(b'\) and \(b''\) be the last component of \([\ostar]\), \(h[\ostar]\)
and \(g\cdot h[\ostar]\), respectively. \begin{align*}
&[a_D](g,\,[a_D](h,\,S_s([\ostar],\,X_1,\,\ldots,\,X_{n+1})))\\
=& [a_D](g,\,S_{\mtt{\AE}(g\ostar) + \mtt{\AE}(\ostar) + s}([a_D](h,\,[\ostar]),\,T(\pm')_{1}+T(R(\sigma'\cdot\sigma^{-1},\,\pm))_{1}+b'+b+X_{\sigma'\cdot\sigma^{-1}(1)},\,\ldots,\,\\
  &T(\pm')_{n+1}+T(R(\sigma'\cdot\sigma^{-1},\,\pm))_{n+1}+b'+b+X_{\sigma'\cdot\sigma^{-1}(n+1)}))\\
= &S_{\mtt{\AE}(g\cdot h\ostar)+\mtt{\AE}(h\ostar)+\mtt{\AE}(h\ostar)+\mtt{\AE}(\ostar)+s}([a_D](g,\,[a_D](h,\,[\ostar])),\,T(\pm'')_{1}+T(R(\sigma''\cdot\sigma'^{-1},\,\pm'))_{1}\\
  &+T(\pm')_{\sigma''\cdot\sigma'^{-1}(1)}+T(R(\sigma'\cdot\sigma^{-1},\,\pm))_{\sigma''\cdot\sigma'^{-1}(1)}+b''+b'+b'+b+X_{\sigma''\cdot\sigma'^{-1}\cdot\sigma'\cdot\sigma^{-1}(1)},\,\ldots,\,\\
  &T(\pm'')_{n+1}+T(R(\sigma''\cdot\sigma'^{-1},\,\pm'))_{n+1}\\
  &+T(\pm')_{\sigma''\cdot\sigma'^{-1}(n+1)}+T(R(\sigma'\cdot\sigma^{-1},\,\pm))_{\sigma''\cdot\sigma'^{-1}(n+1)}+b''+b'+b'+b+X_{\sigma''\cdot\sigma'^{-1}\cdot\sigma'\cdot\sigma^{-1}(n+1)})\\
= &S_{\mtt{\AE}(g\cdot h\ostar)+\mtt{\AE}(\ostar)+s}([a_D](g\cdot h,\,[\ostar]),\,T(\pm'')_{1}+T(R(\sigma''\cdot\sigma'^{-1},\,\pm'))_{1}\\
  &+P(\sigma''\cdot\sigma'^{-1},\,T(\pm'))_{1}+P(\sigma''\cdot\sigma'^{-1},\,T(R(\sigma'\cdot\sigma^{-1},\,\pm)))_{1}+b''+b+X_{\sigma''\cdot\sigma^{-1}(1)},\,\ldots,\,\\
  &T(\pm'')_{n+1}+T(R(\sigma''\cdot\sigma'^{-1},\,\pm'))_{n+1}\\
  &+P(\sigma''\cdot\sigma'^{-1},\,T(\pm'))_{n+1}+P(\sigma''\cdot\sigma'^{-1},\,T(R(\sigma'\cdot\sigma^{-1},\,\pm)))_{n+1}+b''+b+X_{\sigma''\cdot\sigma^{-1}(n+1)})\\
= &S_{\mtt{\AE}(g\cdot h\ostar)+\mtt{\AE}(\ostar)+s}([a_D](g\cdot h,\,[\ostar]),\,T(\pm'')_{1}+T(R(\sigma''\cdot\sigma'^{-1},\,\pm'))_{1}\\
  &+T(Q(\sigma''\cdot\sigma'^{-1},\,\pm'))_{1}+T(Q(\sigma''\cdot\sigma'^{-1},\,R(\sigma'\cdot\sigma^{-1},\,\pm)))_{1}+b''+b+X_{\sigma''\cdot\sigma^{-1}(1)},\,\ldots,\,\\
  &T(\pm'')_{n+1}+T(R(\sigma''\cdot\sigma'^{-1},\,\pm'))_{n+1}\\
  &+T(Q(\sigma''\cdot\sigma'^{-1},\,\pm'))_{n+1}+T(Q(\sigma''\cdot\sigma'^{-1},\,R(\sigma'\cdot\sigma^{-1},\,\pm)))_{n+1}+b''+b+X_{\sigma''\cdot\sigma^{-1}(n+1)})\\
= &S_{\mtt{\AE}(g\cdot h\ostar)+\mtt{\AE}(\ostar)+s}([a_D](g\cdot h,\,[\ostar]),\,\\
  &T(\pm'')_{1}+T(R(\sigma''\cdot\sigma'^{-1},\,0)+Q(\sigma''\cdot\sigma'^{-1},\,R(\sigma'\cdot\sigma^{-1},\,\pm)))_{1}+b''+b+X_{\sigma''\cdot\sigma^{-1}(1)},\,\ldots,\,\\
  &T(\pm'')_{n+1}+T(R(\sigma''\cdot\sigma'^{-1},\,0)+Q(\sigma''\cdot\sigma'^{-1},\,R(\sigma'\cdot\sigma^{-1},\,\pm)))_{n+1}+b''+b+X_{\sigma''\cdot\sigma^{-1}(n+1)})\\
= &S_{\mtt{\AE}(g\cdot h\ostar)+\mtt{\AE}(\ostar)+s}([a_D](g\cdot h,\,[\ostar]),\,\\
  &T(\pm'')_{1}+T(R(\sigma''\cdot\sigma'^{-1},\,R(\sigma'\cdot\sigma^{-1},\,\pm)))_{1}+b''+b+X_{\sigma''\cdot\sigma^{-1}(1)},\,\ldots,\,\\
  &T(\pm'')_{n+1}+T(R(\sigma''\cdot\sigma'^{-1},\,R(\sigma'\cdot\sigma^{-1},\,\pm)))_{n+1}+b''+b+X_{\sigma''\cdot\sigma^{-1}(n+1)})\\
= &S_{\mtt{\AE}(g\cdot h\ostar)+\mtt{\AE}(\ostar)+s}([a_D](g\cdot h,\,[\ostar]),\,T(\pm'')_{1}+T(R(\sigma''\cdot\sigma^{-1},\,\pm))_{1}+b''+b+X_{\sigma''\cdot\sigma^{-1}(1)},\,\ldots,\,\\
  &T(\pm'')_{n+1}+T(R(\sigma''\cdot\sigma^{-1},\,\pm))_{n+1}+b''+b+X_{\sigma''\cdot\sigma^{-1}(n+1)})\\
= &[a_D](g\cdot h,\,S_s([\ostar],\,X_1,\,\ldots,\,X_{n+1}))
\end{align*} where we used proposition \ref{perm_q_t} in the fourth step
and proposition \ref{q_from_r} in the fifth and sixth steps and
proposition \ref{alpha_act} for using \(R\) as an action.

\ndpr

\bgL

\label{sign_wf}

Let \(C\) be a connective family. Let \([\ostar]\in [C]\) and
\(X_1,\,\ldots,\,X_{n+1}\in[\mca{L}_C]\). For any
\(i\in\{1,\,\ldots,\,n+1\},\,g\in G_{\overline\ostar},\,s(S_{\mtt{\AE}}([\ostar],\,X_1,\,\ldots,\,X_{n+1}),\,X_i) = s([a_D](g,\,S_{\mtt{\AE}}([\ostar],\,X_1,\,\ldots,\,X_{n+1})),\,X_i)\).

\ndL

\begin{proof}

Let \(\pm\) and \(\pm'\) be the tonicity signature, \(\sigma\) and
\(\sigma'\) be the permutation and \(b\) and \(b'\) be the last
component of \([\ostar]\) and \([a_D](g,\,[\ostar])\), respectively. In
this proof we use the \(\ent\lres2\ent\) notation for signs.

In the left side, we have
\(s(S_{\mtt{\AE}}([\ostar],\,X_1,\,\ldots,\,X_{n+1}),\,X_i) = \delta_{i,\,n+1} + (1 + \delta_{i,\,{n+1}})\pm_i(\ostar) + \mtt{\AE}(\ostar)\)
by definition \ref{ant_parts} (it can be seen by case analysis on
\(\delta_{i,\,n+1}\)). In the right side we have \begin{align*}
&s([a_D](g,\,S_{\mtt{\AE}}([\ostar],\,X_1,\,\ldots,\,X_{n+1})),\,X_i)\\
=&s(S_{\pm'_{n+1}+\pm_{n+1}+b'+b+\mtt{\AE}(\ostar)}(g[\ostar],\,\\
  &T(\pm')_{1}+T(R(\sigma'\cdot\sigma^{-1},\,\pm))_{1}+b'+b+X_{\sigma'\cdot\sigma^{-1}(1)},\,\ldots,\,\\
  &T(\pm')_{n+1}+T(R(\sigma'\cdot\sigma^{-1},\,\pm))_{n+1}+b'+b+X_{\sigma'\cdot\sigma^{-1}(n+1)}),\,X_i)
\end{align*} We must notice now that \(X_i\) will still only be a
structure ocurrence of itself, so that its sign corresponds to
\begin{align*}
&s([a_D](g,\,S_{\mtt{\AE}}([\ostar],\,X_1,\,\ldots,\,X_{n+1})),\,X_i)\\
=&s(T(\pm')_{\sigma\cdot\sigma'^{-1}(i)}+T(R(\sigma'\cdot\sigma^{-1},\,\pm))_{\sigma\cdot\sigma'^{-1}(i)}+b'+b+X_i,\,X_i)\\
&+ \delta_{\sigma\cdot\sigma'^{-1}(i),\,n+1} + (1+\delta_{\sigma\cdot\sigma'^{-1}(i),\,n+1})\pm'_{\sigma\cdot\sigma'^{-1}(i)} + \pm'_{n+1} + \pm_{n+1}+b'+b+\mtt{\AE}(\ostar)\\
=&T(\pm')_{\sigma\cdot\sigma'^{-1}(i)}+P(\sigma\cdot\sigma'^{-1},\,T(R(\sigma'\cdot\sigma^{-1},\,\pm)))_{i}+b'+b\\
&+ \delta_{\sigma\cdot\sigma'^{-1}(i),\,n+1} + (1+\delta_{\sigma\cdot\sigma'^{-1}(i),\,n+1})\pm'_{\sigma\cdot\sigma'^{-1}(i)} + \pm'_{n+1} + \pm_{n+1}+b'+b+\mtt{\AE}(\ostar)\\
=&T(\pm')_{\sigma\cdot\sigma'^{-1}(i)}+T(Q(\sigma\cdot\sigma'^{-1},\,R(\sigma'\cdot\sigma^{-1},\,\pm)))_{i}\\
&+ \delta_{\sigma\cdot\sigma'^{-1}(i),\,n+1} + (1+\delta_{\sigma\cdot\sigma'^{-1}(i),\,n+1})\pm'_{\sigma\cdot\sigma'^{-1}(i)} + \pm'_{n+1} + \pm_{n+1}+\mtt{\AE}(\ostar)\\
=&T(Q(\sigma\cdot\sigma'^{-1},\,R(\sigma'\cdot\sigma^{-1},\,\pm)))_{i} + \delta_{\sigma\cdot\sigma'^{-1}(i),\,n+1} + \pm_{n+1}+\mtt{\AE}(\ostar)\\
=&T(Q(\sigma\cdot\sigma'^{-1},\,Q(\sigma'\cdot\sigma^{-1},\,\pm)))_{i} +T(Q(\sigma\cdot\sigma'^{-1},\,((1 + \delta_{\sigma'\cdot\sigma^{-1}(j),\,n+1})(1 + \delta_{\sigma'\cdot\sigma^{-1}(n+1),\,n+1}))_j^{n+1}))_i\\
&+ \delta_{\sigma\cdot\sigma'^{-1}(i),\,n+1} + \pm_{n+1}+\mtt{\AE}(\ostar)\\
=&T(\pm)_{i} + P(\sigma\cdot\sigma'^{-1},\,((1 + \delta_{\sigma'\cdot\sigma^{-1}(j),\,n+1})(1 + \delta_{\sigma'\cdot\sigma^{-1}(n+1),\,n+1}) + (1+\delta_{j,\,n+1})(1 + \delta_{\sigma'\cdot\sigma^{-1}(n+1),\,n+1}))_j^{n+1})_i\\
&+ \delta_{\sigma\cdot\sigma'^{-1}(i),\,n+1} + \pm_{n+1}+\mtt{\AE}(\ostar)\\
=&T(\pm)_{i} + (1 + \delta_{i,\,n+1})(1 + \delta_{\sigma'\cdot\sigma^{-1}(n+1),\,n+1}) + (1+\delta_{\sigma\cdot\sigma'^{-1}(i),\,n+1})(1 + \delta_{\sigma'\cdot\sigma^{-1}(n+1),\,n+1})\\
&+ \delta_{\sigma\cdot\sigma'^{-1}(i),\,n+1} + \pm_{n+1}+\mtt{\AE}(\ostar)\\
=&(1+\delta_{i,\,n+1})\pm_i + \pm_{n+1} + \delta_{i,\,n+1}(1 + \delta_{\sigma'\cdot\sigma^{-1}(n+1),\,n+1}) + \delta_{\sigma\cdot\sigma'^{-1}(i),\,n+1}(1 + \delta_{\sigma'\cdot\sigma^{-1}(n+1),\,n+1})\\
&+ \delta_{\sigma\cdot\sigma'^{-1}(i),\,n+1} + \pm_{n+1}+\mtt{\AE}(\ostar)\\
=&(1+\delta_{i,\,n+1})\pm_i + \pm_{n+1} + \delta_{i,\,n+1}(1 + \delta_{\sigma'\cdot\sigma^{-1}(n+1),\,n+1}) + \delta_{\sigma\cdot\sigma'^{-1}(i),\,n+1}\delta_{\sigma'\cdot\sigma^{-1}(n+1),\,n+1} + \pm_{n+1}+\mtt{\AE}(\ostar)\\
=&\delta_{i,\,n+1} + (1+\delta_{i,\,n+1})\pm_{i}+\mtt{\AE}(\ostar)
\end{align*} where we used again proposition \ref{perm_q_t} and the
\((b(j))_j^{n+1}\) notation for some function
\(b:\{1,\,\ldots,\,n+1\}\to\ent\lres2\ent\) is notation for
\((b(1),\,\ldots,\,b(n+1))\).

\ndpr

\bgE

\label{act_seq_alpha_varsigma}

The partial action on sequents given in definition \ref{act_seq_def} for
the action \(\alpha\times\varsigma\) is \begin{equation*}
\begin{matrix}[a_D]((s_1,\,\ldots,\,s_{n+1})\rho,\,S_{\mtt{\AE}}([\ostar],\,X_1,\,\ldots,\,X_{n+1}))\\
\rotatebox[origin=c]{-90}{$\mapsto$}\\
S_{\mtt{\AE}}((s_1,\,\ldots,\,s_{n+1})\rho[\ostar],\,s_{n+1}+s_1+X_{\sigma(n)},\,\ldots,\,s_{n+1}+s_n+X_{\rho(n)},\,s_{n+1}+X_{\rho(n+1)})\end{matrix}\end{equation*}

\ndE

\begin{proof}

Let \(\ostar\) be a connective of arity \(n\), permutation \(\sigma\),
tonicity signature \(\pm\) and action \(\alpha\times\varsigma\). Let
\(\pm'\) be the tonicity signature and \(\sigma'\) be the permutation of
\((s_1,\,\ldots,\,s_{n+1})\sigma[\ostar]\).

We just have to see that
\(\pm' + R(\sigma'\cdot\sigma^{-1},\,\pm) = (s_1,\,\ldots,\,s_{n+1})\),
as we already know that
\(T(s_1,\,\ldots,\,s_{n+1}) = (s_1+s_{n+1},\,\ldots,\,s_n+s_{n+1},\,s_{n+1})\).
We have \(\pm' = (s_1,\,\ldots,\,s_{n+1})+R(\rho,\,\pm)\) and
\(\sigma' = \rho\cdot\sigma\). Therefore
\(\sigma'\cdot\sigma^{-1} = \rho\) and
\(\pm' + R(\sigma'\cdot\sigma^{-1},\,\pm) = (s_1,\,\ldots,\,s_{n+1}) + R(\rho,\,\pm) + R(\rho,\,\pm) = (s_1,\,\ldots,\,s_{n+1})\)
as we wanted.

\ndpr

\hypertarget{atomic-logics-for-connective-families}{%
\section{Atomic Logics for Connective
Families}\label{atomic-logics-for-connective-families}}

\label{sect5}

The definition of families of structural connectives will help us better
understand what are the structures and how to treat them in our calculi.
Through connective families we will be capable to compare the logics
given by both versions of the atomic logics. All this will let us
afterwards redefine Atomic Logics syntax.

For an example on how connective families work and what they are useful
for the reader might want to go to the next section.

\hypertarget{syntax}{%
\subsection{Syntax}\label{syntax}}

We remind that a connective family is composed of:

\begin{itemize}
\tightlist
\item
  A set \(C\),
\item
  A function \(\mtt{sk}:C\to\mbb{C}\),
\item
  A partition of \(C\), \(\mca{A}\),
\item
  For all \(D\in\mca{A}\), a group \(G\) and an action
  \(a:G\times\mbb{C}\to\mbb{C}\).
\end{itemize}

And it satisfies all the conditions imposed by Definition
\ref{conn_fam}.

Then we set the syntax of Atomic Logics to be composed of:

\begin{itemize}
\tightlist
\item
  A Connective Family \(C\),
\item
  Its Family of Structural Connectives \([C]\),
\item
  Its formulas \(\mca{L}_C\),
\item
  Its structures \([\mca{L}_C]\),
\item
  Its sequents \(\mca{S}_C\) and \([\mca{S}_C]\).
\end{itemize}

\hypertarget{an-example-lambek-and-bi-lambek-connectives}{%
\subsubsection*{An Example, Lambek and Bi-Lambek
Connectives}\label{an-example-lambek-and-bi-lambek-connectives}}

We will further comment on the connective Family for Lambek Calculus,
introduced in example \ref{lamb_fam}.

We show in the following calculations that the Lambek Connectives form
an orbit of the \(\alpha\) action; as defined in the precedent chapter,
we will work on the tonicity signatures as vectors and the morphism
\(Q\) in its matrix form:

{\vspace{-2em}\scriptsize
$$\alpha((2\ 3),\,\otimes)
= ((2\ 3)\cdot\id,\,P((2\ 3))\begin{pmatrix}1\\1\\1\end{pmatrix},\,
        \begin{pmatrix}1\\0\\1\end{pmatrix} + Q((2\ 3))\begin{pmatrix}0\\0\\0\end{pmatrix})
= ((2\ 3),\,\begin{pmatrix}1\\1\\1\end{pmatrix},\,
        \begin{pmatrix}1\\0\\1\end{pmatrix} + \begin{pmatrix}1&1&0\\0&1&0\\0&1&1\end{pmatrix}\begin{pmatrix}0\\0\\0\end{pmatrix})
= \rres$$
$$\alpha((1\ 3),\,\otimes)
= ((1\ 3)\cdot\id,\,P((1\ 3))\begin{pmatrix}1\\1\\1\end{pmatrix},\,
        \begin{pmatrix}0\\1\\1\end{pmatrix} + Q((1\ 3))\begin{pmatrix}0\\0\\0\end{pmatrix})
= ((1\ 3),\,\begin{pmatrix}1\\1\\1\end{pmatrix},\,
        \begin{pmatrix}0\\1\\1\end{pmatrix} + \begin{pmatrix}1&0&0\\1&1&0\\1&0&1\end{pmatrix}\begin{pmatrix}0\\0\\0\end{pmatrix})
= \lres$$}

We define the \emph{Bi-Lambek Connective Family} \(BL\) as:

\begin{itemize}
\tightlist
\item
  the set
  \(\{\otimes,\,\rres,\,\lres,\,\oplus,\,\succ,\,\prec\}\sqcup V\),
  where:

  \begin{itemize}
  \tightlist
  \item
    the connective in \(L\) have the same skeletons,
  \item
    the cotensor
    \(\mtt{sk}(\oplus)=(\id,\,(1,\,1,\,1),\,(+,\,+,\,-),\,-)\),
  \item
    the right coresidual
    \(\mtt{sk}(\succ)=((2\ 3),\,(1,\,1,\,1),\,(-,\,+,\,+),\,-)\),
  \item
    the left coresidual
    \(\mtt{sk}(\prec)=((1\ 3),\,(1,\,1,\,1),\,(+,\,-,\,+),\,-)\),
  \end{itemize}
\item
  the partition
  \(\{\{\otimes,\,\rres,\,\lres,\,\oplus,\,\prec,\,\succ\}\}\sqcup\{\{v\}\mid V\}\),
\item
  the action \(\alpha\times\varsigma\) on
  \(\{\otimes,\,\rres,\,\lres,\,\oplus,\,\succ,\,\prec\}\) and the
  action \(\id\) on each \(\{v\}\), for \(v\in V\).
\end{itemize}

Just by using \(\alpha\) we reach a skeleton outside of both families.
Indeed, \[
\alpha((3\ 2\ 1),\,\otimes) = ((3\ 2\ 1)\cdot\id,\,(1,\,1,\,1),\,
        \begin{pmatrix}1\\0\\1\end{pmatrix} + \begin{pmatrix}0&1&0\\1&1&0\\0&1&1\end{pmatrix}\begin{pmatrix}0\\0\\0\end{pmatrix})= ((3\ 2\ 1),\,(1,\,1,\,1),\,
        \begin{pmatrix}0\\1\\1\end{pmatrix})
\] Another case of connective falling outside both families was given
through \(\alpha*\beta\) in example \ref{lambek_ex1}

For this reason we needed structural connectives for the calculi to be
able to act on connectives. We will later use those actions to prove the
display theorem.

We could have also defined the Lambek Connective Family
\(L\cap\alpha(L)\) as:

\begin{itemize}
\tightlist
\item
  the same set \(\{\otimes,\,\rres,\,\lres\}\sqcup V\),
\item
  the same partition
  \(\{\{\otimes,\,\rres,\,\lres\}\}\sqcup\{\{v\}\mid V\}\),
\item
  the action \(\alpha\) on \(\{\otimes,\,\rres,\,\lres\}\) and the
  action \(\id\) on each \(\{v\}\), for \(v\in V\).
\end{itemize}

Then we can check that \(L\lesssim L\cap\alpha(L)\) by:

\begin{itemize}
\tightlist
\item
  the bijection \(\id: L\to L\cap\alpha(L)\),
\item
  the morphism
  \(\iota:\mtt{Sym}(3)\hookrightarrow(\ent\lres2\ent)^3\rtimes\mtt{Sym}(3):\sigma\mapsto ((0,\,0,\,0),\,\sigma)\).
\end{itemize}

We can also check that \(L\lesssim BL\) through the set inclusion and
the identity morphisms.

\hypertarget{semantics}{%
\subsection{Semantics}\label{semantics}}

\label{semantics}

Now we formalize generalized Kripke semantics for all atomic logics,
using its formulas \(\mca{L}_C\) and using an algebra of \(k\)-ary
relations on sets (called structures).

\begin{quote}
This is Aucher's Semantics for Atomic Logics extended to all connective
families.
\end{quote}

\bgD

Let \(C\) be a family of connectives of action \(a\) on the group \(G\).
A \(C\)-model for Atomic Logics is a tuple \(M = (W,\,\mca{R},\,R)\),
where \(W\) is a non-empty set,
\(\mca{R}\subseteq \bigsqcup_{i\in\nat^*}\mca{P}(W^i)\) and a surjective
map \(R:C\to\mca{R}\) such that
\(\forall\ostar,\,\ostar'\in C,\,\overline{\ostar'}=\overline{\ostar}\)
implies \(R(\ostar') = R(\ostar)\) and
\(R(\ostar)\in W^{\sum_{i\in\{1,\,\ldots,\,n+1\}}k_i}\), where
\((k_1,\,\ldots,\,k_{n+1})\) is the type tuple of \(\ostar\).

A \(k\)-assignment is a tuple from \(W^k\). A pointed \(C\)-model of
type \(k\) is a model along with a \(k\)-assignment.

The class of pointed \(C\)-models is noted \(\mca{E}_C\).

\ndD

\begin{quote}
For a relation \(R\subseteq W^N\) we'll define the action of \(\pm\) on
the set \(\mca{P}(W^N)\) by \(+R = R\) and \(-R = W^N-R\) and the action
of \(\mtt{Sym}(n)\) by
\(\sigma R\ w_1\ \ldots\ w_N = R\ w_{\sigma(1)}\ \ldots\ w_{\sigma(N)}\).

For a \(n\)-ary connective \(\ostar\), all \(i\in{1,\,\ldots,\,n}\), all
elements \(a\) and sets \(A\) we define
\(a\pitchfork^{\ostar}_i A = \left\{\bgmat a\in A&\mtt{if }\pm_i(\ostar) = +\\ a\notin A&\mtt{if }\pm_i(\ostar) = -\ndmat\right.\).

This notations will be used in the following definition.
\end{quote}

\bgD

Let \(C\) be a Connective Family and \(M\) be a \(C\)-model. We define
the interpretation function of \(\mca{L}_C\) in \(M\) by
\(\llbracket\,.\,\rrbracket^M:\mca{L}_C\to\bigcup_{k\in\nat^*}W^k\)
inductively:

\begin{itemize}
\tightlist
\item
  For all propositional letters \(p\in C\) of type \(k\),
  \(\llbracket p\rrbracket^m = \left\{\bgmat R(p)&\text{if }\pm(p) = +\\ -R(p)&\text{if }\pm(p) = -\ndmat\right.\).
\item
  For all \(n>0\), all \(n\)-ary connectives
  \((\sigma,\,\pm,\,\text{\AE},\,(k,k_1,\,\ldots,\,k_{n+1}),\,(\pm_1,\,\ldots,\,\pm_{n}))\in C\)
  and all \((\varphi_1,\,\ldots,\,\varphi_{n+1})\in \mca{L}_C^n\):
  \(\llbracket\ostar(\varphi_1,\,\ldots,\,\varphi_n)\rrbracket^M = \{\overline{w}\in W^k\mid\mca{C}^{\ostar}(\llbracket\varphi_{1}\rrbracket,\,\ldots,\,\llbracket\varphi_{n}\rrbracket,\,\overline{w})\}\),
  where \(\mca{C}^{\ostar}\) is called truth condition of the connective
  and it is defined by cases:

  \begin{itemize}
  \tightlist
  \item
    \(\mca{C}^{\ostar}(W_1,\,\ldots,\,W_n,\,\overline{w}) = \forall \overline{w_1}\in W^{k_1},\,\ldots,\,\overline{w_n}\in W^{k_n},\,\overline{w_1}\pitchfork^{\ostar}_1 W_1\vee\ldots\vee\overline{w_n}\pitchfork^{\ostar}_n W_n\vee \pm\sigma R_{\overline\ostar}\ \overline{w_1}\ \ldots\ \overline{w_n}\ \overline{w}\)
    when \(\mtt{\AE}(\ostar) = \forall\).
  \item
    \(\mca{C}^{\ostar}(W_1,\,\ldots,\,W_n,\,\overline{w}) = \exists \overline{w_1}\in W^{k_1},\,\ldots,\,\overline{w_n}\in W^{k_n},\,\overline{w_1}\pitchfork^{\ostar}_1 W_1\wedge\ldots\wedge\overline{w_n}\pitchfork^{\ostar}_n W_n\wedge \pm\sigma R_{\overline\ostar}\ \overline{w_1}\ \ldots\ \overline{w_n}\ \overline{w}\)
    when \(\mtt{\AE}(\ostar) = \exists\).
  \end{itemize}
\end{itemize}

We define the interpretation function of \(\mca{S}_C\) by
\(w\in\llbracket\varphi\vdash\psi\rrbracket^M\) iff if
\(w\in\llbracket\varphi\rrbracket^M\) then
\(w\in\llbracket\psi\rrbracket^M\).

\ndD

\begin{quote}
We can extend this definition to \([C]\)-models by translating
structures into \(\mca{O}(C)\)-formulas.
\end{quote}

\bgD

The satisfaction relation \(\vDash\subseteq\mca{E}_C\times\mca{L}_C\) is
defined
\(\forall (M,\,\overline{w})\in\mca{E}_C,\,\varphi\in\mca{L}_C,\,M,\,\overline{w}\vDash\varphi\)
if and only if \(\overline{w}\in\llbracket\varphi\rrbracket^M\). In that
case, we say that \(\varphi\) is \textit{\textbf{true}} on
\((M,\,\overline{w})\).

\ndD

\begin{quote}
With the formulas \(\mca{L}_C\), the evaluations \(\mca{E}_C\) and the
satisfaction relation \(\vDash\) we have succesfully built the Atomic
Logics associated to the set of connectives \(C\).
\end{quote}

Let two logics \(L = (\mca{L},\mca{E},\,\vDash)\) and
\(L' = (\mca{L}',\,\mca{E},\,\vDash')\) be given over the same class of
models. We say that \(L'\) is \textit{\textbf{at least as expressive}}
as \(L\) if \(\forall\varphi\in\mca{L}\) there is a
\(\varphi'\in\mca{L}'\) such that
\(\forall M\in\mca{E},\,M\vDash\varphi\) if and only if
\(M\vDash\varphi'\).

We say that \(L\) and \(L'\) are \textit{\textbf{equi-expressive}} when
\(L\) is at least as expressive as \(L'\) and viceversa. We say that
\(L\) is \textit{\textbf{strictly less expressive}} than \(L'\) whenever
\(L'\) is at least as expressive as \(L\) but they are not
equi-expressive.

\bgP

\label{bijec_equi}

Let \(C\) and \(C'\) be two connective families. If there is a bijection
\(T:C\to C'\) such that
\(\forall\ostar\in C,\,\overline{\ostar} = \overline{\ostar'}\) if and
only if \(\overline{T(\ostar)} = \overline{T(\ostar')}\) and
\(\mtt{sk}(T(\ostar)) = \mtt{sk}(\ostar)\), \(\mca{L}_C\) and
\(\mca{L}_{C'}\) are equi-expressive.

\ndP

\begin{proof}

Firstly, we observe that the classes of models are the same. We recall
that classes of models correspond to the class of triples
\((W,\,\mca{R},\,R)\) of world sets \(W\), relation sets \(\mca{R}\) and
surjective functions \(R:C\to \mca{R}\), where
\(R(\ostar)\in W^{\sum_{i\in\{1,\,\ldots,\,n+1\}}t_i}\) and the
conditions \(R(a_D(\ostar)) = R(\ostar)\), where \(D\in \mca{A}_C\) and
\(\ostar\in C\). Similarly for \(C'\).

Then \(\mca{E}_C\) are all the models with, for each \(\ostar\in C\),
exactly a \(\sum_{i\in\{1,\ldots,\,n+1\}}k_i\)-ary relation, where
\((k_1,\,\ldots,\,k_{n+1})\) is the type tuple of \(\ostar\). Similarly
\(\mca{E}_C'\) are all the models with, for each \(\ostar\in C'\),
exactly a \(\sum_{i\in\{1,\ldots,\,n+1\}}k_i\)-ary relation, where
\((k_1,\,\ldots,\,k_{n+1})\) is the type tuple of \(\ostar\).

As the bijection \(T\) preserves the partition, the relations of
\(\ostar\) and \(\ostar'\) are equal if and only if those of
\(T(\ostar)\) and \(T(\ostar')\) are. Finally, as \(T\) preserves the
skeletons, \(R_{\overline\ostar}\) and \(R_{\overline{T(\ostar)}}\) are
relations of the same arity (the sum of the type tuples). Therefore,
they form the same class.

We give as translation between frames the identity on worlds and
\(T(R_{\overline\ostar}) = R_{\overline{T(\ostar)}}\) on relationships.

Now we give inductively the translation on formulas
\(T:\mca{L}_C\to\mca{L}_{C'}\):

\begin{itemize}
\tightlist
\item
  For any \(n\)-ary connective \(\ostar\in C\),
  \(\varphi_1,\,\ldots,\,\varphi_n\in\mca{L}_{C}\), we define
  \(T(\ostar(\varphi_1,\,\ldots,\,\varphi_n))\) as
  \(T(\ostar)(T(\varphi_1),\,\ldots,\,T(\varphi_n))\).
\end{itemize}

Let \((\mca{M},\,w)\) be a pointed model. We prove by induction on
\(\varphi\) that
\(\overline{w}\in\llbracket\varphi\rrbracket_{\mca{M}}\) if and only if
\(\overline{w}\in\llbracket T(\varphi)\rrbracket_{\mca{M}}\). The proof
just uses that
\(C^\ostar(\llbracket\varphi_1\rrbracket,\,\ldots,\,\llbracket\varphi_n\rrbracket,\,\overline{w})\)
is only determined by \(\mtt{sk}(\ostar)\) and \(R_{\overline\ostar}\),
that \(\mtt{sk}(\ostar) = \mtt{sk}(T(\ostar))\), that we have defined
the translation on frames through
\(T(R_{\overline\ostar}) = R_{\overline{T(\ostar)}}\) and the induction
step.

\ndpr

\bgP

\label{isom_equi}

Whenever \(C\sim_W C'\), \(\mca{L}_C\) and \(\mca{L}_{C'}\) are
equi-expressive.

\ndP

\begin{proof}

We just use last proposition with \(\phi\) as bijection. The condition
on skeletons correspond to the recursively using equality on actions on
the identity,
\(\mtt{sk}(\phi(\ostar)) = a(1,\,\mtt{sk}(\phi(\ostar))) = a(\varphi(1),\,\mtt{sk}(\phi(\ostar))) = a(1,\,\mtt{sk}(\ostar)) = \mtt{sk}(\ostar)\).

\ndpr

\hypertarget{sequent-calculus}{%
\subsection{Sequent Calculus}\label{sequent-calculus}}

Let \(C\) be a connective family with partition \(\mca{A}\), and for
each \(D\in\mca{A}\), associated groups \(G_D\), actions \(a_D\) and a
morphism
\(\varphi_D:G_D\times(\ent\lres2\ent)^{n+1}\to(\ent\lres2\ent)^{n+1}\).
The rules of our Sequent Calculi are \footnote{On this Calculus we are
  going to make use of proposition \ref{act_seq}.}\footnote{In the
  present calculi we will treat full formulas like structural
  propositional letters non-compatible with \(\beta\) (not the
  subformulas).}:

\begin{quote}
It must be noted that in the following rules we are using Aucher's
syntaxis. The meaning of \(\cdot\) and \(-\) must be understood as the
one in \(\mbb{B}\). Similarly for \(\mtt{\AE}\).
\end{quote}

\begin{rules}
    \textbf{Display Rules}:
    \begin{multicols}{2}
        {\underline{\textbf{The rule dsr1}}:\\
        {\small For any $[D]\in\mca{A}$ and $\forall[\ostar]\in[D]$ of arity $n$, $g\in G_D,\,X_{1},\,\ldots,\,X_{n+1}\in[\mca{L}_C]$,}
        $$\begin{prooftree}
            \hypo{
                S_{\mtt{\AE}}([\ostar],\,X_1,\,\ldots,\,X_{n+1})
            }
            \infer1[dsr1]{
                [a_D](g,\,S_{\mtt{\AE}}([\ostar],\,X_1,\,\ldots,\,X_{n+1}))
        }\end{prooftree}$$
        Where we write $1+X_{\sigma(i)}$ for $*X_{\sigma(i)}$.\\}

        \columnbreak
        \underline{\textbf{The rule dsr2}}:\\
        {\small For any $[D]\in\mca{A}$ non-compatible with the action $\beta$ and $\forall[\ostar]\in[D]$ of arity $n$, $s\in \mbb{B},\,X_{1},\,\ldots,\,X_{n+1}\in[\mca{L}_C]$,}
        $$\begin{prooftree}
            \hypo{
                S_{s}([\ostar](X_1,\,\ldots,\,X_{n}),\,X_{n+1})
            }
            \infer1[dsr2]{
                S_{-s}(*[\ostar](X_1,\,\ldots,\,X_{n}),\,*X_{n+1})
        }\end{prooftree}$$
    \end{multicols}
    \hrule\hspace{1em}

    \textbf{Left and Right Introduction Rules}:
    \begin{multicols}{2}
        \underline{\textbf{The rule $\vdash\ostar$}}:\\
        {\small For any $\ostar\in C$, $X_{1},\,\ldots,\,X_n\in[\mca{L}_C]$, $\varphi_1,\,\ldots,\,\varphi_{n}\in\mca{L}_C$,}
        $$\begin{prooftree}
            \hypo{
                S_{\mtt{\AE}\pm_1(\ostar)}(X_1,\,\varphi_1)
            }\hypo{
                \ldots
            }\hypo{
                S_{\mtt{\AE}\pm_n(\ostar)}(X_n,\,\varphi_n)
            }
            \infer3[\(\vdash\ostar\)]{
                S_{\mtt{\AE}}([\ostar],\,X_1,\,\ldots,\,X_n,\,\ostar(\varphi_1,\,\ldots,\,\varphi_n))
        }\end{prooftree}$$

        \columnbreak
        \underline{\textbf{The rule $\ostar\vdash$}}:\\
        {\small For any $\ostar\in C$, $U\in[\mca{L}_C]$, $\varphi_1,\,\ldots,\,\varphi_{n+1}\in\mca{L}_C$,}
        $$\begin{prooftree}
            \hypo{
                S_{\mtt{\AE}}([\ostar],\,\varphi_1,\,\ldots,\,\varphi_n,\,U)
            }
            \infer1[\(\ostar\vdash\)]{
                S_{\mtt{\AE}}(\ostar,\,\varphi_1,\,\ldots,\,\varphi_n,\,U)
        }\end{prooftree}$$
    \end{multicols}
    \caption*{Rules of Atomic Sequent Calculus, $ASL_C$}
\end{rules}

\bgT

\label{asl_display}

Let \(C\) be a connective family with partition \(\mca{A}\), and for
each \(D\in\mca{A}\), associated groups \(G_D\), actions \(a_D\) and a
morphism
\(\varphi_D:G_D\times(\ent\lres2\ent)^{n+1}\to(\ent\lres2\ent)^{n+1}\).
If for all \(D\in\mca{A}\), \(i\in\{1,\,\ldots,\,n\}\), there is at
least one structural connective \([\ostar]\in [D]\) such that its
permutation \(\sigma\) is in the left coset \((i\ n+1)\mtt{Sym}(n)\),
then the proof system \(ASL_C\) is properly displayable.

In all the other cases \(ASL_C\) still satisfies \(C1\) to \(C7\).

\ndT

\begin{quote}
Apart from \(C8\), we are going to use my syntax in the proof.
\end{quote}

\begin{proof}

Proof of the display property:

We take \(C\) the set of connectives. We take a constituent \(X\) of a
derivable sequent \(S_s(Y,\,Z)\), for sign \(s\) and structures \(Y\),
\(Z\). We will prove that we can derive a sequent where \(X\) is
displayed, from the premiss \(S_s(Y,\,Z)\), by induction on the
structure \(Y\).

\begin{itemize}
\tightlist
\item
  If \(Y\) is a formula and \(X\) is a substructure of \(Y\) the result
  is trivial.
\item
  If \(Y\) is a strict structure and \(X\) is a substructure of \(Y\),
  and we have \(Y = [\ostar](Y_1,\,\ldots,\,Y_n)\) for some
  \([\ostar]\in[D]\) and the sequent \(S_{\mtt{\AE}(\ostar)}(Y,\,Z)\),
  where \(X\) is a substructure of some \(X_i\). Let \(\sigma\) be the
  permutation of \([\ostar]\) and \(\pm\) its tonicity signature. We
  take some \(g\in G_D\) such that
  \(\sigma\cdot\sigma'^{-1}\in(i\ n+1)\mtt{Sym}(n)\) for \(\sigma'\) the
  permutation of \([a_D](g,\,[\ostar])\). Let \(\pm'\) be the tonicity
  signature of \([a_D](g,\,[\ostar])\). The display rule
  \(dsr1\ostar g\) lets us display any structural variable \(Y_i\) as
  follows
  \begin{align*}S_{\pm'_{n+1}\pm_{n+1}}(&[a_D](g,\,[\ostar]),\,\\
  &T(\pm')_{1}+T(R(\sigma'\cdot\sigma^{-1},\,\pm))_{1} + Y_{\sigma\cdot\sigma'^{-1}(1)},\,\ldots,\,\\
  &T(\pm')_{i-1}+T(R(\sigma'\cdot\sigma^{-1},\,\pm))_{i-1} + Y_{\sigma\cdot\sigma'^{-1}(i-1)},\,\\
  &T(\pm')_{i}+T(R(\sigma'\cdot\sigma^{-1},\,\pm))_{i} + Z,\,\\
  &T(\pm')_{i+1}+T(R(\sigma'\cdot\sigma^{-1},\,\pm))_{i+1} + Y_{\sigma\cdot\sigma'^{-1}(i+1)},\,\ldots,\,\\
  &T(\pm')_{n}+T(R(\sigma'\cdot\sigma^{-1},\,\pm))_{n} + Y_{\sigma\cdot\sigma'^{-1}(n)},\,\\
  &T(\pm')_{n+1}+T(R(\sigma'\cdot\sigma^{-1},\,\pm))_{n+1} + Y_i)\end{align*}
  We then use \(dsr2.b\), for
  \(b = T(\pm')_{n+1}+T(R(\sigma'\cdot\sigma^{-1},\,\pm))_{n+1}\).\\
  By induction hypothesis we deduce that, as \(X\) is substructure of
  \(Y_i\), it can be displayed.
\end{itemize}

This shows that the calculus has the display property.

Proof of the 8 conditions:

\begin{itemize}
\item
  \(C1\) \textbf{Preservation of formulas}:\\
  It is enough to check that neither structural variables nor
  non-parametrical formulas disappear. Let \(X_i\) be structures in the
  rules and \(\varphi_i\) formulas in the rules.

  \begin{itemize}
  \tightlist
  \item
    dsr1 and dsr2 only permute the structural variables.
  \item
    \(\vdash\ostar\) maintains both \(X_i\) and \(\varphi_i\) for all
    \(i\in \{1,\,\ldots,\,n\}\).
  \item
    \(\ostar\vdash\) maintains \(\varphi_i\) for all
    \(i\in \{1,\,\ldots,\,n\}\) and \(U\).
  \end{itemize}
\item
  \(C2\) \textbf{Shape-alikeness of parameters}:\\
  Let an inference \(\mtt{Inf}\), falling under a rule \(\mtt{Ru}\),
  have a structure variable \(\ell\) to which we assign the structure
  \(A\) and let \(B\) be a substructure of \(A\), which is a parameter
  of \(\mtt{Inf}\). We define the congruence relation of \(B\) in
  \(\mtt{Inf}\) equivalence class as composed of all the appearances of
  the structure occurrence of \(B\) in \(A\), substituted inside the
  structure variable \(\ell\).

  \begin{quote}
  For example, let us take the rule from \(X\vdash U\) infer
  \(X,\,Y\vdash U\) and the specific inference with the variable \(X\)
  subtituted by \((\varphi,\,\varphi)\), \(Y\) by
  \((\varphi,\,\varphi)\) and \(U\) by \((\varphi,\,\varphi)\), for some
  formula \(\varphi\). This rewrites to from
  \((\varphi_{X,\,1},\,\varphi_{X,\,2})\vdash(\varphi_{U,\,1},\,\varphi_{U,\,2})\)
  infer
  \((\varphi_{X,\,1},\,\varphi_{X,\,2}),\,(\varphi_{Y,\,1},\,\varphi_{Y,\,2})\vdash(\varphi_{U,\,1},\,\varphi_{U,\,2})\).
  For all structure variable \(\ell\) and \(i\in\{1,\,2\}\) all the
  parameters \(\varphi_{\ell,\,i}\) represent the formula \(\varphi\).
  The congruence relation is defined by:

  \begin{itemize}
  \tightlist
  \item
    For all structure variable \(\ell\) and \(i\in\{1,\,2\}\) the
    parameters \(\varphi_{\ell,\,i}\) are only congruent to other
    appearances of \(\varphi_{\ell,\,i}\).
  \item
    For all structure variable \(\ell\) the parameters
    \((\varphi_{\ell,\,1},\,\varphi_{\ell,\,2})\) are only congruent to
    other appearances of \((\varphi_{\ell,\,1},\,\varphi_{\ell,\,2})\).
  \end{itemize}
  \end{quote}
\item
  \(C3\) \textbf{Non-proliferation of parameters}:\\
  It is enough to see that structural variables have only one congruent
  parameter in the conclusion.

  \begin{itemize}
  \tightlist
  \item
    dsr1 and dsr2 only permute the structural variables.
  \item
    \(\vdash\ostar\) has only one parameter congruent to \(X_i\) in the
    conclusion for all \(i\in \{1,\,\ldots,\,n\}\).
  \item
    \(\ostar\vdash\) has only one parameter congruent to \(U\) in the
    conclusion.
  \end{itemize}
\item
  \(C4\) \textbf{Position-alikeness of parameters}:\\
  It is necessary to check for all rules whether the signs (given in
  definition \ref{ant_parts}) of the parameters in each premise and the
  conclusion are equal. Let \(\ostar\in C\) have tonicity signature
  \(\pm\) and quantification sign \(\mtt{\AE}\).

  \begin{itemize}
  \tightlist
  \item
    Each \(X_i\) in \(\vdash\ostar\) has sign
    \(\mtt{\AE}(\ostar)+\pm_i(\ostar)\) both in the premises and in the
    conclusion.
  \item
    The parameter \(U\) has sign \(1+\mtt{\AE}(\ostar)\) in
    \(\ostar\vdash\) both in the premises and in the conclusion.
  \item
    The parameters \(X_i\) have sign
    \(\delta_{i,\,n+1} + (1 + \delta_{i,\,{n+1}})\pm_i(\ostar) + \mtt{\AE}(\ostar)\)
    in \(dsr2\) both in the premises and in the conclusion.
  \item
    Now we want to check it in \(dsr1\). For any
    \(i\in\{1,\,\ldots,\,n+1\}\) and \(g\in G_D\), we need that
    \(s(S_{\mtt{\AE}}([\ostar],\,X_1,\,\ldots,\,X_{n+1}),\,X_i) = s([a_D](g,\,S_{\mtt{\AE}}([\ostar],\,X_1,\,\ldots,\,X_{n+1})),\,X_i)\).
    For it we refer to lemma \ref{sign_wf}.
  \end{itemize}
\item
  \(C5\) \textbf{Display of principal constituents}:\\
  Principal constituents only appear in \(\vdash\ostar\) and
  \(\ostar\vdash\):

  \begin{itemize}
  \tightlist
  \item
    In \(\vdash\ostar\), \(\ostar(\varphi_1,\,\ldots,\,\varphi_n)\) is
    the full antecedent or consequent by definition of \(S\).
  \item
    Also by definition of \(S\) over a connective, the application of
    \(\ostar\) over \(\varphi_1,\,\ldots,\,\varphi_n\) is the full
    antecedent or consequent. Note that none of the \(\varphi_i\) is a
    principal constituent as none of them appears \textbf{as a
    substructure} in the conclusion.
  \end{itemize}
\item
  \(C6\) and \(C7\) \textbf{Closure under substitution for
  parameters}:\\
  When substituting any of the \(X_1,\,\ldots,\,X_n\), \(X\), \(Y\) or
  \(U\) we get the same rule, as the definition of the proofs do not
  depend on the particular structure to be substituted.
\item
  \(C8\) \textbf{Eliminability of principal consituents}:\\
  For any pairs of conclusions in the rules \(X\vdash\varphi\) and
  \(\varphi\vdash Z\), with \(\varphi\) a principal constituent, either
  \(Y\) equals \(X\) or \(Z\) or we can infer \(X\vdash Z\) from the
  premises of both rules, using only the cut rule for strict subformulas
  of \(\varphi\). We check it now:

  Let's take the rules \(\vdash\ostar\) and \(\ostar\vdash\), we will
  let the sequents' order be given by the quantification sign of
  \(\ostar\). We suppose \(S_{\pm_i(\ostar)}(X_i,\,\varphi_i)\) for all
  \(i\leq n\) and
  \(S_{\mtt{\AE}}([\ostar],\,\varphi_1,\,\ldots,\,\varphi_n,\,U)\). To
  prove \(C8\), we want to infer
  \(S_{\mtt{\AE}}([\ostar],\,X_1,\,\ldots,\,X_n,\,U)\).\\
  We will take for every \(i\leq n\) an element \(g_i\in G_D\) such that
  \(g\) displays \(X_i\) in any
  \(S_s([\ostar],\,X_1,\,\ldots,\,X_{n+1})\).

  The proof follows by induction on the number of structural variables
  in\\
  \(\ostar(X_1,\,\ldots,\,X_i,\,\varphi_{i+1},\,\ldots,\,\varphi_n)\):

    {\small $$\begin{prooftree}
        \hypo{
            S_{\mtt{\AE}}([\ostar],\,\varphi_{\tau(1)},\,\ldots,\,\varphi_{\tau(n)},\,U)
        }
        \ellipsis{Induction step with $g_i$ for $i < n$}{
            S_{\mtt{\AE}}([\ostar],\,X_{1},\,\ldots,\,X_{n-1},\,\varphi_n,\,U)
        }
        \infer1[dsr1$\ostar.g_n$]{
            [a_D](g_n,\,S_{\mtt{\AE}}([\ostar],\,X_{1},\,\ldots,\,X_{n-1},\,\varphi_n,\,U))
        }
        \hypo{
            S_{-(\pm_n\cdot\mtt{\AE})([a_D](g_n,\,[\ostar]))}(X_n,\,\varphi_n)}
        \infer1[dsr2$-\pm_n$]{
            S_{\mtt{\AE}([a_D](g_n,\,[\ostar]))}(-\pm_n\varphi_n,\,-\pm_nX_n)
        }
        \infer2[cut $(\varphi_n)$]{
            [a_D](g_n,\,S_{\mtt{\AE}}([\ostar],\,X_{1},\,\ldots,\,X_{n-1},\,X_n,\,U))
        }
        \infer1[dsr1$g_n\ostar.g_n^-$]{
            S_{\mtt{\AE}}([\ostar],\,X_1,\,\ldots,\,X_n,\,U)
    }\end{prooftree}$$}

  where we have used in the cut that \([a_D](g_i,\,.)\) displays the
  \(i\)-essime component and that proposition \ref{act_seq} tells us
  that \(drs1g_n\ostar.g_n^{-}\) and \(drs1\ostar.g_n\) rules are
  inverses and also we have used the notation \(-X = *X\) and
  \(+X = X\).

  The other possible combinations are trivial.
\end{itemize}

\ndpr

\bgD

Let \(\mca{L}\) be the language of formulas in a Sequent proof system
\(\mca{P}\). The identity axiom is the following axiom:

For any \(\varphi\in\mca{L}\), \[\begin{prooftree}
  \infer0[Id]{
  \varphi\vdash\varphi}
\end{prooftree}\]

\ndD

\bgP

For any plain Connective Family \(C\), the identity axiom is admissible
in the atomic calculus \(ASL_C\).

\ndP

\begin{proof}

We will show \(\varphi\vdash\varphi\) by induction on
\(\varphi\in\mca{L}_C\):

\begin{itemize}
\tightlist
\item
  If \(\varphi = p\) with negative quantification:\\
  By \(\vdash p\) we have \([p]\vdash p\) and by \(p\vdash\) ends the
  step.
\item
  If \(\varphi = q\) with positive quantification:\\
  By \(\vdash q\) we have \(q\vdash [q]\) and by \(q\vdash\) ends the
  step.
\item
  If, for any \(n\)-ary connective \(\ostar\in C\),
  \(\varphi = \ostar(\varphi_1,\,\ldots,\,\varphi_n)\):\\
  By IH we have
  \(\forall i \in\{1,\,\ldots,\,n\},\,\varphi_i\vdash\varphi_i = S_{\pm_i(\ostar)}(\varphi_i,\,\varphi_i)\).
  Therefore, by \(\vdash\ostar\), we have
  \(S_{\mtt{\AE}}([\ostar],\,\varphi_1,\,\ldots,\,\varphi_n,\,\ostar(\varphi_1,\,\ldots,\,\varphi_n))\).
  We end the induction step by \(\ostar\vdash\).
\end{itemize}

\ndpr

\bgP

\label{vdash_invertible}

Let \(C\) be a plain connective family. The \(\ostar\vdash\) in
\(ASL_C\) is invertible. We note by \(\ostar\vdash^-\) the resulting
inference.

\ndP

\begin{proof}

\[\begin{prooftree}
    \infer0[Id]{
        S_{\pm_1(\ostar)}(\varphi_1,\,\varphi_1)
    }
    \hypo{
        \ldots
    }
    \infer0[Id]{
        S_{\pm_n(\ostar)}(\varphi_n,\,\varphi_n)
    }
    \infer3[\(\vdash\ostar\)]{
        S_{\mtt{\AE}}([\ostar],\,\varphi_1,\,\ldots,\,\varphi_n,\,\ostar(\varphi_1,\,\ldots,\,\varphi_n))   
    }
    \hypo{
        S_{\mtt{\AE}}(\ostar,\,\varphi_1,\,\ldots,\,\varphi_n,\,U)
    }
    \infer2[cut]{
        S_{\mtt{\AE}}([\ostar],\,\varphi_1,\,\ldots,\,\varphi_n,\,U)
}\end{prooftree}\]

\ndpr

\bgC

\label{struct_form_trans}

We set for strict structures the function
\(s([\ostar](X_1,\,\ldots,\,X_n)) = -\mtt{\AE}(\ostar)\).

Let \(C\) be a plain connective family. We can define a sign-dependent
translation \(\tau:[\mca{L}_C]\to\mca{L}_{\mca{O}(C)}\) such that, for
any formula \(\varphi\), \(\tau(\varphi) = \varphi\) and for any
well-formed strict structure \(X\), \(S_{-s(X)}(X,\,\tau_{s(X)}(X))\).
Therefore, from \(S_{-s(X)}(\tau(X),\,U)\) we can infer
\(S_{-s(X)}(X,\,U)\).

Furthermore, from \(S_{-s(X)}(X,\,U)\) we can infer
\(S_{-s(X)}(\tau(X),\,U)\).

\ndC

\begin{proof}

For any substructure of \(X\) of the form
\([\ostar](X_1,\,\ldots,\,X_n)\) we can check that for all
\(i\in \{1,\,\ldots,\,n\}\),
\(s(X_{i+1}) = -\mtt{\AE}(\ostar)\pm_{i+1}(\ostar)\).

We define by mutual induction two partial translations \(\tau_-\),
\(\tau_+\):

\begin{itemize}
\tightlist
\item
  For any formula \(\varphi\),
  \[\tau_{+}(\varphi) = \tau_{-}(\varphi) = \varphi\]
\item
  For any connective \(\ostar\) and structures \(X_1,\,\ldots,\,X_n\),
  \[\tau_{-\mtt{\AE}(\ostar)}([\ostar](X_1,\,\ldots,\,X_n)) = \ostar(\tau_{-\mtt{\AE}(\ostar)\pm_1(\ostar)}(X_1),\,\ldots,\,\tau_{-\mtt{\AE}(\ostar)\pm_n(\ostar)}(X_n))\]
\end{itemize}

For \(X\) a strict structure we define \(\tau(X)\) as \(\tau_{s(X)}(X)\)
and for \(X\) a formula we define \(\tau(X) = X\). This is well-defined
thanks to \(s\) being well-defined on strict structures.

By induction on \(X\) we see \(S_{-s(X)}(X,\,\tau(X))\).

\begin{itemize}
\tightlist
\item
  If \(X\) is a formula, then \(\tau_-(X) = \tau_+(X) = X\) and both
  sequents are true by \(Id\) Axiom.
\item
  We suppose \(X = [\ostar](X_1,\,\ldots,\,X_n)\). The induction
  hypothesis tells us that for all \(i\in \{1,\,\ldots,\,n\}\),
  \(S_{-s(X_{i})}(X_{i},\,\tau_{s(X_{i})}(X_{i}))\). We then take
  \(\vdash\ostar\) on all those sequents, thanks to
  \(s(X_{i+1}) = -\mtt{\AE}(\ostar)\pm_{i+1}(\ostar)\).\\
  We get a derivation of
  \begin{align*}&S_{\mtt{\AE}}([\ostar],\,X_1,\,\ldots,\,X_n,\,\ostar(\tau_{s(X_1)}(X_1),\,\ldots,\,\tau_{s(X_n)}(X_n)))\\
  &=S_{\mtt{\AE}(\ostar)}(X,\,\ostar(\tau_{-\mtt{\AE}(\ostar)\pm_1(\ostar)}(X_1),\,\ldots,\,\tau_{-\mtt{\AE}(\ostar)\pm_n(\ostar)}(X_n))) = S_{-s(\ostar)}(X,\,\tau_{s(X)}(X))\end{align*}
\end{itemize}

The last bit comes from inductively using the \(\ostar\vdash\), after
displaying all the strict substructures of \(X\).

\ndpr

\begin{quote}
It is trivial seeing that \(\tau(X)\) is a \(J\)-formula whenever \(X\)
is a \(\{[\ostar]\mid\ostar\in J\}\)-structure.

We note \(X\dashv\vdash Y\) whenever \(X\vdash Y\) and \(Y\vdash X\).
\end{quote}

\bgC

Let \(X\) and \(Y\) be two structures such that if \(X\) is strict we
can derive \(S_{s(X)}(\tau(X),\,\tau(Y))\), \(Y\) is well-formed and if
\(X\) is a formula we can derive \(X\dashv\vdash Y\). Then if \(S\) is a
derivable sequent, the sequent \(S[X/Y]\) replacing all appearances of
\(X\) in \(S\) by \(Y\) is also derivable.

\ndC

\begin{proof}

We display each appearance of \(X\) in \(S\), thanks to \(C4\) they will
all be in the same side for strict structures.

\begin{itemize}
\tightlist
\item
  If \(X\) is a formula we get, from displaying, the sequent
  \(S_s(V,\, X)\) for some sign \(s\) and structure \(V\), from cut rule
  on it and \(S_s(X,\,Y)\) we get the result after repeating the
  followed display rules in the reverse order and subtituting each use
  of \(dsr1\ostar.g\) for \(dsr1g\ostar.g^{-1}\). \(drs1\) and \(dsr2\)
  can be applied because the changed structure is a parametric
  occurrence inferences.
\item
  If \(X\) is a strict structure, by the last corollary on the displayed
  structure we can derive \(S_{s(X)}(V,\,\tau(X))\) for some structure
  \(V\). The cut rule with \(S_{s(X)}(\tau(X),\,\tau(Y))\) lets us
  derive \(S_{s(X)}(V,\,\tau(Y))\). Using again last corollary we find
  \(S_{s(X)}(V,\,Y)\). We get the result after repeating the followed
  display rules in the reverse order and subtituting each use of
  \(dsr1\ostar.g\) for \(dsr1g\ostar.g^{-1}\). \(drs1\) and \(dsr2\) can
  be applied because the changed structure is a parametric occurrence
  inferences.
\end{itemize}

\ndpr

\bgP

\label{struct_cons}

Let \(C\) be a set of connectives. The proof system \(ASL_C\) is
conservative on any fragment
\(\{[\ostar]\mid\ostar\in C\}\subseteq J\subseteq[C]\) of structural
\(C\)-connectives.

\ndP

\begin{proof}

We prove it by strong induction on the derivation of any \(J\)-structure
by eliminating the undesired usages of \(dsr1\ostar.g\), for
\([a_D](g,\,[\ostar])\notin J\). We have then to act on some of the
precedent inferences in the proof tree so that we can also derivate the
new proof tree.

Let \([\ostar]\) be the structural connective we are willing to
eliminate. We know that \([\ostar]\) is not in the conclusion of the
derivation because it is not in \(J\). As structural connectives which
are not in the connective family do not have left and right introduction
rules we know that the structure \([\ostar](X_1,\,\ldots,\,X_{n+1})\)
must be a parameter in all rules which are not the very same \(dsr1\)
introducing or eliminating the connective. Therefore, \([\ostar]\) must
have been eliminated by an application of \(dsr1\), let's say
\(dsr1\ostar.g\) for some \(g\in G_D\) such that
\([a_D](g,\,[\ostar])\in J\).

Let us note by \(\pm\), \(\sigma\) and \(\pm'\), \(\sigma'\) the
respective tonicity signature and permutation of \([\ostar]\) and
\([a_D](g,\,[\ostar])\). In case
\(T(\pm')_{n+1}+T(R(\sigma'\cdot\sigma^{-1},\,\pm))_{n+1}=0\), we
rewrite all the congruent appearences of
\(S_s([\ostar],\,X_1,\,\ldots,\,X_{n+1})\) with the conclusion of the
\(dsr1\ostar.g\) inference
\(S_s([\ostar'],\,X_{1}',\,\ldots,\,X'_{n+1})\) up from the root until
we find another congruent parameter \([\ostar](X_1,\,\ldots,\,X_{n})\)
in some inference \(dsr1\ostar''.g'\) or just leaves of the derivation.
In the first case we change \(dsr1\ostar''.g'\) by
\(dsr1\ostar''.(g\cdot g')\), who will have as conclusion
\(S_s([\ostar'],\,X_{1}',\,\ldots,\,X'_{n+1})\), the new tree of
inferences will still be a derivation because the changed constituent is
not playing any role in the inferences (it is parametric) and appears in
the same side of the sequent (so that we can still apply the same
introduction rules, if necessary). The second case is not possible as
the premises are also \(J\)-structures.

Let us note by \(\pm\), \(\sigma\) and \(\pm'\), \(\sigma'\) the
respective tonicity signature and permutation of \([\ostar]\) and
\([a_D](g,\,[\ostar])\). In case
\(T(\pm')_{n+1}+T(R(\sigma'\cdot\sigma^{-1},\,\pm))_{n+1}=1\), we note
the conclusion of the \(dsr1\ostar.g\) inference by
\(S_s([\ostar'],\,X_{1}',\,\ldots,\,*X'_{n+1})\). Now we rewrite all the
congruent appearences of \(S_s([\ostar],\,X_1,\,\ldots,\,X_{n+1})\) by
\(S_{-s}(*[\ostar'],\,X_{1}',\,\ldots,\,X'_n,\,X'_{n+1})\) up from the
root until we find another congruent parameter
\([\ostar](X_1,\,\ldots,\,X_{n})\) in some inference \(dr\ostar''.g'\)
or just leaves of the derivation. In the first case we change
\(dr\ostar''.g'\) by \(dr\ostar''.(g\cdot g')\), who will have as
conclusion \(S_{-s}(*[\ostar'],\,X_{1}',\,\ldots,\,X'_n,\,X'_{n+1})\),
the new tree of inferences will still be a derivation because the
changed constituent is not playing any role in the inferences (it is
parametric) and appears in the same side of the sequent (so that we can
still apply the same introduction rules, if necessary). The second case
is not possible as the premises are also \(J\)-structures.

This must end because each step reduces at least by one the proof
length. Furthermore, while we have some \([\ostar]\notin J\) we can
repeat this process, so that thanks to the derivation being finite we
know that the resulting derivation will not have any
\([\ostar]\notin J\).

\ndpr

\begin{quote}
We already know \(ASL_C\) to be conservative on formulas thanks to
subformula property. We can check with a quick look that for any
connective family \(C\), the \(C\)-fragment of \(ASL_{\mca{O}(C)}\)
corresponds to \(ASL_C\), as they use the same structural connectives.
\end{quote}

\bgD

Let \(C\) be a connective family and \(\ostar\) be a connective. We now
denote by \(\overline\ostar\) the class of \(\ostar\) in the equivalence
given by the orbits (and not the partition of \(C\)). We define
inductively, on formulas \(\varphi_1,\,\ldots,\,\varphi_n\),
\(V_a(\ostar(\varphi_1,\,\ldots,\,\varphi_n)) := \{\overline\ostar\}\cup V_a(\varphi_1)\cup\ldots\cup V_a(\varphi_n)\).
Let \([\ostar]\) be a structural connective. We define inductively, on
structures \(X_1,\,\ldots,\,X_n\),
\(V_a([\ostar](X_1,\,\ldots,\,X_n)) := \{\overline\ostar\}\cup V_a(X_1)\cup\ldots\cup V_a(X_n)\)
and for \(\ostar\) non-compatible with \(\beta\),
\(V_a(*[\ostar](X_1,\,\ldots,\,X_n)) = V_a([\ostar](X_1,\,\ldots,\,X_n))\).

We define the \textit{\textbf{residuated Craig Interpolation property}}
as the Craig Interpolation property, but using the predicate \(V_a\)
instead of \(V\).

\ndD

\begin{quote}
Note that \(V_a(X)\subseteq\mca{O}(C)\).

It is trivial to see that this property is more restrictive than Craig
Interpolation, by considering variables as \(0\)-ary connectives. As
Craig Interpolation refers to variables and predicates this result does
not come as a surprise, for each orbit represents a different relation
in the semantics.
\end{quote}

\bgT

Atomic Logics have the residuated Craig's Interpolation property.

\ndT

\begin{proof}

We will show that for all structures \(X\) and \(Y\), if for some sign
\(s\) the sequent \(S_{s}(X,\,Y)\) is derivable then there is a formula
\(\varphi\) such that:

\begin{itemize}
\tightlist
\item
  if \(X\) is a formula, \(S_s(X,\,\varphi)\) and \(S_s(\varphi,\,Y)\)
  are derivable and \(V(\varphi)\subseteq V(X)\cap V(Y)\).
\item
  if \(X = [\ostar](X_1,\,\ldots,\,X_n)\),
  \(S_{\mtt{\AE}(\ostar)\pm_i(\ostar)}(X_i,\,\varphi)\) and
  \(S_{\mtt{\AE}(\ostar)}([\ostar],\,X_1,\,\ldots,\,X_{i-1},\,\varphi,\,X_{i+1},\,\ldots,\,X_n,\,Y)\)
  are derivable and
  \(V(\varphi)\subseteq V(X_i)\cap (\bigcup_{j\neq i} V(X_j)\cup V(Y)\cup\{\overline\ostar\})\),
  for any \(i\in\{1,\,\ldots,\,n\}\).
\end{itemize}

We proceed by strong induction on the derivation of \(S_s(X,\,Y)\).

\begin{itemize}
\item
  If the last inference was \(\vdash\ostar\), we can assume that the
  property is satisfied for each
  \(S_{-\mtt{\AE}(\ostar)\pm_i(\ostar)}(\varphi_i,\,X_i)\) and
  \(Y = \ostar(\varphi_1,\,\ldots,\,\varphi_n)\). Then the induction
  hypothesis lets us deduce that there is some formula \(\psi_i\) such
  that \(V(\psi_i)\subseteq V(\varphi_i)\cap V(X_i)\),
  \(S_{-\mtt{\AE}(\ostar)\pm_i(\ostar)}(\varphi_i,\,\psi_i)\) and
  \(S_{-\mtt{\AE}(\ostar)\pm_i(\ostar)}(\psi_i,\,X_i)\) are derivable.
  Therefore,
  \(V(\psi_i)\subseteq V(X_i)\cap(\bigcup_{j\neq i} V(X_j)\cup V(Y))\subseteq V(X_i)\cap(\bigcup_{j\neq i} V(X_j)\cup V(Y)\cup\{\overline\ostar\})\),
  \(S_{\mtt{\AE}(\ostar)\pm_i(\ostar)}(X_i,\,\psi_i)\) and
  \(S_{\mtt{\AE}(\ostar)\pm_i(\ostar)}(\psi_i,\,\varphi_i)\), so that by
  applying \(\vdash\ostar\) we can also derive
  \(S_{\mtt{\AE}(\ostar)}([\ostar],\,X_1,\,\ldots,\,X_{i-1},\,\psi_i,\,X_{i+1},\,\ldots,\,X_n,\,Y)\).

  For the other direction of the sequent,
  \(S_{-\mtt{\AE}(\ostar)}(\ostar(\psi_1,\,\ldots,\,\psi_n),\,[\ostar](X_1,\,\ldots,\,X_n))\)
  by \(\vdash\ostar\),
  \(S_{-\mtt{\AE}(\ostar)}(\ostar(\varphi_1,\,\ldots,\,\varphi_n),\,\ostar(\psi_1,\,\ldots,\,\psi_n))\)
  by \(\vdash\ostar\) and \(\ostar\vdash\) and
  \(V(\ostar(\psi_1,\,\ldots,\,\psi_n))\subseteq\{\overline\ostar\}\cup(V(\varphi_1)\cap V(X_1))\cup\ldots\cup(V(\varphi_n)\cap V(X_n)) \subseteq V(Y)\cap (\{\overline\ostar\}\cup V(X_1)\cup\ldots\cup V(X_n))\).
\item
  If the last inference was \(\ostar\vdash\), we can assume that the
  property is satisfied for
  \(S_{\mtt{\AE}}([\ostar],\,\varphi_1,\,\ldots,\,\varphi_n,\,U)\). Then
  there is some \(\psi_i\) such that
  \(V(\psi_i)\subseteq V(\varphi_i)\cap(\bigcup_{j\neq i}V(\varphi_j)\cup V(U)\cup\{\overline\ostar\})\),
  for any sign \(s_i\), \(S_{s_i}(\varphi_i,\,\psi_i)\) and
  \(S_{\mtt{\AE}}([\ostar],\,\varphi_1,\,\ldots,\,\varphi_{i-1},\,\psi_i,\,\varphi_{i+1},\,\ldots,\,\varphi_n,\,U)\).

  Then it must be that an application of \(dr\) or an application of
  some \(\vdash\ostar'\) precedes \(\ostar\vdash\). In the second case
  we must have \(\ostar' = \ostar\) and so we can see that
  \(\ostar\in V(U)\) and the desired formula is
  \(\ostar(\psi_1,\,\ldots,\,\psi_n)\). Therefore,
  \(S_{\mtt{\AE}}([\ostar],\,\varphi_1,\,\ldots,\,\varphi_n,\,U)\) must
  be preceded by other sequents with the outmost structural connective
  still in the orbit of \([\ostar]\). As the derivation is finite, there
  must be an application of \(\vdash\ostar'\) preceding the sequent. By
  the argument we just saw, we know that
  \(\overline{[\ostar]} =\overline{[\ostar']}\), so that
  \(\overline\ostar\in V(U)\) and
  \(V(\ostar(\psi_1,\,\ldots,\,\psi_n))\subseteq V(\ostar(\varphi_1,\,\ldots,\,\varphi_n))\cap V(U)\).

  To prove it for the other direction of the sequent we must study by
  cases. If \(U\) is a formula, then we can use the same interpolant for
  the conclusion as in the premise. If \(U\) is a not a formula then we
  can use for each component of the structural connective the same
  interpolant for each component as in the premise.
\item
  If the last inference was some \(dr\) it comes trivially from doing
  case analysis on whether we are working with formulas or not and using
  the appropriate induction step from the premise.
\end{itemize}

The theorem corresponds to the case where both sides are formulas.

\ndpr

\begin{quote}
On boolean atomic logics we will have to go back to the usual Craig
Interpolation result.
\end{quote}

\hypertarget{atomic-logics-for-alphabeta}{%
\section{\texorpdfstring{Atomic Logics for
\(\alpha*\beta\)}{Atomic Logics for \textbackslash alpha*\textbackslash beta}}\label{atomic-logics-for-alphabeta}}

\label{sect6}

In this section we get back into Guillaume Aucher's action by using the
families \((\alpha*\beta)(\ostar)\).

\bgD

Let \(\ostar\) be a connective. Let \(V\) be a set of propositional
letters.

We define the connective families \((\alpha*\beta)(\ostar,\,V)\),
\((\alpha\times\varsigma)(\ostar,\,V)\),
\((\alpha\circ\delta)(\ostar,\,V)\), \(\alpha(\ostar,\,V)\),
\(\beta(\ostar,\,V)\) and \(\delta(\ostar,\,V)\) as the families
\(\mca{O}(\{\hat{\ostar_a}\}\sqcup V_a)\), where \(a\) ranges on the
previously mentioned actions.

We define the connective families \((\alpha*\beta)(\ostar)\),
\((\alpha\times\varsigma)(\ostar)\), \((\alpha\circ\delta)(\ostar)\),
\(\alpha(\ostar)\), \(\beta(\ostar)\) and \(\delta(\ostar)\) as the
families \(\mca{O}(\{\hat{\ostar_a}\}\sqcup V_\id)\), where \(a\) ranges
on the previously mentioned actions.

\ndD

\hypertarget{expressivity-of-the-different-residuation-and-negation-families}{%
\subsection{Expressivity of the Different Residuation and Negation
Families}\label{expressivity-of-the-different-residuation-and-negation-families}}

We recall theorem \ref{quot_sign}:

\setcounter{countT2}{3}

\begin{The2}

The group $(\ent\lres2\ent)^{n+1}\rtimes\mtt{Sym}(n+1)$ previously defined is isomorphic to $(\ent\lres2\ent*\mtt{Sym}(n+1))\lres(\ent\lres2\ent*\mtt{Sym}(n+1))_x$ through a morphism $\varphi$, commuting both actions ($\alpha''\circ\varphi = \alpha$).

\end{The2}

Now we conclude:

\bgC

\label{semi_free_equi}

Let \(C\) be a set of connectives. We have that
\(\mca{O}(C_{\alpha*\beta})\sim_W\mca{O}(C_{\alpha\times\varsigma})\).
Therefore the logics \(\mca{L}_{\mca{O}(C_{\alpha*\beta})}\) and
\(\mca{L}_{\mca{O}(C_{\alpha\times\varsigma})}\) are equi-expressive.

In particular
\(\forall\ostar\in C,\,(\alpha*\beta)(\ostar)\sim_W(\alpha\times\varsigma)(\ostar)\).

\ndC

\begin{proof}

On theorem \ref{quot_sign} we have already given the required morphism
on the identity bijective function, along with the equality.

This proves the result
\(\forall\ostar\in C,\,(\alpha*\beta)(\ostar)\simeq(\alpha\times\varsigma)(\ostar)\).
As we're working on discrete families, we can extend the bijection to
the sets \(\mca{O}(C_{\alpha*\beta})\) and
\(\mca{O}(C_{\alpha\times\varsigma})\) with the union of all of the
bijections while respecting the morphisms and the conditions.

\ndpr

\begin{quote}
The converse of propositions \ref{isom_equi} and \ref{bijec_equi} is
false. The following proposition shows it, as the cardinals of
\(\delta(\ostar,\,V)\) and \((\alpha\times\varsigma)(\ostar,\,V)\) are
very diferent. It also shows the need for distinguishing the families'
actions on all the connectives, including the propositional letters.
\end{quote}

\bgP

\label{dual_equi}

Let \(V\) be the propositional letters set and let \(\ostar\) be a
connective.

\begin{enumerate}
\def\labelenumi{\arabic{enumi}.}
\tightlist
\item
  For all
  \(C,\,C'\in\{(\alpha\times\varsigma)(\ostar),\,(\alpha\circ\delta)(\ostar),\,\alpha(\ostar),\,\beta(\ostar),\,\delta(\ostar),\,(\alpha\times\varsigma)(\ostar,\,V),\,(\alpha\circ\delta)(\ostar,\,V),\)
  \(\alpha(\ostar,\,V),\,\beta(\ostar,\,V),\,\delta(\ostar,\,V)\}\) we
  have \(\mca{E}_C = \mca{E}_{C'}\).
\item
  For any pair of families \(C\), \(C'\) in
  \(\{(\alpha\times\varsigma)(\ostar,\,V),\,(\alpha\circ\delta)(\ostar,\,V),\,\beta(\ostar,\,V),\,\delta(\ostar,\,V)\}\),
  \(\mca{L}_C\) and \(\mca{L}_{C'}\) are equi-expressive.
\item
  The logic over \(\mca{L}_{(\alpha\times\varsigma)(\ostar)}\) is at
  least as expressive as the logics over
  \(\mca{L}_{(\alpha\circ\delta)(\ostar)}\) and
  \(\mca{L}_{\beta(\ostar)}\). The logic over
  \(\mca{L}_{(\alpha\circ\delta)(\ostar)}\) is at least as expressive as
  the logics over \(\mca{L}_{\delta(\ostar)}\).
\item
  If we take a connective \(\ostar\) of arity \(n>0\), the logic over
  \(\mca{L}_{(\alpha\times\varsigma)(\ostar)}\) is strictly more
  expressive than the logic over
  \(\mca{L}_{(\alpha\circ\delta)(\ostar)}\).
\end{enumerate}

\ndP

\begin{proof}

Let \(n\) be the arity of \(\ostar\). Let \((k_1,\,\ldots,\,k_{n+1})\)
be the type tuple of \(\ostar\).

Firstly, we observe that the classes of models are the same. All classes
of models correspond to the class of triples \((W,\,\mca{R},\,R)\) of
world sets \(W\), relation sets \(\mca{R}\) and surjective functions
\(R\), where \(R(\ostar)\in W^{\sum_{i\in\{1,\,\ldots,\,n+1\}}t_i}\) and
the conditions \(R(a_V(p)) = R(p)\) and \(R(a(\ostar)) = R(\ostar)\),
where \(p\in V\). These are all the models with exactly a
\(\sum_{i\in\{1,\ldots,\,n+1\}}k_i\)-ary relation and, for each
\(p\in V\) of type signature \((k)\), a \(k\)-ary relation on \(W\).
Therefore, they form the same class.

We give now four translations by induction:

\begin{itemize}
\tightlist
\item
  \(\phi_1:\mca{L}_{\delta(\ostar)}\to\mca{L}_{(\alpha\times\varsigma)(\ostar)},\,\mca{L}_{\delta(\ostar,\,V)}\to\mca{L}_{(\alpha\times\varsigma)(\ostar,\,V)}\)
  such that:

  \begin{itemize}
  \tightlist
  \item
    For \(p\) propositional letter, \(\phi_1(p) = p\).
  \item
    For any \(s\in\mbb{B}\),
    \[\phi_1(\delta(s,\,\ostar)(\varphi_1,\,\ldots,\,\varphi_n)) = ((\alpha\times\varsigma)(((0,\,\ldots,\,0,\,s),\,\id),\,\ostar))(\phi_1(\varphi_1),\,\ldots,\,\phi_n(\varphi_n))\]
  \end{itemize}
\item
  We note
  \(\phi_2:\mca{L}_{\beta(\ostar)}\to\mca{L}_{(\alpha\times\varsigma)(\ostar)}\),
  \(\phi_3:\mca{L}_{\delta(\ostar)}\to\mca{L}_{(\alpha\circ\delta)(\ostar)}\)
  and
  \(\phi_4:\mca{L}_{(\alpha\circ\delta)(\ostar)}\to\mca{L}_{(\alpha\times\varsigma)(\ostar)}\).
  This translations are trivialized when using
  \((\alpha*\beta)(\ostar)\) instead of \(\alpha\times\varsigma\),
  thanks to corollary \ref{semi_free_equi}.
\end{itemize}

We need to check that \(\forall i\in\{1,\,2,\,3\}\),
\(\llbracket\phi_i(\varphi)\rrbracket = \llbracket\varphi\rrbracket\),
which are routinary computations.

We now show the reciprocal translations
\(\psi_1:\mca{L}_{(\alpha\times\varsigma)(\ostar)}\to\mca{L}_{\delta(\ostar)}\)
and \(\psi_2:\mca{L}_{\beta(\ostar)}\to\mca{L}_{\delta(\ostar)}\).

We define two translations by mutual induction \(\psi_{1,\,+}\),
\(\psi_{1,\,-}\).

\begin{itemize}
\tightlist
\item
  For any atom connective \(p\), \(\psi_{1,\,+}(p) = p\) and
  \(\psi_{1,\,-}(p) = \delta(-,\,p)\).
\item
  For any residuation of \(\ostar\), \[
  \psi_{1,\,t}(((s_1,\,\ldots,\,s_{n+1})\sigma\ostar)(\varphi_1,\,\ldots,\,\varphi_n))
  = \delta(t_{n+1},\,\ostar)(
      \psi_{1,\,t_1}(\varphi_1),\,\ldots,\,\psi_{1,\,t_n}(\varphi_n))
  \] Where
  \(\forall i\in\{1,\,\ldots,\,n\},\,t_i =\pm_i((s_1,\,\ldots,\,s_{n+1})\sigma\ostar)\cdot\pm_i(\ostar)\cdot t\).
\end{itemize}

We now show that
\(\forall t\in\mbb{B},\,\llbracket\psi_{1,\,t}(\varphi)\rrbracket_{\delta(\ostar)} = t\llbracket\varphi\rrbracket_{(\alpha\times\varsigma)(\ostar)}\)
in some model \(M\). We proceed by induction on \(\varphi\).

\begin{itemize}
\item
  For any propositional letter \(p\) it is trivial from
  \(R(\delta(s,\,p)) = R(sp)\).
\item
  For any
  \(s_1,\,\ldots,\,s_{n+1}\in\mbb{B},\,\sigma\in\mtt{Sym}(n+1)\),
  \begin{align*}
  &\llbracket\psi_{1,\,t}((s_1,\,\ldots,\,s_{n+1})\sigma\ostar(\varphi_1,\,\ldots,\,\varphi_n))\rrbracket_{\delta(\ostar)}\\
  =&\llbracket\delta(
      \pm((s_1,\,\ldots,\,s_{n+1})\sigma\ostar)\cdot\pm(\ostar)\cdot t,\,\ostar)\\
  &\phantom{\llbracket\delta}(
      \psi_{1,\,\pm_1((s_1,\,\ldots,\,s_{n+1})\sigma\ostar)\cdot\pm_1(\ostar)\cdot t}(\varphi_1),\,\ldots,\,
      \psi_{1,\,\pm_n((s_1,\,\ldots,\,s_{n+1})\sigma\ostar)\cdot\pm_n(\ostar)\cdot t}(\varphi_n))\rrbracket_{\delta(\ostar)}\\
  =&\{\overline{w}\in W^k\mid\mca{C}^{\delta(
      \pm((s_1,\,\ldots,\,s_{n+1})\sigma\ostar)\cdot\pm(\ostar)\cdot t,\,\ostar)}\\
    &\phantom{\{\overline{w}\in W^k\mid\mca{C}}(
        \llbracket\psi_{,\,1\pm_1((s_1,\,\ldots,\,s_{n+1})\sigma\ostar)\cdot\pm_1(\ostar)\cdot t}(\varphi_1)\rrbracket_{\delta(\ostar)},\,
        \ldots,\,\\
    &\phantom{\{\overline{w}\in W^k\mid\mca{C}(}
        \llbracket\psi_{1,\,\pm_n((s_1,\,\ldots,\,s_{n+1})\sigma\ostar)\cdot\pm_n(\ostar)\cdot t}(\varphi_n)\rrbracket_{\delta(\ostar)},\,
        \overline{w})\}\\
    % Using HI and the fact that relations R are the same:
  =&\{\overline{w}\in W^k\mid\mca{C}^{\varsigma(0,\,\ldots,\,0,\,
      \pm((s_1,\,\ldots,\,s_{n+1})\sigma\ostar)\cdot\pm(\ostar)\cdot t,\,\ostar)}\\
    &\phantom{\{\overline{w}\in W^k\mid\mca{C}}(
        (\pm_1((s_1,\,\ldots,\,s_{n+1})\sigma\ostar)\cdot\pm_1(\ostar)\cdot t)\llbracket\varphi_1\rrbracket_{(\alpha\times\varsigma)(\ostar)},\,\ldots,\,\\
    &\phantom{\{\overline{w}\in W^k\mid\mca{C}(}
        (\pm_n((s_1,\,\ldots,\,s_{n+1})\sigma\ostar)\cdot\pm_n(\ostar)\cdot t)\llbracket\varphi_n\rrbracket_{(\alpha\times\varsigma)(\ostar)},\,
        \overline{w})\}\\
  % Detallar molt més aquest pas
  =&\{\overline{w}\in W^k\mid\mca{C}^{(s_1,\,\ldots,\,s_{n+1})\sigma\ostar}(
        \llbracket\varphi_1\rrbracket_{(\alpha\times\varsigma)(\ostar)},\,
        \ldots,\,
        \llbracket\varphi_n\rrbracket_{(\alpha\times\varsigma)(\ostar)},\,
        \overline{w})\}\\
  =&\llbracket(s_1,\,\ldots,\,s_{n+1})\sigma\ostar(\varphi_1,\,\ldots,\,\varphi_n)\rrbracket
  \end{align*}

  The penultimate step comes from the fact that
  \(\mca{C}^{s_i(\pm,\,\ostar')}(W_1,\,\ldots,\,W_n,\,\overline{w})\) if
  and only if not
  \(\mca{C}^{\ostar'}(W_1,\,\ldots,\,\pm W_i,\,\ldots,\,W_n,\,\overline{w})\)
  and the previous ones come from \(\beta\) being an action.
\end{itemize}

We define two translations by mutual induction \(\psi_{2,\,+}\),
\(\psi_{2,\,-}\).

\begin{itemize}
\tightlist
\item
  For any atom connective \(p\), \(\psi_+(p) = p\) and
  \(\psi_-(p) = \delta(-,\,p)\).
\item
  For any residuation of \(\ostar\), \[
  \psi_t((s\ostar)(\varphi_1,\,\ldots,\,\varphi_n))
  = \delta(t_{n+1},\,\ostar)(
      \psi_{t_1}(\varphi_1),\,\ldots,\,\psi_{t_n}(\varphi_n))
  \] Where
  \(\forall i\in\{1,\,\ldots,\,n\},\,t_i = \pm_i(s\ostar)\cdot\pm_i(\ostar)\cdot t = s\cdot t\).
\end{itemize}

The proof of
\(\llbracket\psi_{2,\,s}(\varphi)\rrbracket = s\llbracket\varphi\rrbracket\)
proceeds similarly to the previous one.

We end by giving counterexamples in order to prove 4.

We will use the Proposition 5 from \emph{Generalized Keisler Theorems
for First-order Logic and Protologics} \cite{aucher_keisler}. We show
the case where \(\mtt{\AE} = \forall\).

For \(p_1,\,\ldots,\,p_n\) propositional letters, the formula
\(-\ostar(p_1,\,\ldots,\,p_n)\in\mca{L}_{(\alpha\times\varsigma)(\ostar)}\)
is not definable in \(\mca{L}_{(\alpha\circ\delta)(\ostar)}\). In order
to see it we will show that there are two bisimilar
\(\alpha\circ\delta\)-models, \(\mca{M}\) and \(\mca{N}\), such that
\(\mca{M}\vDash-\ostar(p_1,\,\ldots,\,p_n)\) and not
\(\mca{N}\vDash-\ostar(p_1,\,\ldots,\,p_n)\).

The models are:

\begin{itemize}
\tightlist
\item
  \(\mca{M}\): \(\{w_1,\,\ldots,\,w_n,\,v_1,\,\ldots,\,v_n\}\),\\
  \(R_{\overline\ostar} = \{(w_1,\,\ldots,\,w_n,\,w_1),\,\ldots,\,(w_1,\,\ldots,\,w_n,\,w_n),\,(w_1,\,\ldots,\,w_n,\,v_1),\,\ldots,\,(w_1,\,\ldots,\,w_n,\,v_n)\}\)
  and\\
  \(V(p_1) = \{w_1,\,\ldots,\,w_n,\,v_2,\,\ldots,\,v_3\},\,\ldots,\,V(p_n) = \{w_1,\,\ldots,\,w_n,\,v_1,\,\ldots,\,v_{n-1}\}\).
\item
  \(\mca{N}\): \(\{w_1',\,\ldots,\,w_n'\}\),\\
  \(R_{\overline\ostar} = \{(w_1',\,\ldots,\,w_n',\,w_1'),\,\ldots,\,(w_1',\,\ldots,\,w_n',\,w_n')\}\)
  and\\
  \(V(p_1) = \ldots = V(p_n) = \{w_1',\,\ldots,\,w_n'\}\).
\end{itemize}

The automatic bisimulation for atomic logics is
\(\{(v_1,\,w_1'),\,\ldots,\,(v_n,\,w_n'),\,(w_1,\,w_1'),\,\ldots,\,(w_n,\,w_n')\}\).

We can also check \(\mca{M}\vDash -\ostar(p_1,\,\ldots,\,p_n)\) and
\(\mca{N}\vDash \ostar(p_1,\,\ldots,\,p_n)\).

Checking it is routinary.

\ndpr

\begin{quote}
The orbit of \(\alpha\times\varsigma\) corresponds to the \(n\)-ary
connectives of fixed sign \(\pm(\ostar)\mtt{\AE}(\ostar)\). This
proposition tells us that whenever we have the negation of all the
propositional letters in the connective family (as for the boolean case
that we introduce in subsection \ref{boolean}) we do not need to work
with all the \(2^{n+1}(n+1)!\) connectives of
\((\alpha\times\varsigma)(\ostar,\,V)\) but we have enough with the
\(2(n+1)!\) connectives of \((\alpha\circ\delta)(\ostar,\,V)\) or even
with the 2 connectives of \(\delta(\ostar,\,V)\) to express the same.

We will focus on the families \((\alpha\times\varsigma)(\ostar)\),
\((\alpha\circ\delta)(\ostar)\), \(\beta(\ostar)\) and
\(\delta(\ostar)\) with the objective of keeping the maximal possible
generality, as logics do not always have negation available on all the
variables. When needed, propositional letters with non-trivial actions
will be explicitly introduced.

In what follows \(\beta(\ostar)\), \(\delta(\ostar)\) families will not
be of our interest as the display theorem \ref{asl_display} does not
apply to them. In spite of that, it is worthy noting that the
translations into \(\delta(\ostar,\,V)\) corresponds to the negation
normal form for boolean families.
\end{quote}

\hypertarget{sequent-calculus-1}{%
\subsection*{Sequent Calculus}\label{sequent-calculus-1}}

\bgE

Let \(C\) be a connective family such that \(\forall D\in\mca{A}\)
whenever \(\ostar\in D\) are not propositional letters then
\(a_D = \alpha\times\varsigma\) is satisfied. Then we can derive the
following display rules, written in Aucher's notation along with signs
in \(\mbb{B}\).

\begin{rules}
        {\underline{\textbf{The rule dsr1}}:\\
        {\small For any $[D]\in\mca{A}$ compatible with the action $\alpha\times\varsigma$ and $\forall[\ostar]\in[D]$ of arity $n$, $s_1,\,\ldots,\,s_{n+1}\in \ent\lres2\ent,\,\sigma\in\mtt{Sym}(n+1),\,X_{1},\,\ldots,\,X_{n}\in[\mca{L}_C]$,}
        $$\begin{prooftree}
            \hypo{
                S_{\mtt{\AE}}([\ostar],\,X_1,\,\ldots,\,X_{n+1})
            }
            \infer1[dsr1]{
                S_{\mtt{\AE}}((s_1,\,\ldots,\,s_{n+1})\sigma[\ostar],\,s_{n+1}s_1X_{\sigma(1)},\,\ldots,\,s_{n+1}s_nX_{\sigma(n)},\,s_{n+1}X_{\sigma(n+1)})}
        \end{prooftree}$$
        Where we write $-X_{\sigma(i)}$ for $*X_{\sigma(i)}$.\\}

        \underline{\textbf{The rule dsr2}}:\\
        {\small For any $[D]\in\mca{A}$ and $\forall[\ostar]\in[D]$ of arity $n$, $s\in \mbb{B},\,X_{1},\,\ldots,\,X_{n}\in[\mca{L}_C]$,}
        $$\begin{prooftree}
            \hypo{
                S_{s}([\ostar](X_1,\,\ldots,\,X_{n}),\,X_{n+1})
            }
            \infer1[dsr2]{
                S_{-s}(*[\ostar](X_1,\,\ldots,\,X_{n}),\,*X_{n+1})
        }\end{prooftree}$$
\end{rules}

\ndE

\begin{proof}

Follows from example \ref{act_seq_alpha_varsigma}.

\ndpr

\bgE

\label{GGL_aucher}

We can see that the preceding proof system is sound and complete with
respect the original one by using Corollary \ref{semi_free_equi}.

Let \(C\) be a connective family with action \(\alpha*\beta\) on any
non-variable connective. As the variables, which are not compatible with
\(\beta\), can't access \(dsr1\) we can use \(ASL_{C}\) to get back the
Aucher's display rules.

Below I write the original display rules, whose proof system (alongside
with the introduction rules \(\ostar\vdash\), \(\vdash\ostar\)) we call
\(GGL^0_C\).

\begin{rules}
        {\underline{\textbf{The rule dr1}}:\\
        {\small For any $[D]\in\mca{A},\,\forall[\ostar]\in[D]$ of arity $n$, $i< n+1,\,X_{1},\,\ldots,\,X_{n+1}\in[\mca{L}_C]$,}
        $$\begin{prooftree}
            \hypo{
                S_{\mtt{\AE}}([\ostar],\,X_1,\,\ldots,\,X_{n+1})
            }
            \infer1[dr1]{
                S_{\mtt{\AE}}((i\ n+1)[\ostar],\,X_{1},\,\ldots,\,X_{i-1},\,X_n,\,X_{i+1},\,\ldots,\,X_{n},\,X_{i})
        }\end{prooftree}$$\\}

        \underline{\textbf{The rule dr2}}:\\
        {\small For any $s\in \mbb{B},\,X,\,Y\in[\mca{L}_C]$,}
        $$\begin{prooftree}
            \hypo{
                S_{s}(X,\,Y)
            }
            \infer1[dr2]{
                S_{-s}(*X,\,*Y)
        }\end{prooftree}$$
\end{rules}

\begin{proof}

Let \(\ostar\in C\). Any rule \(dr1\ostar.i\) can be translated into
\(dsr1\ostar.(i\ n+1)\). Any rule \(dr2\) can be translated into
\(dsr2.1\).

We now want to translate the rule \(dsr1\) into the original Aucher's
system.

We can decompose any \(\sigma\in\mtt{Sym}(n)\) in transpositions, all
having \(n+1\) in the support, by recursively multiplying \(\sigma\)
with \((n+1\ \sigma^-(n+1))\) on the left, which will take
\(\sigma^-(n+1)\) out from the support of \(\sigma\). Then \(\sigma\)
will equal the permutations needed to reach the identity in the reverse
order.

Therefore, we can write any \((\sigma,\,b)\) as a sequence of
\(((a_1\ n+1),\,0)\cdots((a_l\ n+1),\,0)\cdot(\id,\,b)\). The rule
\(dsr1\ostar.(\sigma,\,b)\) of \(ASL_C\) will be a composition of rules
\(dr1\ostar.(a_i\ n+1)\) and, if \(b = 1\), \(dr2\) by using proposition
\ref{act_seq}.

Like this we will be able to translate the applications of the rule
\(dsr1\ostar.(\sigma_1,\,b_1)\cdots(\sigma_k,\,b_k)\). Before, let us
call the previously defined \(a_j\) for each \(\sigma_i\) as
\(a_{i,\,j}\). Indeed, the derivation is a composition of rules
\(dr1\ostar.(a_{1,\,1}\ n+1),\,\ldots,\,dr1\ostar.(a_{1,\,l_1})\), dr2
if \(b_1 = 1\),
\(dr1\ostar.(a_{2,\,1}\ n+1),\,\ldots,\,dr1\ostar.(a_{2,\,l_2})\), dr2
if \(b_2 = 1\), etc. We finish the recursion when we reach the
application of rule \(dr1\ostar.(a_{k,\,l_k})\), followed by \(dr2\) if
\(b_k = 1\).

\ndpr

\begin{figure}\begin{center}

\includegraphics{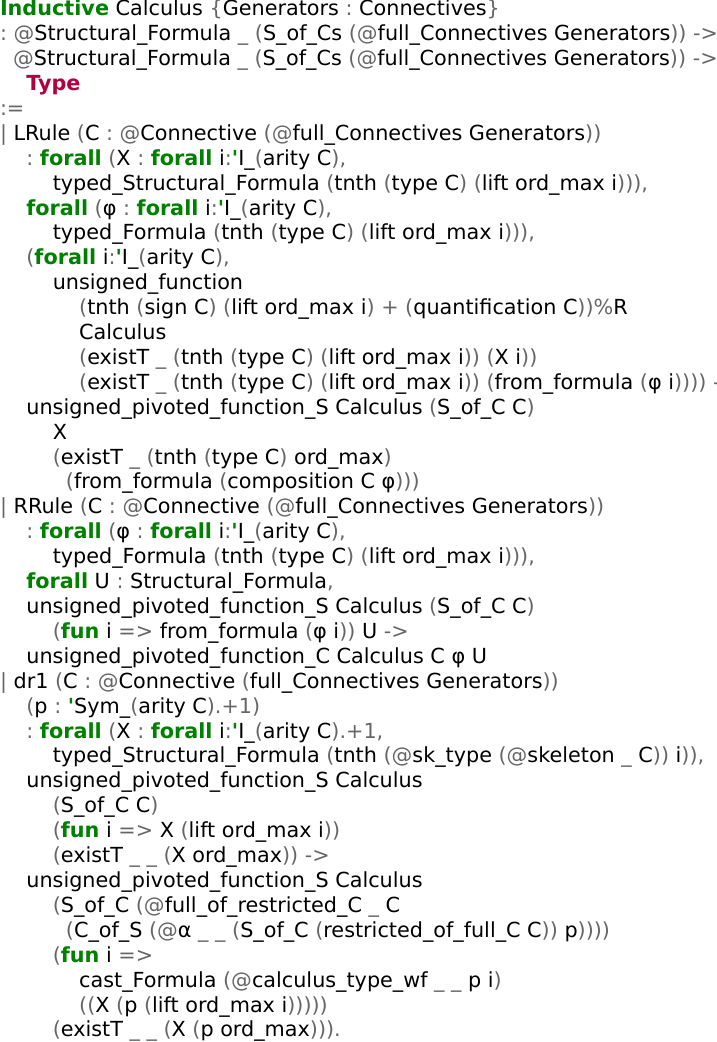}

\caption{
    Implementation of $GGL_C^0$ without $dr2$.\\
    The unsigned functions are called $S$ in the report, the \texttt{cast\_Formulas} are to cast formulas between equal types and the \texttt{calculus\_type\_wf} just checks that the type tuple of the residuation is the permutation of components.
}

\end{center}\end{figure}

\ndE

\bgT

Let \(C\) be a connective family with action \(\alpha*\beta\) on any
non-variable connective. The proof system \(GGL_C^0\) is sound and
complete with respect to the atomic logics semantics.

\ndT

\begin{proof}

Theorem 7.10 of \emph{Display and Hilbert Calculi for Atomic Logics}.
\cite{aucher_disp_hil}

\ndpr

\bgP

\label{factorial_orbits}

Let \(n\) be the arity of a connective \(\ostar\).

In the calculus \(ASL_{(\alpha*\beta)(\ostar)}\) there are, for any
\(n\)-ary connective \(\ostar\), \(2^{n+1}(n+1)\) sets of \(n!\)
connectives, where each pair of connectives in each set is provably
equivalent by permutations (\emph{i.e.} for any two \(\ostar\),
\(\ostar'\) there is some \(\sigma\in\mtt{Sym}(n)\) such that
\(\ostar(\varphi_0,\,\ldots,\,\varphi_{n_i-1})\dashv\vdash\ostar'(\varphi_{\sigma(0)},\,\ldots,\,\varphi_{\sigma(n_i-1)})\)).

Therefore, by the last corollary we can replace each appearance of
\(\ostar,\,\ostar'\) in the same set and \([\ostar],\,[\ostar']\) in any
derivable sequent to get a new derivable sequent.

\ndP

\begin{proof}

The sets correspond to the \(2^{n+1}(n+1)\) orbits of
\(\{+\sigma\in\mbb{B}*\mtt{Sym}(n+1)\mid\sigma(n) = n\}\) over the
connectives forming the \(\mca{O}(\ostar)\).

\begin{itemize}
\tightlist
\item
  For \(\otimes\), \(\sigma\otimes\) (where \(\sigma(n) = n\) and
  therefore \(\mtt{\AE}(\sigma\ostar) = \mtt{\AE}(\ostar)\)):\\
  \[\begin{prooftree}
        \infer0[Id]{
            S_{\mtt{\AE}(\ostar)}(\ostar(\varphi_0,\,\ldots,\,\varphi_{n-1}),\,\ostar(\varphi_0,\,\ldots,\,\varphi_{n-1}))
        }
        \infer1[\(\ostar\vdash^{-}\)]{
            S_{\mtt{\AE}(\ostar)}([\ostar](\varphi_0,\,\ldots,\,\varphi_{n-1}),\,\ostar(\varphi_0,\,\ldots,\,\varphi_{n-1}))
        }
        \infer1[dsr1$\ostar.\sigma$]{
            S_{\mtt{\AE}(\sigma\ostar)}([\sigma\ostar](\varphi_{\sigma(0)},\,\ldots,\,\varphi_{\sigma(n-1)}),\,\ostar(\varphi_0,\,\ldots,\,\varphi_{n-1}))
        }
        \infer1[\(\sigma\otimes\vdash\)]{
            S_{\mtt{\AE}(\ostar)}(\sigma\ostar(\varphi_{\sigma(0)},\,\ldots,\,\varphi_{\sigma(n-1)}),\,\ostar(\varphi_0,\,\ldots,\,\varphi_{n-1}))
        }
    \end{prooftree}\] The reciprocal comes from the same derivation on
  \(\sigma\ostar\vdash^{-}\), \(dsr1\vdash\ostar.\vdash^-\) and
  \(\sigma^-(\sigma\ostar)\vdash\).
\end{itemize}

\ndpr

\begin{quote}
This shows that for each coset of \(\mtt{Sym}(n+1)\) by the subgroup
\(\mtt{Sym}(n)\) we need at most a single representative, which greatly
reduces the number of connectives and display rules. At the same time,
for what we commented in the proof of the display properties, we need at
least a representative of each residuation class
\((i\ n+1)\mtt{Sym}(n)\) in the family structural connectives.
\end{quote}

\hypertarget{lambek-on-display}{%
\subsubsection*{Lambek on Display}\label{lambek-on-display}}

The rules of the Display Calculus for Lambek Calculus are:

\begin{quote}
From Rajeev Goré's \emph{How to display your favourite substructural
logic} \cite{gore_display}.\\
We have noted the connectives as follows (article's connectives in the
left):

\begin{itemize}
\tightlist
\item
  \(;\) as \([\otimes]\),
\item
  \(>\) as \([\rres]\),
\item
  \(<\) as \([\lres]\).
\end{itemize}
\end{quote}

\begin{rules}
\begin{multicols}{2}
    \underline{\textbf{Axiom}}:
    
    $$\begin{prooftree}
        \hypo{
            }\infer1[Id]{
            p\vdash p
    }\end{prooftree}
    $$

    \columnbreak
    \underline{\textbf{Residuation}}:
    $$\begin{prooftree}
        \hypo{
            X\vdash Z[\lres] Y}
        \infer[double]1[dp1]{
            X[\otimes] Y\vdash Z
        }
        \infer[double]1[dp2]{
            Y\vdash X[\rres] Z
    }\end{prooftree}$$
\end{multicols}
\hrulefill
\begin{multicols}{3}
    \underline{\textbf{For $\otimes$}}:
    $$\begin{prooftree}
        \hypo{
            X_1\vdash \varphi_1
        }\hypo{
            X_2\vdash\varphi_2
        }
        \infer2[\(\vdash\otimes\)]{
            X_1[\otimes] X_2\vdash \varphi_1\otimes\varphi_2
    }\end{prooftree}$$
    $$\begin{prooftree}
        \hypo{
            \varphi_1[\otimes] \varphi_2\vdash X
        }
        \infer1[\(\otimes\vdash\)]{
            \varphi_1\otimes\varphi_2\vdash X
    }\end{prooftree}$$

    \columnbreak
    \underline{\textbf{For $\rres$}}:
    $$\begin{prooftree}
        \hypo{
            X_1\vdash\varphi_1
        }\hypo{
            \varphi_2\vdash X_2
        }
        \infer2[\(\vdash\rres\)]{
            \varphi_1\rres\varphi_2\vdash X_1[\rres] X_2
    }\end{prooftree}$$
    $$\begin{prooftree}
        \hypo{
            X\vdash\varphi_1[\rres] \varphi_2
        }\infer1[\(\rres\vdash\)]{
            X\vdash\varphi_1\rres\varphi_2
    }\end{prooftree}$$

    \columnbreak
    \underline{\textbf{For $\lres$}}:
    $$\begin{prooftree}
        \hypo{
            \varphi_1\vdash X_1
        }\hypo{
            X_2\vdash\varphi_2
        }
        \infer2[\(\vdash\lres\)]{
            \varphi_1\lres\varphi_2\vdash X_1[\lres] X_2
    }\end{prooftree}$$
    $$\begin{prooftree}
        \hypo{
            X\vdash\varphi_1[\lres] \varphi_2
        }\infer1[\(\lres\vdash\)]{
            X\vdash\varphi_1\lres\varphi_2
    }\end{prooftree}$$
\end{multicols}
\caption*{Rules of $DL_{\mtt{Lambek}}$}
\end{rules}

\begin{figure}\begin{center}

\includegraphics{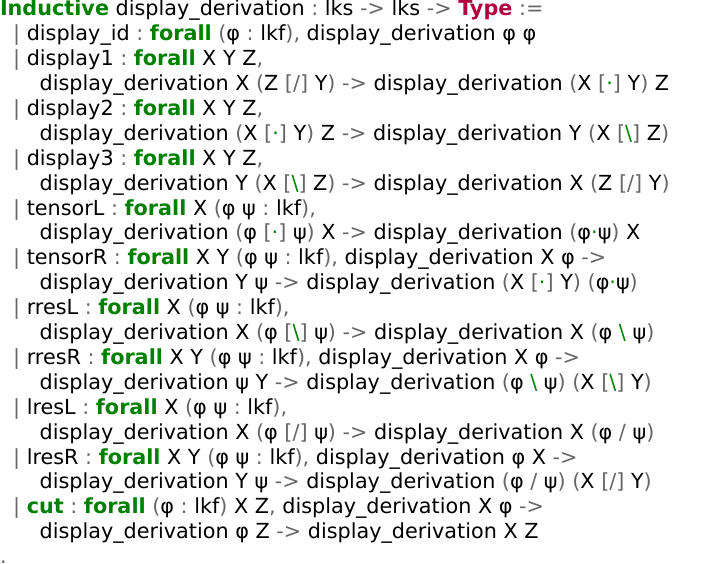}

\caption{Implementation of Goré's Non-Associative Lambek Display Calculus, we call it $DL_{\mtt{Lambek}}$}

\end{center}\end{figure}

\newpage

The rules of the calculus derived from atomic logics for
\(\otimes_{\alpha}\) are:

\begin{rules}
\small
\begin{multicols}{2}
    \underline{\textbf{For any $p = (\id,\,\pm,\,k,\,\pm)\in L$}}:\\
    $$\begin{prooftree}
        \hypo{
            [p]\vdash Z
        }
        \infer1[dsr1$p.\id$]{
            [p]\vdash Z
    }\end{prooftree}$$

    \vspace{1em}\hrulefill

    \vspace{1em}
    \underline{\textbf{For any $q = (\id,\,-\pm,\,k,\,\pm)\in L$}}:\\
    $$\begin{prooftree}
        \hypo{
            X\vdash [q]
        }
        \infer1[dsr1$q.\id$]{
            X\vdash [q]
    }\end{prooftree}$$

    \columnbreak
    \underline{\textbf{For $\otimes\in L$}}:\\
    $$\begin{prooftree}
        \hypo{
            X[\otimes] Y\vdash Z
        }
        \infer[double]1[dsr1$\otimes.\id$]{
            X[\otimes] Y\vdash Z
        }
        \infer[double]1[dsr1$\otimes.(2\ 3)$]{
            Y\vdash X[\rres] Z
        }
        \infer[double]1[dsr1$\otimes.(1\ 3)$]{
            X\vdash Z[\lres] Y
        }
        \infer[double]1[dsr1$\otimes.(1\ 2)$]{
            Y[\otimes_{\mtt{rev}}] X\vdash Z
        }
        \infer[double]1[dsr1$\otimes.(1\ 2\ 3)$]{
            Y\vdash Z[\lres_{\mtt{rev}}] X
        }
        \infer[double]1[dsr1$\otimes.(3\ 2\ 1)$]{
            X\vdash Y[\rres_{\mtt{rev}}] Z
    }\end{prooftree}$$

\end{multicols}

\caption*{Display Rules of $ASL_{\mtt{Lambek}}$}
\end{rules}

\begin{rules}
\begin{multicols*}{2}
\underline{\textbf{For any $p = (\id,\,\pm,\,k,\,\pm)\in L$}}:
\begin{multicols*}{2}
    $$\begin{prooftree}
        \infer0[\(\vdash p\)]{
            [p]\vdash p
    }\end{prooftree}$$

    \columnbreak
    $$\begin{prooftree}
        \hypo{
            [p]\vdash X
        }
        \infer1[\(p\vdash\)]{
            p\vdash X
    }\end{prooftree}$$
\end{multicols*}

\columnbreak
\underline{\textbf{For any $q = (\id,\,-\pm,\,k,\,\pm)\in L$}}:
\begin{multicols*}{2}
    $$\begin{prooftree}
        \infer0[\(\vdash q\)]{
            q\vdash [q]
    }\end{prooftree}$$

    \columnbreak
    $$\begin{prooftree}
        \hypo{
            X\vdash [q]
        }
        \infer1[\(q\vdash\)]{
            X\vdash q
    }\end{prooftree}$$
\end{multicols*}
\end{multicols*}
\hrulefill

\vspace{0.5em}
\underline{\textbf{For $\otimes\in L$}}:
    {$$\begin{prooftree}
        \hypo{
            X_1\vdash \varphi_1
        }\hypo{
            X_2\vdash\varphi_2
        }
        \infer2[\(\vdash\otimes\)]{
            X_1[\otimes] X_2\vdash \varphi_1\otimes\varphi_2
    }\end{prooftree}$$
    $$\begin{prooftree}
        \hypo{
            \varphi_1[\otimes] \varphi_2\vdash X
        }
        \infer1[\(\otimes\vdash\)]{
            \varphi_1\otimes\varphi_2\vdash X
    }\end{prooftree}$$}
\hrulefill

\begin{multicols*}{2}
    \underline{\textbf{For $\rres\in L$}}:
    $$\begin{prooftree}
        \hypo{
            X_1\vdash\varphi_1
        }\hypo{
            \varphi_2\vdash X_2
        }
        \infer2[\(\vdash\rres\)]{
            \varphi_1\rres\varphi_2\vdash X_1[\rres] X_2
    }\end{prooftree}$$
    $$\begin{prooftree}
        \hypo{
            X\vdash\varphi_1[\rres] \varphi_2
        }
        \infer1[\(\rres\vdash\)]{
            X\vdash\varphi_1\rres\varphi_2
    }\end{prooftree}$$

    \columnbreak
    \underline{\textbf{For $\lres\in L$}}:
    $$\begin{prooftree}
        \hypo{
            \varphi_1\vdash X_1
        }\hypo{
            X_2\vdash\varphi_2
        }
        \infer2[\(\vdash\lres\)]{
            \varphi_1\lres\varphi_2\vdash X_1[\lres] X_2
    }\end{prooftree}$$
    $$\begin{prooftree}
        \hypo{
            X\vdash\varphi_1[\lres] \varphi_2
        }
        \infer1[\(\lres\vdash\)]{
            X\vdash\varphi_1\lres\varphi_2
    }\end{prooftree}$$
\end{multicols*}
\caption*{Left and Right Introduction Rules of $ASL_{\mtt{Lambek}}$}
\end{rules}

\begin{quote}
As an example, we show the calculations for the rule \(dsr1.(2\ 3)\).
\[\begin{prooftree}
\hypo{
S_{\mtt{\AE}}([\otimes],\,X,\,Y,\,Z)
= S_{\exists}(X[\otimes]Y,\,Z)
= X[\otimes] Y\vdash Z
}
\infer1[dsr1$\otimes.(2\ 3)$]{
S_{\mtt{\AE}}([\rres],\,X,\,Z,\,Y)
= S_{\forall}(X[\rres] Z,\,Y)
= Y\vdash X[\rres] Z
}\end{prooftree}\]
\end{quote}

We will characterize the relationship between the three Calculi.

\bgT

\label{lambek_equi}

Let \(L\) be the Lambek Connective Family:

\begin{itemize}
\tightlist
\item
  If \(X,\,Y\in[\mca{L}_L]\) do not contain \([p]\) or \(\ast\),
  \(X\vdash_{ASL_L}Y\) in the \(\{[\ostar]\mid\ostar\in L\}\)-fragment
  of \(ASL_L\) if and only if \(X\vdash_{DL_{\mtt{Lambek}}}Y\).
\item
  If \(\varphi,\,\psi\in\mca{L}_L\), \(\varphi\vdash_{ASL_L}\psi\) if
  and only if \(\varphi\vdash_{DL_{\mtt{Lambek}}}\psi\) if and only if
  \(\varphi\vdash_{\mbb{FL}}\psi\).
\end{itemize}

\ndT

\begin{proof}

We use proposition \ref{struct_cons} to reduce the number of structural
connective appearing in the \(ASL_L\)-derivation.

We begin by showing equivalence between display calculi, ignoring
structural propositional connectives.\\
The display rules are exactly the same for the set of connectives
\(\otimes\), \(\rres\) and \(\lres\) and their restricted rule-set.
Also, the Introduction rules for binary connectives. The systems only
differ by the Axiom Id in \(DL\) and the left and right introduction
rules for the propositional variables in \(C\) in \(ASL_L\).\\
We have already shown that Id is derivable in \(ASL_L\). Now
transforming any derivation \(X\vdash_{DL_{\mtt{Lambek}}} Y\) into a
derivation \(X\vdash_{ASL_L} Y\) is trivial.

We want to see that if \(X\vdash_{ASL_L} Y\) does not contain any
\([p]\) or \([q]\), then we can get a derivation
\(X\vdash_{DL_{\mtt{Lambek}}} Y\). We proceed by induction on the
derivation. By what we commented before, the last inference must be some
\(p\vdash\), \(q\vdash\), \(dsr1\otimes\), \(\otimes\vdash\),
\(\rres\vdash\), \(\lres\vdash\), \(\vdash\otimes\), \(\vdash\rres\) or
\(\vdash\lres\). In the first two cases the premise has \([p]\) or
\([q]\) in some sequent, so that we can't apply the induction step. We
have a derivable premise \([p]\vdash X\) (where \(X\) has not any
\([p']\) nor \([q']\)) and therefore, on the derivation, they must have
been introduced by \(\vdash p\) or \(\vdash q\). In either case we
delete all \(dsr1p.\id\) and \(dsr1q.\id\) and we change the usages of
\(\vdash p\) or \(\vdash q\) for premises \([p]\vdash p\) and
\(q\vdash [q]\), respectively. Now we change all the usages of
\(p\vdash\) by premises \(p\vdash p\) and \(q\vdash q\). The new
derivation doesn't have any appearance of \([p]\) nor \([q]\) and can
therefore be translated into \(DL_{\mtt{Lambek}}\). Note how the
premises \(p\vdash p\), \(q\vdash q\) now can be changed into instances
of the \(\id\) axiom, so that we get from the derivation a proof of
\(X\vdash_{DL_{\mtt{Lambek}}}Y\).

For all \(X,\,Y\in[\mca{L}_L]\) not containing \([p]\) or \(\ast\) a
proof of \(X\vdash_{DL}Y\) implies \(\tau_-X\vdash_{\mbb{FL}}\tau_+ Y\)
and for all \(\varphi,\,\psi\in\mca{L}_L\) a proof of
\(\varphi\vdash_{\mbb{FL}}\psi\) implies \(\varphi\vdash_{DL}\psi\) are
found in the Coq code. By using Corollary \ref{struct_form_trans} and
the first point, we get the second point.

\ndpr

\begin{figure}\begin{center}

\includegraphics{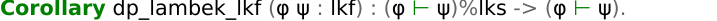}
\includegraphics{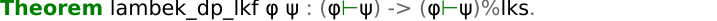}

\caption{We have proven the second point of Theorem \ref{lambek_equi} in Coq}

\end{center}\end{figure}

We can also compute its semantics.

The interpretation for Lambek Calculus' formulas given by the Atomic
Logics is defined inductively as follows:

\bgD

For the set \(W\), the relation \(R\) and model \(M = (W,\,\{R\})\):

\begin{itemize}
\tightlist
\item
  \(\llbracket \varphi\otimes\psi\rrbracket := \{w\in W\mid\exists u,\,v\in W,\,u\in\llbracket\varphi\rrbracket\wedge v\in\llbracket\psi\rrbracket\wedge R\ u\ v\ w\}\).
\item
  \(\llbracket \varphi\rres\psi\rrbracket := \{w\in W\mid\forall u,\,v\in W,\,u\notin\llbracket\varphi\rrbracket\vee v\in\llbracket\psi\rrbracket\vee R\ u\ w\ v\}\).
\item
  \(\llbracket \varphi\lres\psi\rrbracket := \{w\in W\mid\forall u,\,v\in W,\,u\in\llbracket\varphi\rrbracket\vee v\notin\llbracket\psi\rrbracket\vee R\ w\ v\ u\}\).
\end{itemize}

\ndD

\begin{quote}
This semantics are equal to the semantics given in definition
\ref{routley-meyer}.
\end{quote}

As we commented before the Bi-Lambek calculus lacks some connectives
from the full \(\alpha\times\varsigma\)-orbit. We list them below in
Aucher's syntaxis:

\begin{itemize}
\tightlist
\item
  \(\otimes_{\mtt{rev}}\):
  \(((1\ 2),\,+,\,\exists,\,1,\,(1,\,1),\,(+,\,+))\).
\item
  \(\lres_{\mtt{rev}}\):
  \(((1\ 2\ 3),\,-,\,\forall,\,1,\,(1,\,1),\,(+,\,-))\).
\item
  \(\rres_{\mtt{rev}}\):
  \(((3\ 2\ 1),\,-,\,\forall,\,1,\,(1,\,1),\,(-,\,+))\).
\item
  \(\oplus_{\mtt{rev}}\):
  \(((1\ 2),\,-,\,\forall,\,1,\,(1,\,1),\,(+,\,+))\).
\item
  \(\prec_{\mtt{rev}}\):
  \(((1\ 2\ 3),\,+,\,\exists,\,1,\,(1,\,1),\,(+,\,-))\).
\item
  \(\succ_{\mtt{rev}}\):
  \(((3\ 2\ 1),\,+,\,\exists,\,1,\,(1,\,1),\,(-,\,+))\).
\end{itemize}

\hypertarget{boolean-connectives-and-hilbert-calculus}{%
\subsection{Boolean Connectives and Hilbert
Calculus}\label{boolean-connectives-and-hilbert-calculus}}

\label{boolean}

\bgD

The Boolean Connectives set \(B\) is formed of exactly the following
non-variables connectives, for each \(k\in\nat\):

\begin{multicols}{2}
\begin{itemize}
\item $\bot_k$: $(\id,\,k,\,+,\,-)$.
\item $\top_k$: $(\id,\,k,\,-,\,-)$.
\end{itemize}
\end{multicols}
\vspace{-2em}\begin{multicols}{2}
\begin{itemize}
\item $\neg_k$: $(\id,\,k,\,(-,\,+),\,-)$.
\item $\rightarrow_k$: $((2\ 3),\,(k,\,k,\,k),\,(-,\,+,\,-),\,-)$.
\end{itemize}
\end{multicols}
\vspace{-2em}\begin{multicols}{2}
\begin{itemize}
\item $\wedge_k$: $(\id,\,(k,\,k,\,k),\,(+,\,+,\,+),\,-)$.
\item $\vee_k$: $(\id,\,(k,\,k,\,k),\,(+,\,+,\,-),\,-)$.
\end{itemize}
\end{multicols}

Furthermore we add a non-empty set of variables of output type \(k\) and
all possible tonicities.

\begin{quote}
By abuse of notation we avoid the type subindex.
\end{quote}

\begin{itemize}
\tightlist
\item
  The partition is
  \(\{\{\bot,\,\top\},\,\{\neg\},\,\{\wedge,\,\vee,\,\rightarrow\}\}\).
\item
  The actions are \(\alpha\times\varsigma\) on the connectives and
  \(\id\) on the variables.
\end{itemize}

We note \([\wedge] = [\vee] =\ ,\,\),
\(\alpha((1\ 2),\,[\wedge]) = \alpha((1\ 2),\,[\vee]) =\ ,_{\mtt{rev}}\),
\(\alpha((1\ 3),\,[\wedge]) = [\leftarrow]\),
\(\varsigma((+,\,+,\,-),\,[\leftarrow]) = [\leftarrow']\),
\(\alpha((1\ 2),\,[\leftarrow]) = [\leftarrow_{\mtt{rev}}]\),
\(\alpha\times\varsigma((+,\,+,\,-)(1\ 2),\,[\leftarrow]) = [\leftarrow'_{\mtt{rev}}]\),
\(\varsigma((+,\,+,\,-),\,[\rightarrow]) = [\rightarrow']\),
\(\alpha((1\ 2),\,[\rightarrow]) = [\rightarrow_{\mtt{rev}}]\),
\(\alpha\times\varsigma((+,\,+,\,-)(1\ 2),\,[\rightarrow]) = [\rightarrow'_{\mtt{rev}}]\)
and \([\top] = [\bot] = I\). This poses no problem as the different
structural connectives associated can be distinguished by whether they
appear in the antecedent or the consequent.

Furthermore, we define
\([\neg'] = \varsigma((+,\,+,\,-),\,[\neg]) = \ast\) and
\([\neg] = \ast\). This is well-defined because \(*[\neg] = [\neg']*\),
so that \(**X\) has a unique meaning depending on whether it appears in
the antecedent or the consequent, and
\(-\neg\varphi\dashv\vdash\varphi\) and
\(\neg'-\varphi\dashv\vdash\varphi\).

We add the following rules to \(GGL_B\):

\begin{rules}\textbf{Structural Rules}:
    \begin{multicols}{3}
        \underline{\textbf{The rule CI$\vdash$}}:\\
        {\small For any $X,\,Y,\,U\in[\mca{L}_C]$,}
        $$\begin{prooftree}
            \hypo{
                (X,\,Y)\vdash U
            }
            \infer1[CI\(\vdash\)]{
                (Y,\,X)\vdash U
        }\end{prooftree}$$

        \columnbreak
        \underline{\textbf{The rule K$\vdash$}}:\\
        {\small For any $X,\,Y,\,U\in[\mca{L}_C]$,}
        $$\begin{prooftree}
            \hypo{
                X\vdash U
            }
            \infer1[K\(\vdash\)]{
                (X,\,Y)\vdash U
        }\end{prooftree}$$
    
        \columnbreak
        \underline{\textbf{The rule WI$\vdash$}}:\\
        {\small For any $X,\,U\in[\mca{L}_C]^n$,}
        $$\begin{prooftree}
            \hypo{
                (X,\,X)\vdash U
            }
            \infer1[WI\(\vdash\)]{
                X\vdash U
        }\end{prooftree}$$
    \end{multicols}
    \begin{multicols}{2}
        \underline{\textbf{The rule IWI$\vdash$}}:\\
        {\small For any $X,\,U\in[\mca{L}_C]^n$,}
        $$\begin{prooftree}
            \hypo{
                (X,\,I)\vdash U
            }
            \infer1[IWI\(\vdash\)]{
                X\vdash U
        }\end{prooftree}$$
    
        \columnbreak
        \underline{\textbf{The rule $dr2$}}:\\
        {\small For any $X,\,Y,\,Z\in[\mca{L}_C]^n$,}
        $$\begin{prooftree}
            \hypo{
                (X,\,Y)\vdash Z
            }
            \infer[double]1[\(dr2\)]{
                X\vdash (Z,\,\ast Y)
        }\end{prooftree}$$
    \end{multicols}
\end{rules}

A Boolean Connective Family \(C\) is a Connective Family
\(C := B\sqcup C'\), for some Connective Family \(C'\), where we
consider its proof system \(GGL_C\) as the rules of \(GGL_B\) along with
the rules of \(GGL^0_{C}\) from example \ref{GGL_aucher}.

\ndD

\begin{quote}
Using CI\(\vdash\), \(dsr1\) and \(dr2\) we can show:

\begin{itemize}
\tightlist
\item
  \((X,\,Y)\vdash U\) if and only if \((X,_{\mtt{rev}}\,Y)\vdash U\).
\item
  \(U\vdash (X,\,Y)\) if and only if \(U\vdash (X,_{\mtt{rev}}\,Y)\).
\item
  \((X,\,\ast Y)\vdash U\) if and only if
  \(X[\leftarrow']_{\mtt{rev}}Y\vdash U\) if and only if
  \(Y[\rightarrow'] X\vdash U\)\\
  \(\phantom{(X,\,\ast Y)\vdash U}\) if and only if
  \(Y[\rightarrow']_{\mtt{rev}}X\vdash U\) if and only if
  \(X[\leftarrow'] Y\vdash U\).
\item
  \(U\vdash (X,\,\ast Y)\) if and only if
  \(U\vdash X[\leftarrow]_{\mtt{rev}}Y\) if and only if
  \(U\vdash Y[\rightarrow] X\)\\
  \(\phantom{U\vdash (X,\,\ast Y)}\) if and only if
  \(U\vdash Y[\rightarrow]_{\mtt{rev}}X\) if and only if
  \(U\vdash X[\leftarrow] Y\).
\end{itemize}
\end{quote}

\bgP

Let \(C\) be a Boolean Connective Family with actions on non-variable
connectives \(\alpha*\beta\). The rule \(dsr2\) is eliminable in
\(GGL_C\).

\ndP

\begin{proof}

Let \(X\) be a structure \([\ostar](X_1,\,\ldots,\,X_n)\) with
\(\overline[\ostar]\) non-compatible with \(\beta\).

\begin{center}\begin{prooftree}
  \hypo{X\vdash Y}
  \infer1[K\(\vdash\)]{(X,\,I)\vdash Y}
  \infer1[CI\(\vdash\)]{(I,\,X)\vdash Y}
  \infer1[dr2]{(I,\,*Y)\vdash *X}
  \infer1[IWI\(\vdash\)]{*Y\vdash *X}
\end{prooftree}\end{center}

\ndpr

\begin{quote}
We note \(I\vdash X\) as \(\vdash X\) and \(X\vdash I\) as \(X\vdash\).
\end{quote}

Let \(C\) be a Boolean Connective Family with action \(\alpha*\beta\) on
non-variable connectives, we will now define a new calculus for it on
formulas. The rules for Atomic Hilbert's Calculus are found in figure
10.

\begin{figure}
\begin{rules}
    \textbf{Axioms}:
    \begin{multicols}{4}
        \underline{\textbf{The axiom A0}}:\\
        {\small\phantom{AAAAAAAAAAAAAAAAAAAAAA,}}
        $$\begin{prooftree}
        \infer0[A0]{
            \vdash\top
        }\end{prooftree}$$

        \columnbreak
        \underline{\textbf{The axiom A0'}}:\\
        {\small\phantom{AAAAAAAAAAAAAAAAAAAAAA,}}
        $$\begin{prooftree}
        \infer0[A0']{
            \vdash\neg\bot
        }\end{prooftree}$$
        
        \columnbreak
        \underline{\textbf{The axiom A1}}:\\
        {\small For any $\varphi\in \mca{L}_C$,}
        $$\begin{prooftree}
        \infer0[A1]{
            \vdash\varphi\to(\varphi\wedge\varphi)
        }\end{prooftree}$$
        
        \columnbreak
        \underline{\textbf{The axiom A2}}:\\
        {\small For any $\varphi,\,\psi\in \mca{L}_C$,}
        $$\begin{prooftree}
        \infer0[A2]{
            \vdash(\varphi\wedge\psi)\rightarrow\varphi
        }\end{prooftree}$$
    \end{multicols}

    \begin{multicols}{2}
        \underline{\textbf{The axiom A3}}:\\
        {\small For any $\varphi,\,\psi,\,\rho\in \mca{L}_C$,}
        $$\begin{prooftree}
        \infer0[A3]{
            \vdash(\varphi\to\psi)\to(\neg(\psi\wedge\rho)\to\neg(\rho\wedge\varphi))
        }\end{prooftree}$$

        \columnbreak
        \underline{\textbf{The axiom A4}}:\\
        {\small For any $n\in\nat,\,\ostar\in C$ of arity $n$, $i\in \{1,\,\ldots,\,n\},\,\varphi_i\in \mca{L}_C$,}
        $$\begin{prooftree}
        \infer0[A4]{
            \vdash\neg\ostar(\varphi_1,\,\ldots,\,\varphi_n)\leftrightarrow-\ostar(\varphi_1,\,\ldots,\,\varphi_n)
        }\end{prooftree}$$
    \end{multicols}

        \underline{\textbf{The axiom A7}}:\\
        {\small For any $n\in\nat,\,\ostar\in C$ of arity $n$, $i\in \{1,\,\ldots,\,n\},\,j\in \{1,\,\ldots,\,n\},\,\varphi_i\in\mca{L}_C$ such that $\mtt{\AE}(\ostar) = \exists$,}
        $$\begin{prooftree}
        \infer0[A7]{
            \vdash\ostar(\varphi_1,\,\ldots,\,(j\ n+1)(\ostar)(\varphi_1,\,\ldots,\,\varphi_n),\,\ldots\varphi_n)\to\varphi_j
        }\end{prooftree}$$

        \underline{\textbf{The axiom A8}}:\\
        {\small For any $n\in\nat,\,\ostar\in C$ of arity $n$, $i\in \{1,\,\ldots,\,n\},\,j\in \{1,\,\ldots,\,n\},\,\varphi_i\in\mca{L}_C$ such that $\mtt{\AE}(\ostar) = \forall$,}
        $$\begin{prooftree}
        \infer0[A8]{
            \vdash\varphi_j\to\ostar(\varphi_1,\,\ldots,\,(j\ n+1)(\ostar)(\varphi_1,\,\ldots,\,\varphi_n),\,\ldots\varphi_n)
        }\end{prooftree}$$
    \caption*{Axioms of Hilbert Calculus of Atomic Logics}
\end{rules}

\begin{rules}
    \textbf{Inference Rules}:
    \begin{multicols}{2}
        \underline{\textbf{The rule R3}}:\\  
        {\small For any $n\in\nat,\,\ostar\in C$ of arity $n$, $i,\,j\in \{1,\,\ldots,\,n\},\,\varphi_i,\,\psi_j\in\mca{L}_C$ such that $\pm_j(\ostar)\mtt{\AE}(\ostar) = \exists$,}
        {\footnotesize
        $$\begin{prooftree}
          \hypo{
            \vdash\varphi_j\to\psi_j
          }
          \infer1[R3]{
            \vdash\ostar(\varphi_1,\,\ldots,\,\varphi_j,\,\ldots,\,\varphi_n)\to\ostar(\varphi_1,\,\ldots,\,\psi_j,\,\ldots,\,\varphi_n)
        }\end{prooftree}$$}

        \underline{\textbf{The rule R4}}:\\  
        {\small For any $n\in\nat,\,\ostar\in C$ of arity $n$, $i\in \{1,\,\ldots,\,n\},\,\varphi_i\in\mca{L}_C,\,j\in \{1,\,\ldots,\,n\}$ such that $\pm_j(\ostar)\mtt{\AE}(\ostar) = \forall$,}
        {\footnotesize
        $$\begin{prooftree}
          \hypo{
            \vdash\psi_j\to\varphi_j
          }
          \infer1[R4]{
            \vdash\ostar(\varphi_1,\,\ldots,\,\varphi_j,\,\ldots,\,\varphi_n)\to\ostar(\varphi_1,\,\ldots,\,\psi_j,\,\ldots,\,\varphi_n)
        }\end{prooftree}$$}
    \end{multicols}
    
        \underline{\textbf{The rule MP}}:\\  
        {\small For any $\varphi,\,\psi\in\mca{L}_C$,}
        $$\begin{prooftree}
        \hypo{
            \vdash\varphi
        }\hypo{
            \vdash\varphi\to\psi
        }
        \infer2[MP]{
            \vdash\psi
        }\end{prooftree}$$
        
    \caption*{Rules of Hilbert Calculus of Atomic Logics}
\end{rules}
\caption{Hilbert Calculus of Atomic Logics, $AHL_C$}
\label{hilbert}
\end{figure}

\begin{quote}
Some other rules can be derived from \(AHL_C\):
\end{quote}

\begin{rules}
        {\small For any $n\in\nat,\,\ostar\in C$ of arity $n$, $i\in \{1,\,\ldots,\,n\},\,\varphi_i\in\mca{L}_C,\,j\in \{1,\,\ldots,\,n\}$ such that $\mtt{\AE}(\ostar) = \forall$,}
        \begin{multicols}{2}
        {\small If $\pm_j(\ostar) = +$,}
        $$\begin{prooftree}
          \hypo{
            \vdash\varphi_j
          }
          \infer1[R1]{
            \vdash\ostar(\varphi_1,\,\ldots,\,\varphi_j,\,\ldots,\,\varphi_n)
        }\end{prooftree}$$

        \columnbreak
        {\small If $\pm_j(\ostar) = -$,}
        $$\begin{prooftree}
          \hypo{
            \vdash\neg\varphi_j
          }
          \infer1[R2]{
            \vdash\ostar(\varphi_1,\,\ldots,\,\varphi_j,\,\ldots,\,\varphi_n)
        }\end{prooftree}$$
        \end{multicols}

        {\small For any $n\in\nat,\,\ostar\in C$ of arity $n$, $i\in \{1,\,\ldots,\,n\},\,j\in \{1,\,\ldots,\,n\},\,\varphi_i,\,\varphi_j'\in\mca{L}_C$ such that $\mtt{\AE}(\ostar) = \exists$ and $\pm(\ostar)_j = +$,}
        $$\begin{prooftree}
        \infer0[A5]{
            \vdash\ostar(\varphi_1,\,\ldots,\,\varphi_j\vee\varphi_j',\,\ldots,\,\varphi_n)\to(\ostar(\varphi_1,\,\ldots,\,\varphi_j,\,\ldots,\,\varphi_n)\vee\ostar(\varphi_1,\,\ldots,\,\varphi'_j,\,\ldots,\,\varphi_n))
        }\end{prooftree}$$

        {\small For any $n\in\nat,\,\ostar\in C$ of arity $n$, $i\in \{1,\,\ldots,\,n\},\,j\in \{1,\,\ldots,\,n\},\,\varphi_i,\,\varphi_j'\in\mca{L}_C$ such that $\mtt{\AE}(\ostar) = \exists$ and $\pm(\ostar)_j = -$,}
        $$\begin{prooftree}
        \infer0[A6]{
            \vdash\ostar(\varphi_1,\,\ldots,\,\varphi_j\wedge\varphi_j',\,\ldots,\,\varphi_n)\to(\ostar(\varphi_1,\,\ldots,\,\varphi_j,\,\ldots,\,\varphi_n)\vee\ostar(\varphi_1,\,\ldots,\,\varphi'_j,\,\ldots,\,\varphi_n))
        }\end{prooftree}$$
\end{rules}

\bgT

\label{hilbert_seq}

Let \(C\) be a Boolean Connective Family. Let \(X,\,Y\) be
\(\{[\ostar]\mid\ostar\in C\}\)-structures, the sequent \(X\vdash Y\) is
derivable in \(GGL_C\) if and only if the formula
\(\vdash \tau(X\vdash Y)\) is derivable in \(AHL_C\), where
\(\tau(X\vdash Y) = \tau_-(X)\rightarrow\tau_+(Y)\).

\ndT

\begin{proof}

Let \(X\) and \(Y\) be structures. We begin by taking a derivation in
the \(\{[\ostar]\mid\ostar\in C\}\)-fragment of \(ASL_C\) of
\(X\vdash_{ASL_C}Y\) and proving \(\vdash\tau(X\vdash Y)\) in \(AHL_C\),
we will proceed by induction on the derivation.

\begin{quote}
The usages of \(dr1\) refer to the rule \(dsr1\) for \(\alpha\), in
accordance to what was observed in example \ref{GGL_aucher}, and the
usages of \(dr2\) refer to the one presented in \(GGL_C\).
\end{quote}

\begin{itemize}
\item
  If the last inference of the derivation was some of the structural
  rules, the derivation corresponds to a propositional logic lemma.
\item
  If the last inference of the derivation was dr2 and we have structures
  \(X\), \(Y\) and \(Z\) such that by induction hypothesis we can
  derivate \(\vdash\tau((X,\,Y)\vdash Z)\), this is
  \(\vdash\tau_-(X)\wedge\tau_-(Y)\to\tau_+(Z)\), then by a
  propositional logics' lemma we get
  \(\vdash\tau_-(X)\to\tau_+(Z)\vee\neg\tau_-(Y)\). Using \(A4\) and
  doing case analysis on \(Y\) we know that
  \(\vdash\tau_+(*Y)\to\neg\tau_-(Y)\) so that by \(R4\) we get
  \(\vdash(\tau_+(Z)\vee\neg\tau_-(Y))\to\tau_+(Z)\vee\tau_+(*Y)\).
  Transitivity of \(\rightarrow\) lets us finish with
  \(\vdash\tau(X\vdash(Z,\,*Y))\).
\item
  If the last inference of the derivation was dsr1 and
  \(\mtt{\AE}(\ostar) = \exists\), then by induction hypothesis we have
  \(\vdash\ostar(\tau_{-\pm_1(\ostar)}(X_1),\,\ldots,\,\tau_{-\pm_1(\ostar)}(X_1))\to\tau_+(X_{n+1})\).
  We will prove
  \(\vdash\tau(\alpha*\beta((b_1,\,\sigma_1)\cdots(b_l,\,\sigma_l),\,[\ostar](X_1,\,\ldots,\,X_{n+1})))\)
  in the Hilbert System. It suffices to show it true \(\alpha\) for
  \((i\ n+1)\), as consecutive applications of the rule \(dsr1\) over
  \(\sigma\) and \(\sigma'\) are the same as the rule \(dsr1\) over
  \(\sigma'\circ\sigma\) and \(\beta\) for \(-\) for the same reason. If
  we take \(\sigma\in\mtt{Sym}(n)\) and we decompose it in
  transpositions, all having \(n+1\) in the support (the algorithm
  consists on recursively multiplying \(\sigma\) by
  \((n+1\ \sigma^-(n+1))\) on the left, which will take
  \(\sigma^-(n+1)\) out from the support of \(\sigma\), \(\sigma\) will
  equal the permutations needed to reach the identity in the reverse
  order), we can extend the inference to any permutation. Composing both
  the \(\alpha\) and \(\beta\) cases we get the derivation for
  \(\alpha*\beta\).

  Let \(j\in \{1,\,\ldots,\,n\}\) and we further assume
  \(\pm_j(\ostar) = +\), by induction hypothesis we have
  \begin{equation*}
  \vdash\ostar(\tau_{-\pm_1(\ostar)}(X_1),\,\ldots,\,\tau_{-\pm_n(\ostar)}(X_n))\to \tau_+(X_{n+1})
  \end{equation*} Therefore, by R3, we have \begin{multline}\label{eq:1}
    \vdash(j\ n+1)\ostar(\tau_{-\pm_1(\ostar)}(X_1),\,\ldots,\,\ostar(\tau_{-\pm_1(\ostar)}(X_1),\,\ldots,\,\tau_{-\pm_1(\ostar)}(X_1)),\,\ldots,\,\tau_{-\pm_n(\ostar)}(X_n))\to\\
    (j\ n+1)\ostar(\tau_{-\pm_1(\ostar)}(X_1),\,\ldots,\,\tau_+(X_{n+1}),\,\ldots,\,\tau_{-\pm_n(\ostar)}(X_n))
  \end{multline} We conclude by transitivity of \(\rightarrow\) on the
  axiom \(A8\) \begin{multline*}
  \vdash\tau_{-\pm_j(\ostar)}(X_j)\to\\
  (j\ n+1)\ostar(\tau_{-\pm_1(\ostar)}(X_1),\,\ldots,\,\ostar(\tau_{-\pm_1(\ostar)}(X_1),\,\ldots,\,\tau_{-\pm_1(\ostar)}(X_1)),\,\ldots,\,\tau_{-\pm_n(\ostar)}(X_n))
  \end{multline*} and the sequent on derivation \ref{eq:1}.

  Whenever \(\pm_j(\ostar) = -\) we repeat the same procedure, but using
  \(A7\) and \(R4\). Whenever \(\mtt{\AE}(\ostar) = \forall\) we repeat
  the same procedure but using, respectively, \(A7\) and \(R3\) and
  \(A8\) and \(R4\).
\item
  We will prove now \(\vdash\ostar\), assuming that we can translate all
  sequents from the hypothesis. We assume
  \(\mtt{\AE}(\ostar) = \exists\). We proceed by induction on \(k\) over
  \begin{equation*}
    \vdash \ostar(\tau_{-\pm_1}(X_1),\,\ldots,\,\tau_{-\pm_n}(X_n))\to\ostar(\varphi_1,\,\ldots,\,\varphi_k,\,\tau_{-\pm_{k+1}}(X_{k+1}),\,\ldots,\,\tau_{-\pm_n}(X_n))
  \end{equation*} The base case \(k = 0\) is the propositional logic
  result \(\vdash A\to A\). The induction step consists on using
  transitivity of \(\rightarrow\) after:

  \begin{itemize}
  \tightlist
  \item
    when \(\pm_{k+1} = +\), \(R3\) for the first induction hypothesis on
    \(X_{k+1}\vdash\varphi_{k+1}\)
    (\(\vdash \tau_{+}(X_{k+1})\to\varphi_{k+1}\)).
  \item
    when \(\pm_{k+1} = -\), \(R4\) for the first induction hypothesis on
    \(\varphi_{k+1}\vdash X_{k+1}\)
    (\(\vdash\varphi_{k+1}\to\tau_{-}(X_{k+1})\)).
  \end{itemize}

  When \(\mtt{\AE}(\ostar) = \forall\) we inverse the applications of
  \(R3\) and \(R4\).

  The proof for the action of \(\beta\) corresponds to eliminating
  \(dr2\).
\item
  Applying induction hypothesis to the premises of \(\ostar\vdash\) we
  can trivially infer the conclusions, as they are the same sequents.
\end{itemize}

We now want to prove that if we can derive \(\vdash \tau(X\vdash Y)\) in
\(AHL_C\) then we can derive \(\vdash \tau(X\vdash Y)\) in \(ASL_C\). We
proceed by induction on the Hilbert Calculi derivation.

\begin{quote}
We want to highlight the following derivations
  \begin{multicols}{2}
    $$\begin{prooftree}
    \hypo{\varphi\vdash\psi}
    \infer1[K\(\vdash\)]{\varphi,\,I\vdash\psi}
    \infer1[CI\(\vdash\)]{I,\,\varphi\vdash\psi}
    \infer1[dsr2\(\wedge.(2\ 3)\)]{I\vdash\varphi[\rightarrow]\psi}
    \infer1[\(\rightarrow\vdash\)]{I\vdash\varphi\rightarrow\psi}
    \end{prooftree}$$

    \columnbreak
    $$\begin{prooftree}
    \hypo{I\vdash\varphi\rightarrow\psi}
    \infer1[\(\rightarrow\vdash^-\)]{I\vdash\varphi[\rightarrow]\psi}
    \infer1[dsr2\(\rightarrow.(2\ 3)\)]{I,\,\varphi\vdash\psi}
    \infer1[CI\(\vdash\)]{\varphi,\,I\vdash\psi}
    \infer1[IWI\(\vdash\)]{\varphi\vdash\psi}
    \end{prooftree}$$
  \end{multicols}
Therefore $\varphi\vdash_{ASL_L}\psi$ can be derived if and only if $\vdash_{ASL_L}\varphi\rightarrow\psi$ can be derived.
\end{quote}

\begin{itemize}
\item
  The inference rules which only use booleans connectives are lemmas
  from propositional logic.
\item
  We want to prove axiom \(A4\). First we prove
  \(\vdash\neg\ostar(\varphi_1,\,\ldots,\,\varphi_n)\to-\ostar(\varphi_1,\,\ldots,\,\varphi_n)\)
  and then
  \(\vdash-\ostar(\varphi_1,\,\ldots,\,\varphi_n)\to\neg\ostar(\varphi_1,\,\ldots,\,\varphi_n)\).
  With both sequents, by using \(\vdash\wedge\) and \(IWI\vdash\) we
  derive \(A4\).

  We suppose first \(\mtt{\AE}(\ostar) = \exists\). To prove
  \(\vdash\neg\ostar(\varphi_1,\,\ldots,\,\varphi_n)\to-\ostar(\varphi_1,\,\ldots,\,\varphi_n)\)
  we begin with the identities
  \(\varphi_1\vdash\varphi_1,\,\ldots,\,\varphi_n\vdash\varphi_n\) and
  applying \(\vdash\ostar\), we infer
  \([\ostar](\varphi_1,\,\ldots,\,\varphi_n)\vdash\ostar(\varphi_1,\,\ldots,\,\varphi_n)\).
  Now by dsr2 we have
  \(\ast \ostar(\varphi_1,\,\ldots,\,\varphi_n)\vdash\ast [\ostar](\varphi_1,\,\ldots,\,\varphi_n)\).
  We then use \(\neg\vdash\) on the left and \(-\ostar\vdash\) on the
  right, to find
  \(\neg\ostar(\varphi_1,\,\ldots,\,\varphi_n)\vdash-\ostar(\varphi_1,\,\ldots,\,\varphi_n)\).
  We end with K\(\vdash\), \(dsr1\wedge\) and \(\rightarrow\vdash\).

  Now to prove
  \(\vdash-\ostar(\varphi_1,\,\ldots,\,\varphi_n)\to\neg\ostar(\varphi_1,\,\ldots,\,\varphi_n)\)
  we begin with the identity
  \(-\ostar(\varphi_1,\,\ldots,\,\varphi_n)\vdash-\ostar(\varphi_1,\,\ldots,\,\varphi_n)\).
  Using the inverse inference of \(-\ostar\vdash\) (Proposition
  \ref{vdash_invertible}) and then \(\neg\vdash\) we deduce
  \(-\ostar(\varphi_1,\,\ldots,\,\varphi_n)\vdash\neg\ostar(\varphi_1,\,\ldots,\,\varphi_n)\).
  We end with K\(\vdash\), \(dsr1\wedge\) and \(\rightarrow\vdash\).

  If \(\mtt{\AE}(\ostar) = \forall\) we invert the proof's order.
\item
  Now the proof of A7.\\
  It corresponds to dr1 after applying the inverse of
  \((j\ n+1)\ostar\vdash\) on the right of
  \((j\ n+1)\ostar(\varphi_1,\,\ldots,\,\varphi_n)\vdash(j\ n+1)\ostar(\varphi_1,\,\ldots,\,\varphi_n)\).
  We end by using dr2.
\item
  The proof of A8. It corresponds to dr1 after applying the inverse of
  \((j\ n+1)\ostar\vdash\) on the left of
  \((j\ n+1)\ostar(\varphi_1,\,\ldots,\,\varphi_n)\vdash(j\ n+1)\ostar(\varphi_1,\,\ldots,\,\varphi_n)\).
  We end by using dr2.
\item
  The proof of R3.\\
  It corresponds to the inverse of \(\rightarrow\vdash\), dr2,
  \(\vdash\ostar\), dr2 and \(\rightarrow\vdash\) in that order.
\item
  Finally, the proof of R4.\\
  It corresponds to the inverse of \(\rightarrow\vdash\), dr2,
  \(\vdash\ostar\), dr2 and \(\rightarrow\vdash\) in that order.
\end{itemize}

\ndpr

\begin{quote}
We now write the proof of \(R1\) and \(A5\), which will also illustrate
how \(R2\) and \(A6\) are proven.

\underline{Proof of $R1$}:\\
We have \(\vdash\varphi_j\), by using \(MP\) on axiom \(A8\)\\
we get
\(\vdash\ostar(\varphi_1,\,\ldots,\,(j\ n+1)\ostar(\varphi_1,\,\ldots,\,\varphi_n),\,\ldots,\,\varphi_n)\).\\
Furthermore, by Rule \(R3\) on the propositional logic's result
\(\vdash(j\ n+1)\ostar(\varphi_1,\,\ldots,\,\varphi_n)\to\top\)\\
we get
\(\vdash\ostar(\varphi_1,\,\ldots,\,(j\ n+1)\ostar(\varphi_1,\,\ldots,\,\varphi_n),\,\ldots,\,\varphi_n)\to\ostar(\varphi_1,\,\ldots,\,\top,\,\ldots,\,\varphi_n)\).\\
The rule \(MP\) gets us
\(\vdash\ostar(\varphi_1,\,\ldots,\,\top,\,\ldots,\,\varphi_n)\).\\
We use rule \(R3\) on \(\vdash\top\to\varphi_j\) to get
\(\vdash\ostar(\varphi_1,\,\ldots,\,\top,\,\ldots,\,\varphi_n)\to\ostar(\varphi_1,\,\ldots,\,\varphi_j,\,\ldots,\,\varphi_n)\).\\
We end by using again the rule \(MP\).

\underline{Proof of $A5$}:\\
We have by axiom \(A8\) that
\(\vdash\varphi_j\to(j\ n+1)\ostar(\varphi_1,\,\ldots,\,\ostar(\varphi_1,\,\ldots,\,\varphi_j,\,\ldots,\,\varphi_n),\,\ldots,\,\varphi_n)\)\\
and
\(\vdash\varphi_j'\to(j\ n+1)\ostar(\varphi_1,\,\ldots,\,\ostar(\varphi_1,\,\ldots,\,\varphi_j',\,\ldots,\,\varphi_n),\,\ldots,\,\varphi_n)\).\\
We get by \(R3\) on the propositional logics results\\
\(\vdash\ostar(\varphi_1,\,\ldots,\,\varphi_j,\,\ldots,\,\varphi_n)\to\ostar(\varphi_1,\,\ldots,\,\varphi_j,\,\ldots,\,\varphi_n)\vee\ostar(\varphi_1,\,\ldots,\,\varphi_j',\,\ldots,\,\varphi_n)\)\\
and
\(\vdash\ostar(\varphi_1,\,\ldots,\,\varphi_j',\,\ldots,\,\varphi_n)\to\ostar(\varphi_1,\,\ldots,\,\varphi_j,\,\ldots,\,\varphi_n)\vee\ostar(\varphi_1,\,\ldots,\,\varphi_j',\,\ldots,\,\varphi_n)\)
that \begin{align*}
\vdash
(j\ n+1)\ostar&(
\varphi_1,\,
\ldots,\,
\ostar(\varphi_1,\,\ldots,\,\varphi_j,\,\ldots,\,\varphi_n),\,
\ldots,\,
\varphi_n)\\
&\to(j\ n+1)\ostar(
\varphi_1,\,
\ldots,\,
\ostar(
\varphi_1,\,
\ldots,\,
\varphi_j,\,
\ldots,\,
\varphi_n)
\vee\ostar(
\varphi_1,\,
\ldots,\,
\varphi_j',\,
\ldots,\,
\varphi_n),\,
\ldots,\,
\varphi_n)
\end{align*} and \begin{align*}
\vdash
(j\ n+1)\ostar&(
\varphi_1,\,
\ldots,\,
\ostar(\varphi_1,\,\ldots,\,\varphi_j',\,\ldots,\,\varphi_n),\,
\ldots,\,
\varphi_n)\\
&\to(j\ n+1)\ostar(
\varphi_1,\,
\ldots,\,
\ostar(
\varphi_1,\,
\ldots,\,
\varphi_j,\,
\ldots,\,
\varphi_n)
\vee\ostar(
\varphi_1,\,
\ldots,\,
\varphi_j',\,
\ldots,\,
\varphi_n),\,
\ldots,\,
\varphi_n)
\end{align*} By transitivity of \(\rightarrow\) we get
\[\vdash\varphi_1\to(j\ n+1)\ostar(
\varphi_1,\,
\ldots,\,
\ostar(
\varphi_1,\,
\ldots,\,
\varphi_j,\,
\ldots,\,
\varphi_n)
\vee\ostar(
\varphi_1,\,
\ldots,\,
\varphi_j',\,
\ldots,\,
\varphi_n),\,
\ldots,\,
\varphi_n)\] and \[\vdash\varphi_1'\to(j\ n+1)\ostar(
\varphi_1,\,
\ldots,\,
\ostar(\varphi_1,\,\ldots,\,\varphi_j,\,\ldots,\,\varphi_n)
\vee\ostar(\varphi_1,\,\ldots,\,\varphi_j',\,\ldots,\,\varphi_n),\,
\ldots,\,
\varphi_n)\] Propositional logic results give us
\[\vdash(\varphi_j\vee\varphi_j')\to
(j\ n+1)\ostar(
\varphi_1,\,
\ldots,\,
\ostar(\varphi_1,\,\ldots,\,\varphi_j,\,\ldots,\,\varphi_n)
\vee\ostar(\varphi_1,\,\ldots,\,\varphi_j',\,\ldots,\,\varphi_n),\,
\ldots,\,
\varphi_n)\] By rule \(R3\), \begin{align*}\footnotesize
\vdash\ostar(\varphi_1,\,&\ldots,\,(\varphi_j\vee\varphi_j'),\,\ldots,
\,\varphi_n)\\\footnotesize
\to(\ostar(&\varphi_1,\,\ldots,\,\\\footnotesize
&(j\ n+1)\ostar(
\varphi_1,\,
\ldots,\,\ostar(\varphi_1,\,\ldots,\,\varphi_j,\,\ldots,\,\varphi_n)
\vee\ostar(\varphi_1,\,\ldots,\,\varphi_j',\,\ldots,\,\varphi_n),\,
\ldots,\,
\varphi_n),\,
\ldots,\,\varphi_n))
\end{align*} Axiom \(A7\) tells us \begin{align*}\footnotesize
\vdash&\ostar(\varphi_1,\,\ldots,\,(j\ n+1)\ostar(
\varphi_1,\,
\ldots,\,
\ostar(\varphi_1,\,\ldots,\,\varphi_j,\,\ldots,\,\varphi_n)
\vee\ostar(\varphi_1,\,\ldots,\,\varphi_j',\,\ldots,\,\varphi_n),\,
\ldots,\,
\varphi_n),\,\ldots,\,\varphi_n)\\\footnotesize
&\to
(\ostar(\varphi_1,\,\ldots,\,\varphi_j,\,\ldots,\,\varphi_n)
\vee\ostar(\varphi_1,\,\ldots,\,\varphi_j',\,\ldots,\,\varphi_n))
\end{align*} We finish the derivation by transitivity of
\(\rightarrow\).
\end{quote}

\newpage

\hypertarget{conclusions-and-future-work}{%
\section{Conclusions and Future
Work}\label{conclusions-and-future-work}}

The context in which this report is written is the discovery by Aucher
that residuation can be explained with two group actions. In this report
we see that using a class of actions we can still prove the display
theorem for their Atomic Logics, reinforcing Aucher's findings. Through
this framework left and right introduction rules are already given
(although it may sometimes be necessary to work with a composition of
Atomic Logics rules, called then Molecular Logics), reducing the problem
of defining a logics into extending the proof system with structural
rules in an approach satisfying Do\v sen principle. Some of the problems
we found in the process of formalizing Atomic Logics were related to the
lack of a precise definition for connective families, which we just
introduced.

In this report we have provided a finitary version of Aucher's action by
using the semi-direct product of the actions
\(\alpha:\mtt{Sym}(n+1)\to\mtt{GL}(\ent\lres2\ent)^{n+1}\) and
\(\varsigma:(\ent\lres2\ent)^{n+1}\to\mtt{GL}(\ent\lres2\ent)^{n+1}\).
We have also given an explicit formulation for his version of the
action. We constructed a new framework where we can compare different
possible versions of actions on connective families. This includes
definitions for connectives and structures by algebraic means. In
proposition \ref{act_seq} we have given a result which makes it possible
to construct a new sequent calculus and prove the display theorem
\ref{asl_display} in arbitrary connective families (not necessarily
using the \(\alpha*\beta\) action). Other results in the report are the
proof of redundancy of 4 Hilbert System Rules, the proof of equivalence
between the Hilbert and the Sequent Systems for Atomic Logics, the proof
of Craig Interpolation on Atomic Logics and the proof of proposition
\ref{dual_equi}. This last proposition gives a contraexample to the
reciprocal of proposition \ref{bijec_equi} and it shows the need for
explicitly giving the partitions on propositional letters. The author
personal interpretation on this report is that an algebraic
interpretation of Atomic Logics properties has made possible a better
understanding on the matter, and provided necessary changes on the
theory for its implementation on a proof-assistant.

There are some topics this report has not covered. Namely, the
construction of a connective family with action on
\((\ent\lres2\ent)^{n+1}\times\ent\lres(n+1)\ent\), whose logic is
equi-expressive to the one given by \(\alpha\times\varsigma\) by using
proposition \ref{factorial_orbits}, greatly reducing so the complexity
of the resulting Atomic Logics, a study on how to fully embed boolean
negation structural connective \(*\) inside connective families, a study
on how to properly add structural rules to the proof systems and
properly formalize a wider range of logics with it. Possibly this last
point will be best covered in a proof-checker software, where we might
profit from formalizing on it proposition \ref{act_seq} before
establishing a choice for the actions we are interested in for
Connective Families. Specially relevant would be to have a survey on
those properties in proof systems we are interested in proving within
the general context of fragments of first-order logics. Some open
questions are whether it is possible to have connective families more
expressive than those of action \(\alpha*\beta\), how does adding
structural connectives change a connective families and sufficient and
necessary conditions for completeness of the newly presented \(ASL_C\)
Calculi with regard to the Kripke semantics.

\printbibliography[heading=bibintoc, title={References}]

\hypertarget{some-comments-on-the-internship}{%
\section*{Some Comments on the
Internship}\label{some-comments-on-the-internship}}
\addcontentsline{toc}{section}{Some Comments on the Internship}

The Master Thesis research has mainly developed in two halves. On the
first half we planned to write down in Coq the proof system and
semantics of Atomic Logics and then going on proving the display
theorem, consistency with regard to the semantics and some results on
those lines. For it I chose the \emph{mathcomp} library as an
appropriate context for introducing Atomic Logics, both for its usage of
group theory and for the common goal of keeping everything computable.
At this point the formalization showed several difficulties, which did
not seem natural in the problem we were working on. Indeed, the action
presented in Guillaume Aucher's articles required a lot of dependent
typing for its definition and furthermore it used a free product of
groups, which is an infinite groups. This second point was a huge
problem, as mathcomp does not work with infinite groups. For this reason
I tried to work separately with the display calculus using only the rule
\(dr1\), resulting on the
\href{https://github.com/The-busy-man/Universal-FO-Logics/commit/fe6a68645030adee9f51bae8d10e212f9300004b}{\textit{proof of action} commit},
where I provide a version of the proof system without dr2. By this point
I had already found a function which did the same as Guillaume Action
but without having to resort to its definition as an action and I was
planning to go on into adding \(dr2\) without even having to use the
free product.

At that moment another problem showed up: I had been using full orbits
and representing structures as copies of those full orbits. To have the
most flexibility in the code we required for the connective set the
possibility of being strict subsets of the full orbit. This point was
quite worrying, as it was not clear to me what were the connective
families. On this second half I begun to work on a precise formulation
of connective families and structures. During this process I realised
that the residuation action \(\alpha\) in Atomic Logics was acting as a
linear transformation on tonicity signatures. Guillaume pointed out to
me that he had already expected some parallelism with geometry in his
first article on the topic\cite{aucher_gaggle}, so that this observation
stays in line with the expectations. This geometrical interpretation
turned out to be a useful insight, as it led me to use the semi-direct
product of groups for providing a new action giving the same orbits to
Aucher's while being defined on a finite group.



\section*{Supplementary material}

# Appendix: $\alpha$ is an action

In this appendix the theorem sk_$\alpha$_is_action found in is presented. This constitutes a proof of Proposition 5. This original proof is quite long and Proposition 11 along with Proposition 12 constitute an alternative to it.

The code is written in the *SSReflect* extension of Coq and uses *mathcomp* libraries.

```coq
Theorem sk_α_is_action {n} : is_action [set: 'Sym_n.+1] (@ska_Residuation n).
Proof.
  rewrite /ska_Residuation/sk_Residuation.
  apply: is_total_action => [C|C p1 p2];
    case: C => C Heq;
    apply sa_inj => /=;
    case: C Heq => n0 p s q t Heq.
    f_equal.
    - rewrite -permP => x.
      by rewrite permE /= cast_permE permE /= cast_ordKV.
    - rewrite /= cast_permE permE /= cast_ordKV eq_refl /=.
      apply eq_from_tnth => i.
      by rewrite tnth_map cast_permE permE /= cast_ordKV tnth_ord_tuple.
    - by rewrite cast_permE permE /= cast_ordKV eq_refl.
    - apply eq_from_tnth => x.
      rewrite tnth_map /=.
      f_equal.
      apply (@perm_inj _ (cast_perm (f_equal S (eqP Heq)) 1)).
      by rewrite cast_permE permE /= cast_ordKV tnth_ord_tuple.
  case H1 : (cast_perm (f_equal succn (eqP Heq)) (p1 * p2) ord_max != ord_max);
    last (move: H1 => /eqP H1);
    (case H2 : (cast_perm (f_equal succn (eqP Heq)) p1 ord_max != ord_max);
      last (move: H2 => /eqP H2);
      (case H3 : (cast_perm (f_equal succn (eqP Heq)) p2 ord_max != ord_max);
        last move : H3 => /eqP H3);
      (f_equal;
        first 1 [
        by rewrite /= -mulgA cast_perm_morphM |
        apply eq_from_tnth => x;
        rewrite !tnth_map !tnth_ord_tuple;
          case H6 : (x != ord_max);
          first (
          case H4 : (cast_perm (f_equal succn (eqP Heq)) (p1 * p2) x != ord_max);
```

```
        last (move: H4 => /eqP H4);
        (case H5 : (cast_perm (f_equal succn (eqP Heq)) p2 x != ord_max);
          last (move: H5 => /eqP H5)));
        last (move: H6 => /eqP H6)];
    last 1 [apply eq_from_tnth => x; rewrite !tnth_map !tnth_ord_tuple;
    by rewrite !cast_permE !permE /= cast_ordK];
    try (rewrite !cast_permE permE /= in H1 H2 H3;
        by rewrite -[p2 _](cast_ordK (f_equal succn (eqP Heq))) H3 H2 eq_refl in H1);
    try (rewrite -[in RHS]H1 in H4;
        by rewrite (perm_inj H4) eq_refl in H6);
    try (rewrite !cast_permE !permE /= in H1 H2 H3;
        rewrite -[p2 _](cast_ordK (f_equal succn (eqP Heq))) H3 in H1;
        by rewrite H1 eq_refl in H2);
    try (rewrite !cast_permE !permE /= in H1 H2 H3;
        move: H1; rewrite -[in RHS]H2 => /cast_ord_inj/perm_inj H1;
        by rewrite H1 cast_ordKV eq_refl in H3);
    try (rewrite !cast_permE !permE /= in H1 H2 H3;
        by rewrite !cast_permE !permE /= !cast_ordK)
  )).
- + rewrite !cast_permE !permE /= in H1 H2 H3 H4 H5.
  repeat (rewrite !cast_permE !permE /= !cast_ordK).
  rewrite !H1 /=.
  move: H1 (H4) => /negbTE H1 /negbTE H4'.
  rewrite !H1 !H4' /=.
  have H7 : (forall p' : 'Sym_n.+1, cast_ord (f_equal succn (eqP Heq))
    (p' (cast_ord (esym (f_equal succn (eqP Heq))) x)) ==
  cast_ord (f_equal succn (eqP Heq))
    (p' (cast_ord (esym (f_equal succn (eqP Heq))) ord_max)) = false).
    intros. apply/negbTE/eqP => /cast_ord_inj/perm_inj/cast_ord_inj/eqP.
    apply/negP. exact: H6.
  rewrite -{2}compM -{1}[p1 (p2 (cast_ord _ x))]compM.
  rewrite !H7.
  have H8 : cast_ord (f_equal succn (eqP Heq))
    (p1 (p2 (cast_ord (esym (f_equal succn (eqP Heq))) ord_max))) ==
  cast_ord (f_equal succn (eqP Heq))
    (p1 (cast_ord (esym (f_equal succn (eqP Heq))) ord_max)) = false.
    apply/negbTE/eqP => /cast_ord_inj/perm_inj/eqP. apply/negP.
    exact: H3.
  have H9 : cast_ord (f_equal succn (eqP Heq))
    (p1 (p2 (cast_ord (esym (f_equal succn (eqP Heq))) x))) ==
  cast_ord (f_equal succn (eqP Heq))
    (p1 (cast_ord (esym (f_equal succn (eqP Heq))) ord_max)) = false.
    apply/negbTE/eqP => /cast_ord_inj/perm_inj/eqP. apply/negP.
    exact: H5.
  rewrite !H8.
```

```
  move: H5 => /negbTE H5.
  rewrite !H5 !cast_ordK /= !H4 !H4' !H9.
  rewrite !addTb !mulgE !addNb !negbK !addbN.
  rewrite addbA -[_(+)(_ _ (_ _ (p1 (p2 _)))(+)_]addbA.
  rewrite [_ _ (_ _ (p1 (p2 _)))(+)_ _ (_ _ (p1 (cast_ord _ _)))]addbC addbA.
  by rewrite addbb addFb.
+ rewrite !cast_permE !permE /= in H1 H2 H3 H4 H5.
  repeat (rewrite !cast_permE !permE /= !cast_ordK).
  move: (H1) (H4) => /negbTE H1' /negbTE H4'.
  have H7 : (cast_ord (f_equal succn (eqP Heq))
    (p1 (p2 (cast_ord (esym (f_equal succn (eqP Heq))) x))) ==
   cast_ord (f_equal succn (eqP Heq))
    (p1 (p2 (cast_ord (esym (f_equal succn (eqP Heq))) ord_max))) = false).
    apply/negbTE/negP => /eqP/cast_ord_inj/perm_inj/perm_inj/cast_ord_inj.
    exact/eqP.
  have H9 : cast_ord (f_equal succn (eqP Heq))
     (p1 (p2 (cast_ord (esym (f_equal succn (eqP Heq))) ord_max))) ==
   cast_ord (f_equal succn (eqP Heq))
     (p1 (cast_ord (esym (f_equal succn (eqP Heq))) ord_max)) = false.
    apply/negbTE/eqP => /cast_ord_inj/perm_inj/eqP. apply/negP.
    exact: H3.
  rewrite H1' H4' H7 H9 /=.
  rewrite -{1}[p2 (_ _ x)](cast_ordK (f_equal succn (eqP Heq))) H5.
  by rewrite !addTb !mulgE !addNb -addbC.
+ rewrite !cast_permE !permE /= in H1 H2 H3 H4 H5.
  repeat (rewrite !cast_permE !permE /= !cast_ordK).
  have H7 : cast_ord (f_equal succn (eqP Heq))
     (p1 (p2 (cast_ord (esym (f_equal succn (eqP Heq))) ord_max))) ==
   cast_ord (f_equal succn (eqP Heq))
     (p1 (cast_ord (esym (f_equal succn (eqP Heq))) ord_max)) = false.
    apply/negbTE/negP => /eqP/cast_ord_inj/perm_inj.
    exact/eqP.
  have H8 : cast_ord (f_equal succn (eqP Heq))
     (p1 (p2 (cast_ord (esym (f_equal succn (eqP Heq))) x))) ==
   cast_ord (f_equal succn (eqP Heq))
     (p1 (cast_ord (esym (f_equal succn (eqP Heq))) ord_max)) = false.
    apply/negbTE/negP => /eqP/cast_ord_inj/perm_inj.
    exact/eqP.
  have H10 : cast_ord (f_equal succn (eqP Heq))
        (p2 (cast_ord (esym (f_equal succn (eqP Heq))) x)) ==
      cast_ord (f_equal succn (eqP Heq))
        (p2 (cast_ord (esym (f_equal succn (eqP Heq))) ord_max)) = false.
    apply/negbTE/negP => /eqP/cast_ord_inj/perm_inj/cast_ord_inj.
    exact/eqP.
  move: (H1) (H5) => /negbTE H1' /negbTE H5'.
```

```
      rewrite !H1 !H1' !H5' !H7 !H10 !cast_ordK !H8 !H4 eq_refl /=.
      by rewrite !addTb !mulgE !addNb negbK addbC addbA addbb addFb.
  + rewrite !cast_permE !permE /= in H1 H2 H3 H4 H5.
      rewrite -[p2 _](cast_ordK (f_equal succn (eqP Heq))) H5 in H4.
      by rewrite H4 eq_refl in H2.
  + rewrite !cast_permE !permE /= in H1 H2 H3.
      repeat (rewrite !cast_permE !permE /= !cast_ordK).
      move: (H1) (H2) (H3) => /negbTE H1' /negbTE H2' /negbTE H3'.
      have H9 : cast_ord (f_equal succn (eqP Heq))
         (p1 (p2 (cast_ord (esym (f_equal succn (eqP Heq))) ord_max))) ==
        cast_ord (f_equal succn (eqP Heq))
         (p1 (cast_ord (esym (f_equal succn (eqP Heq))) ord_max)) = false.
      apply/negbTE/eqP => /cast_ord_inj/perm_inj/eqP. apply/negP.
        exact: H3.
      rewrite !H6 !H1 !H3 !H1' !H3' !H9 !eq_refl /= !H2' /=.
      rewrite !addbA !mulgE !addNb negbK !addbN addbA.
      by rewrite [(_ (+) _) as X in (X (+) _ _ ord_max)]addbC !addbA addbb addFb.
- rewrite tnth_map. rewrite tnth_ord_tuple. rewrite cast_permE permE /= in H1.
  rewrite compM tpermD; last first.
    - apply/eqP => /perm_inj/esym. apply/eqP. exact: H3.
    - rewrite eq_sym !cast_permE /= cast_ordK. exact: H1.
  rewrite !cast_permE permE /= cast_ordK H1.
  rewrite !mulgA !mulgE [in RHS]/= ![in _ (+) _ (_ _ (p1 (cast_ord _ _)))]addbC.
  by rewrite !addbA addbb addFb [_ (+) −]addbC.
- + rewrite !cast_permE !permE /= in H1 H2 H3 H4 H5.
    repeat (rewrite !cast_permE !permE /= !cast_ordK).
    have H9 : cast_ord (f_equal succn (eqP Heq))
       (p1 (p2 (cast_ord (esym (f_equal succn (eqP Heq))) x))) ==
      cast_ord (f_equal succn (eqP Heq))
       (p1 (cast_ord (esym (f_equal succn (eqP Heq))) ord_max)) = false.
      apply/negbTE/eqP => /cast_ord_inj/perm_inj/eqP. apply/negP.
        exact: H5.
    have H7 : (cast_ord (f_equal succn (eqP Heq))
      (p1 (p2 (cast_ord (esym (f_equal succn (eqP Heq))) x))) ==
     cast_ord (f_equal succn (eqP Heq))
      (p1 (p2 (cast_ord (esym (f_equal succn (eqP Heq))) ord_max))) = false).
      apply/negbTE/negP => /eqP/cast_ord_inj/perm_inj/perm_inj/cast_ord_inj.
        exact/eqP.
    move: (H5) (H4) => /negbTE H5' /negbTE H4'.
    rewrite !H4' !H9 !H7 /=.
    by rewrite -[(p2 _) in LHS](cast_ordK (f_equal succn (eqP Heq))) H3.
  + rewrite -[in RHS]H5 in H3.
    apply perm_inj in H3.
    by rewrite H3 eq_refl in H6.
  + rewrite !cast_permE !permE /= in H1 H2 H3 H4 H5.
```

```
    repeat (rewrite !cast_permE !permE /= !cast_ordK).
    rewrite !H4 eq_refl /=.
    by rewrite -[(p2 _) in LHS](cast_ordK (f_equal succn (eqP Heq))) H3.
  + rewrite !cast_permE !permE /= in H1 H2 H3 H4 H5.
    repeat (rewrite !cast_permE !permE /= !cast_ordK).
    rewrite -[(p2 _)](cast_ordK (f_equal succn (eqP Heq))) H5 in H4.
    by rewrite H4 eq_refl in H2.
  + rewrite !H6 !H3 !H1 !H2.
    rewrite !cast_permE !permE /= in H1 H2 H3.
    repeat (rewrite !cast_permE !permE /=).
    move: (H1) (H2) => /negbTE H1' /negbTE H2'.
    rewrite H1' H2' !eq_refl /=.
    by rewrite -[(p2 _) in LHS](cast_ordK (f_equal succn (eqP Heq))) H3.
- rewrite cast_permE permE /=.
  rewrite cast_permE /= in H3.
  rewrite -[p2 _](cast_ordK (f_equal succn (eqP Heq))) H3.
  by rewrite cast_permE.
- + rewrite !cast_permE !permE /= in H1 H2 H3 H4 H5.
    repeat (rewrite !cast_permE /=).
    move: (H4) (H5) => /negbTE H4' /negbTE H5'.
    have H8 : cast_ord (f_equal succn (eqP Heq))
      (p2 (cast_ord (esym (f_equal succn (eqP Heq))) x)) ==
          cast_ord (f_equal succn (eqP Heq))
            (p2 (cast_ord (esym (f_equal succn (eqP Heq))) ord_max)) = false.
      apply/negbTE/eqP => /cast_ord_inj/perm_inj/eqP.
      exact/negP.
    have H9 : cast_ord (f_equal succn (eqP Heq))
      (p1 (p2 (cast_ord (esym (f_equal succn (eqP Heq))) x))) ==
      cast_ord (f_equal succn (eqP Heq))
      (p1 (p2 (cast_ord (esym (f_equal succn (eqP Heq))) ord_max))) = false.
      apply/negbTE/eqP => /cast_ord_inj/perm_inj/perm_inj/eqP.
      exact/negP.
    by rewrite !H4' !H5' !H8 !H9 !cast_ordK.
  + rewrite !cast_permE !permE /= in H1 H2 H3 H4 H5.
    by rewrite -[p2 _](cast_ordK (f_equal succn (eqP Heq))) H5 H2 eq_refl in H4.
  + rewrite !cast_permE !permE /= in H1 H2 H3 H4 H5.
    rewrite -[in RHS]H4 in H2.
    move: H2 => /cast_ord_inj/perm_inj H2.
    by rewrite -H2 cast_ordKV eq_refl in H5.
  + rewrite !cast_permE !permE /= in H1 H2 H3.
    repeat (rewrite !cast_permE !permE /=). rewrite !cast_ordK.
    move: (H1) (H3) => /negbTE H1' /negbTE H3'.
    by rewrite !H6 !H1 !H3 !H1' !H3' !eq_refl !H2.
- f_equal. f_equal.
  by rewrite tnth_map /= tnth_ord_tuple !cast_permE !permE /= cast_ordK.
```

```
- + rewrite !cast_permE !permE /= in H1 H2 H3 H4 H5.
  repeat (rewrite !cast_permE !permE /=). rewrite !cast_ordK.
  move: (H5) (H4) => /negbTE H5' /negbTE H4'.
  have H7 : cast_ord (f_equal succn (eqP Heq))
    (p1 (p2 (cast_ord (esym (f_equal succn (eqP Heq))) x))) ==
    cast_ord (f_equal succn (eqP Heq))
    (p1 (cast_ord (esym (f_equal succn (eqP Heq))) ord_max)) = false.
    apply/negbTE/eqP => /cast_ord_inj/perm_inj/eqP.
    exact/negP.
  have H8 : cast_ord (f_equal succn (eqP Heq))
    (p2 (cast_ord (esym (f_equal succn (eqP Heq))) x)) ==
         cast_ord (f_equal succn (eqP Heq))
           (p2 (cast_ord (esym (f_equal succn (eqP Heq))) ord_max)) = false.
    apply/negbTE/eqP => /cast_ord_inj/perm_inj/eqP.
    exact/negP.
  rewrite !H1 eq_refl /= !H5' !H8 !cast_ordK !H4' /= !H7.
  by rewrite mulgA !mulgE addbb addFb.
+ rewrite !cast_permE !permE /= in H1 H2 H3 H4 H5.
  repeat (rewrite !cast_permE !permE /=). rewrite !cast_ordK.
  by rewrite !H1 eq_refl /= -[p2 _](cast_ordK (f_equal succn (eqP Heq))) H5.
+ rewrite !cast_permE !permE /= in H1 H2 H3.
  repeat (rewrite !cast_permE !permE /=). rewrite !cast_ordK.
  move: (H2) (H3) => /negbTE H2' /negbTE H3'.
  rewrite !H6 !H1 !H3 !H3' !eq_refl !H2 !H2' /=.
  by rewrite mulgA !mulgE addbb addFb.
- rewrite cast_permE permE /= in H1.
  rewrite tnth_map tnth_ord_tuple !cast_permE !permE /= cast_ordK.
  by rewrite H1 eq_refl /= mulgA !mulgE addbb addFb.
Qed.
```



\end{document}